\documentclass[12pt]{article}

\usepackage{amsmath}
\usepackage{graphicx}
\usepackage{enumerate}
\usepackage{natbib}
\usepackage{url} % not crucial - just used below for the URL

\usepackage{amssymb}
\usepackage{amsthm}
\usepackage{color}
\usepackage{caption}
\usepackage{booktabs}
\usepackage{bbm}
\usepackage{multicol}
\usepackage{multirow}
\usepackage{rotating}
\usepackage{pdflscape}

\usepackage[utf8x]{inputenc}
\usepackage{amsfonts,amstext,amssymb, mathtools, comment}
\usepackage[margin=3cm]{geometry}
\usepackage{tikz}
\usepackage[ruled,vlined]{algorithm2e}
\usepackage{rotating}
\usepackage{xcolor}

\newcommand{\p}{{\mathbb P}}

\newcommand{\e}{{\mathbb E}}
\newcommand{\D}{{\mathrm d}}

\renewcommand{\a}{{\alpha}}
\newcommand{\oa}{\overline\alpha}

\newcommand{\br}[1]{\langle #1\rangle}
\newcommand{\norm}[1]{\| #1\|}

\newcommand{\ind}[1]{\mbox{\rm\large  1}_{\{#1\}}}
\newcommand{\R}{\mathbb R}

\newcommand{\convd}{\stackrel{\rm d}{\to}}

\newcommand{\bs}{\boldsymbol}
\newcommand{\X}{\bs X}

\renewcommand{\u}{{\bs u}}
\newcommand{\z}{{\bs z}}
\newcommand{\x}{{\bs x}}
\newcommand{\y}{{\bs y}}
\newcommand{\vv}{{\bs v}}

\newcommand{\V}{V}

\newtheorem{proposition}{Proposition}
\newtheorem{corollary}{Corollary}
\newtheorem{lemma}{Lemma}
\newtheorem{example}{Example}

\begin{document}

\title{\LARGE\bf Distributionally Robust Halfspace Depth}
%% Alternative title{Selecting a depth statistic}
\author{Jevgenijs Ivanovs \\    {\small Department of Mathematics, Aarhus University}\\
    Pavlo Mozharovskyi \\    {\small LTCI, T{\'e}l{\'e}com Paris, Institut Polytechnique de Paris}}
\date{April 30, 2024} 
\maketitle

\vspace{5ex}

\begin{abstract}
Tukey's halfspace depth can be seen as a stochastic program and as such it is not guarded against optimizer's curse, so that a limited training sample may easily result in a poor out-of-sample performance.
We propose a generalized halfspace depth concept relying on the recent advances in distributionally robust optimization, where every halfspace is examined using the respective worst-case distribution in the Wasserstein ball of radius $\delta\geq 0$ centered at the empirical law. This new depth can be seen as a smoothed and regularized classical halfspace depth which is retrieved as $\delta\downarrow 0$. It inherits most of the main properties of the latter and, additionally, enjoys various new attractive features such as continuity and strict positivity beyond the convex hull of the support. We provide numerical illustrations of the new depth and its advantages, and develop some fundamental theory. In particular, we study the upper level sets and the median region including their breakdown properties.
\end{abstract}

%%\noindent{\it AMS 2010 subject classifications:} Primary 62H05, 62H30; secondary 62-07.
%%\indent\\


\section{Introduction}
The earliest and arguably most popular depth concept for multivariate data is the halfspace depth introduced by~\cite{tukey}, see also~\cite{donoho_gasko} and~\cite{MoslerM22} for a survey on data depth for multivariate data mentioning numerous applications.
For a random vector $\X\in\R^d, d\geq 1$ with the law $\p$ it is defined by means of the stochastic program
\[D_0(\z|\p)=\inf_{\u\in\R^d:\|\u\|=1}\p(\br{\u,\X-\z}\geq 0)\]
for all points $\z\in\R^d$, where  $\norm{\u}=\sqrt{\br{u,u}}$ is the Euclidean norm.
In words, one considers the probability mass in the closed halfspace $\{\x\in\R^d:\u^\top(\x-\z)\geq 0\}$ and minimizes over all possible directions~$\u$.
One of the major problems in applications is that the law $\p$ is never known exactly, and so it must be inferred from data.
However, optimization in the model calibrated to a given data set often results in a rather poor out-of-sample performance due to an optimistically-biased estimate.
This \emph{optimizer's curse} is a well-known phenomenon in stochastic programming~\citep{esfahani_kuhn}, and it may roughly be compared to overfitting. In view of its importance for motivation, we elaborate on this question in the below following subsection.

\subsection{Optimizer's curse}

%%For a data set $\datX=\{\x_1,...,\x_n\}$ in $\mathbb{R}^d$ (in a slight abuse of notation, since ties can appear in real data), consider Example~2 from~\cite{KuhnENSA19}. According to ??? and for $d=1$, $\frac{1}{n}\sum_{i=1}^n \ind{\br{\u,\x_i-\z}\ge 0}$ is a consistent estimator of $\p(\br{\u,\X-\z}\geq 0)$, for any $\u$.

%%Lemma and three box-plot diagrams.

For a fixed $\u:\|\u\|=1$, denote ${\overline p}(\u,\z|\p) = \p(\br{\u,\X-\z}\geq 0)$, with its unbiased estimator ${\overline p}(\u,\z|\p_n) = \frac{1}{n}\sum_{i=1}^n\ind{\br{\u,\X_i-\z}\geq 0}$ based on $n$ random vectors $X_1$, ..., $X_n$ stemming from $\p$. This notation yields finite-sample version of the halfspace depth:
\[
	D_0(\z|\p_n)=\inf_{\u\in\R^d:\|\u\|=1} {\overline p}(\u,\z|\p_n)\,.
\]

With $D_0(\z|\p)$ being indeed the (unknown) quantity of interest, for the estimator $D_0(\z|\p_n)$ one can write~\citep[see also][for connection to the machine learning literature]{KuhnENSA19}):
\begin{align*}
	\e( D_0(\z|\p_n) ) &= \e\bigl( \inf_{\u\in\R^d:\|\u\|=1} \frac{1}{n}\sum_{i=1}^n\ind{\br{\u,\X_i-\z}\geq 0} \bigr) \\
	&\le \inf_{\u\in\R^d:\|\u\|=1} \e\bigl( \frac{1}{n}\sum_{i=1}^n\ind{\br{\u,\X_i-\z}\geq 0} \bigr) \\
	&= \inf_{\u\in\R^d:\|\u\|=1} {\overline p}(\u,\z|\p) = D_0(\z|\p)\,.
\end{align*}

\begin{figure}[h!]
\begin{tabular}{ccc}
\quad Normal & \quad Skewed normal & \quad Exponential \\
\includegraphics[width=0.3\textwidth,trim=0 0.5cm 1cm 1.5cm,clip=true]{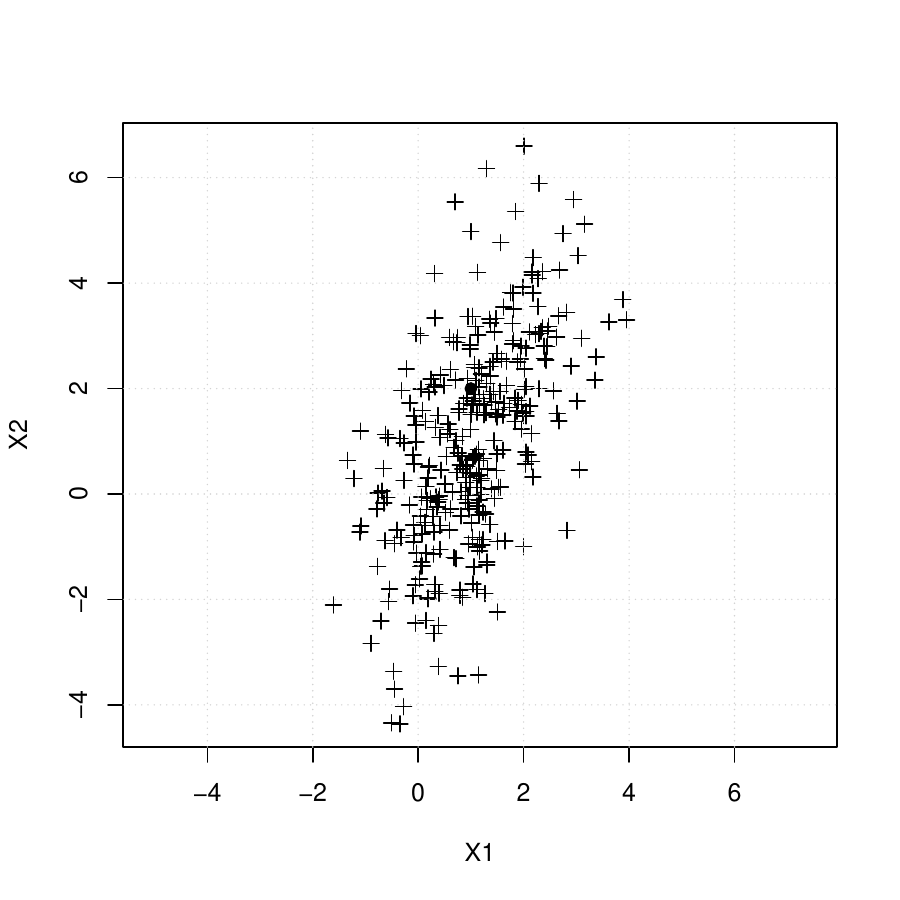} & \includegraphics[width=0.3\textwidth,trim=0 0.5cm 1cm 1.5cm,clip=true]{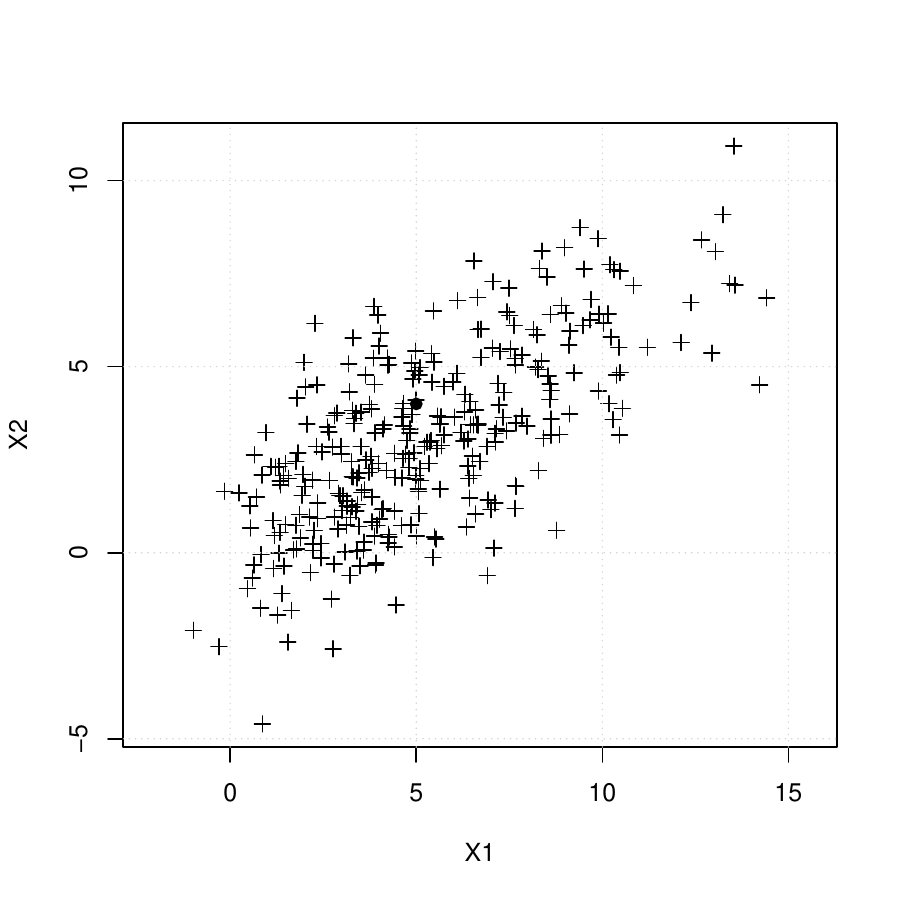} & \includegraphics[width=0.3\textwidth,trim=0 0.5cm 1cm 1.5cm,clip=true]{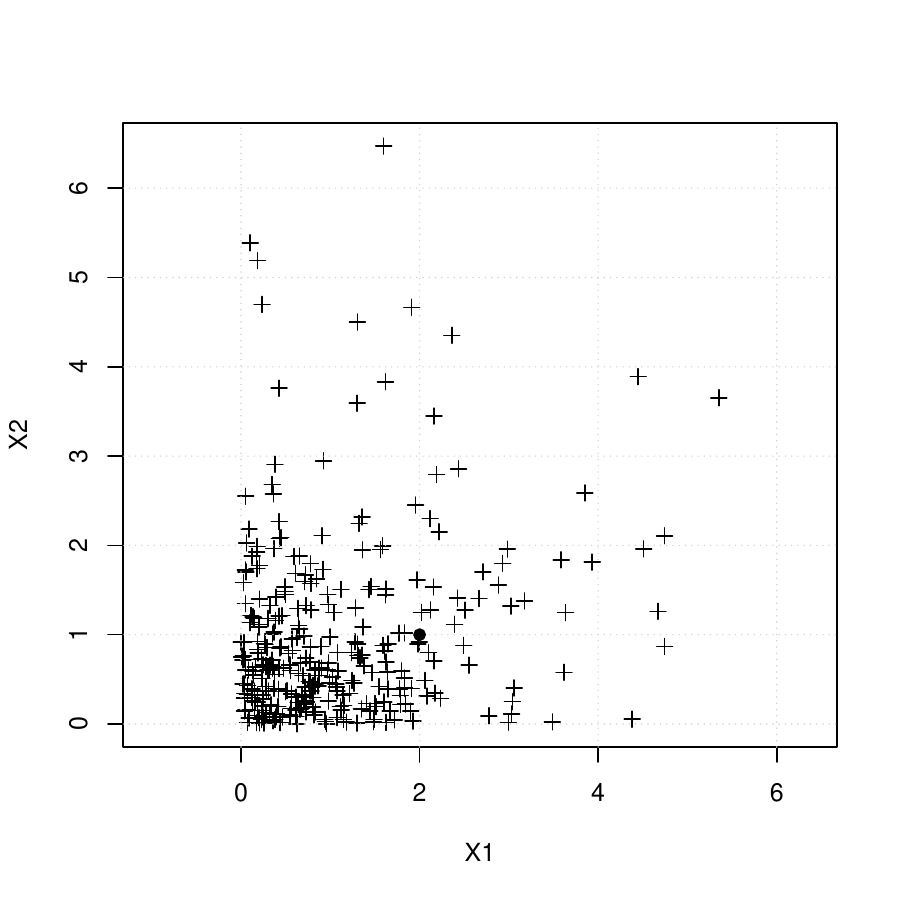} \\
\includegraphics[width=0.3\textwidth,trim=0 0.5cm 1cm 1.5cm,clip=true]{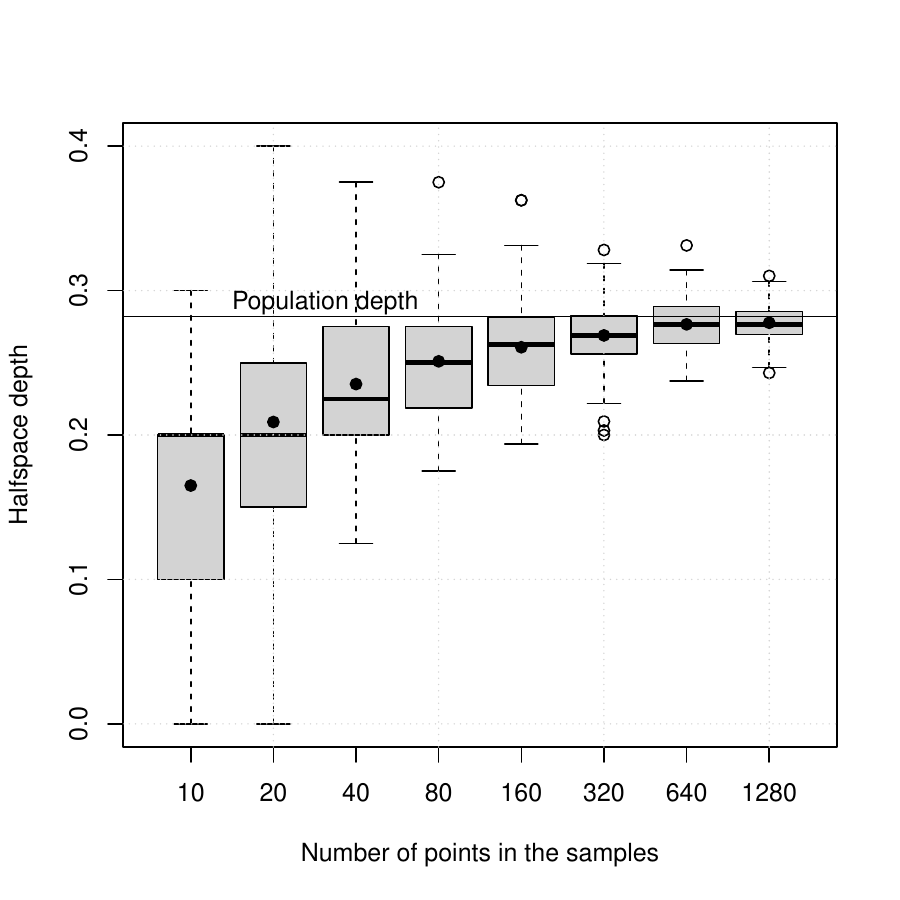} & \includegraphics[width=0.3\textwidth,trim=0 0.5cm 1cm 1.5cm,clip=true]{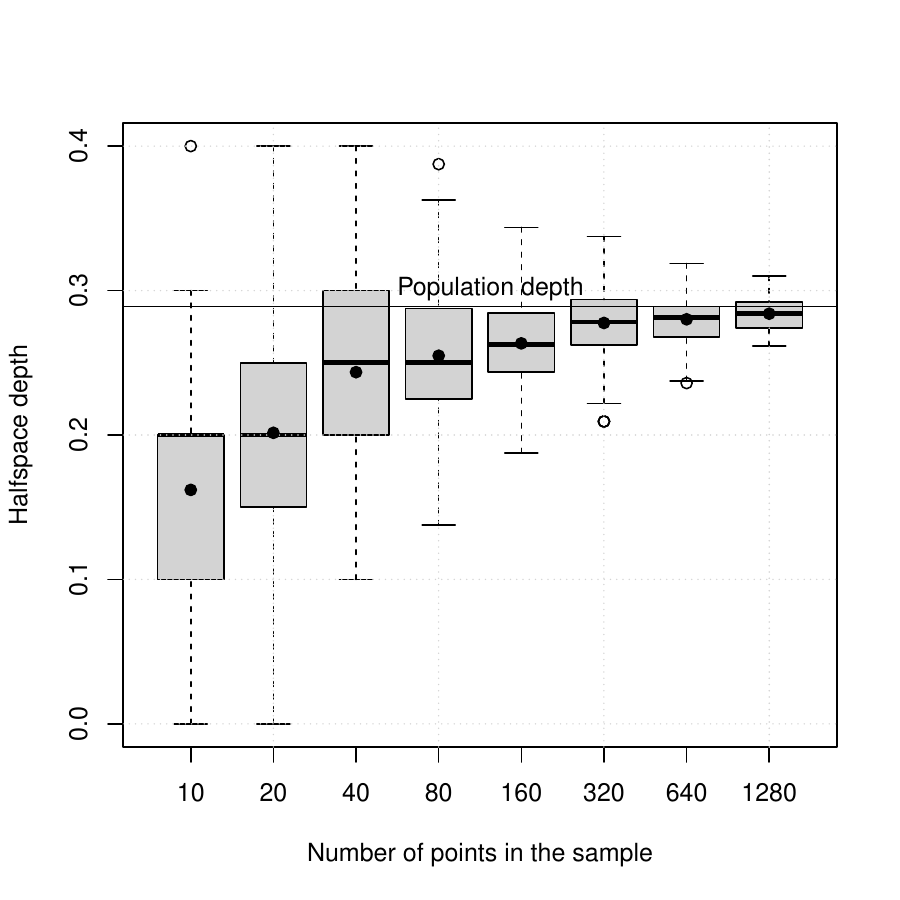} & \includegraphics[width=0.3\textwidth,trim=0 0.5cm 1cm 1.5cm,clip=true]{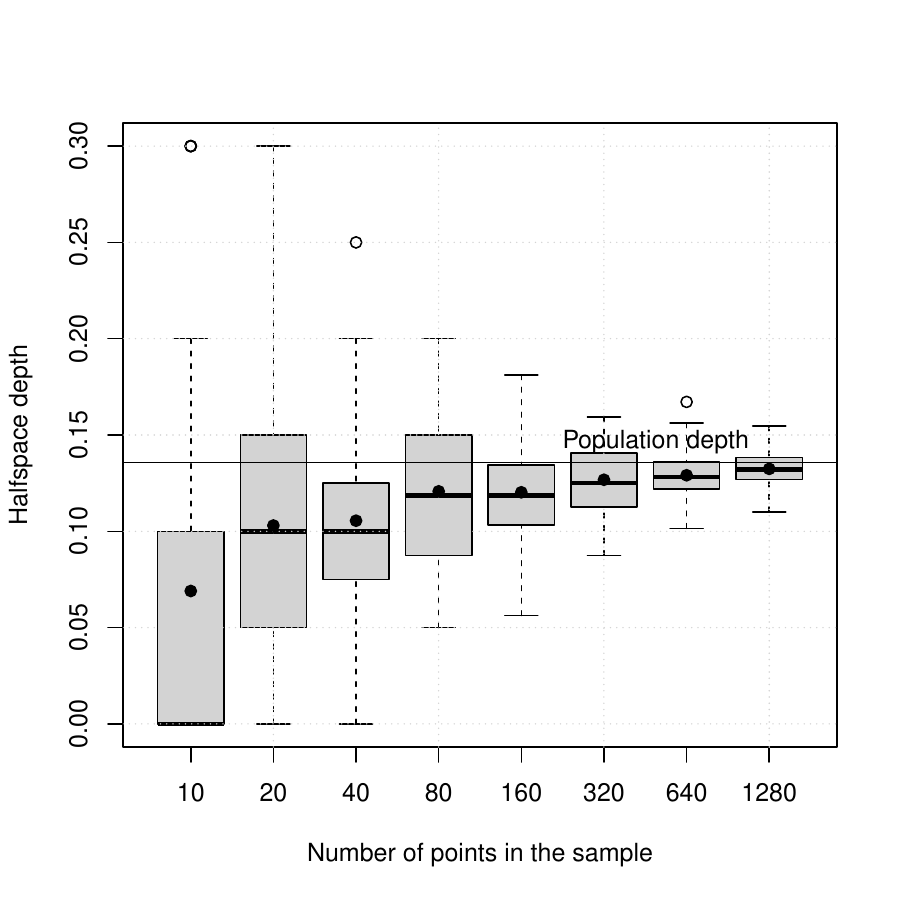}
\end{tabular}
\caption{For normal, skewed normal, and exponential distributions: $320$ bi-variate points (pluses) and $\z$ as point (top row), and the corresponding boxplots of the empirical halfspace depth (over $100$ random samples) with point inside a boxplot indicating mean empirical depth (bottom row).}\label{fig:intro_curse}
\end{figure}

While the inequality holds because the infimum inside the expectation can adapt to the finite sample, it illustrates that the empirical halfspace depth is biased downwards. In Figure~\ref{fig:intro_curse}, we plot empirical halfspace depth for three setting: correlated normal distribution $\mathcal{N}\bigl((1, 1)^\top,\begin{pmatrix} 1 & 1 \\ 1 & 4 \end{pmatrix}\bigr)$, skewed normal distribution with the same covariance matrix and skewness parameter $(5, 1.25)^\top$, and standard exponential distribution with independent marginals. One observes the strong tendency of underestimating the population halfspace depth, with the difference vanishing with growing sample size $n$. It is further noteworthy that, in view of the weak projection property introduced by~\cite{Dyckerhoff04}~\citep[see also][for computational advances based on the latter theoretical result]{CuestaAlbertosNR08,DyckerhoffMN21}, halfspace depth minimization is not necessarily optimal for a finite sample; this probably (additionally) justifies its success in practice.

In this article, we explore at length this phenomenon, and propose a novel notion of data depth.

\subsection{Proposal}

An increasingly popular way of dealing with the optimizer's curse is to explicitly allow for some ambiguity in the approximate model~$\p_n$ and to find the decision optimizing the worst-case expected cost, see~\cite{esfahani_kuhn, scarf1958min}. The worst case is computed over the ambiguity set of laws, which commonly is a ball of probability measures centered at~$\p_n$ with respect to some dissimilarity or distance.
The empirical law~$\p_n$ corresponding to the given training sample is a basic example in this context, and this meaning of $\p_n$ is used in the rest of this paper. %finite number of independent observations of~$X$.
Thus a natural candidate is the Wasserstein (earth mover's) distance between two laws on $\R^d$ (see, for example,~\cite{rachev_ruschendorf}):
\[d_W(\p,\p')=\inf\big\{\int\norm{\x-\y}\Pi(\D \x,\D \y): \Pi\in\mathcal P_{\p,\p'}\big\}\in[0,\infty],\]
where $\mathcal P_{\p,\p'}$ is the set of all joint probability laws on $\R^d\times\R^d$ (equipped with Borel $\sigma$-algebra) with the marginal laws $\p$ and~$\p'$, respectively.
We refer to \cite{ML_earth1,ML_earth2,ML_earth3} for applications of this distance to machine learning tasks. 
%It is noted that we do not exclude laws with infinite first moment, and so we work with an extended distance taking values in $[0,\infty]$.
 % (couplings of the associated random vectors).
%Note that we use the Euclidean distance as the cost function $\norm{X'-X}$ as the cost function associated to transforming $X'$ into~$X$. %, which requires the assumption $\e \norm{X}^2<\infty$.

According to the above described framework we arrive at the following minimax problem for an arbitrary law~$\p$ (an empirical law is just one example)
\begin{equation}\label{eq:minimax}
	D_\delta(\z|\p)=\inf_{\u\in\R^d:\|\u\|=1}\left\{\sup_{\p':d_W(\p',\p)\leq \delta}\p'(\br{\u,\X-\z}\geq 0)\right\}\in[0,1],
\end{equation}
where for every direction $\u$ the worst case probability measure $\p'$ is considered. 
Importantly, the inner problem in~\eqref{eq:minimax} has a simple explicit solution, which follows from the theory of distributional model risk, see~\cite{wozabal2012framework,blanchet_distributional} and references therein. 
We argue that this is a natural way to define a smoothed and regularized version of the halfspace depth.
Normally the focus is on the relative values of the depth function and thus the fact that the robustified halfspace depth in~\eqref{eq:minimax} is at least as large as the traditional depth is of little concern.  

\begin{figure}[h!]
\begin{tabular}{cc}
Traditional Tukey depth $D_0(\cdot|\p_n)$ \qquad & Robust Tukey depth $D_{0.1}(\cdot|\p_n)$ \\
\includegraphics[height=0.42\textwidth,trim=0.1cm 0 0.75cm 1.5cm,clip=true]{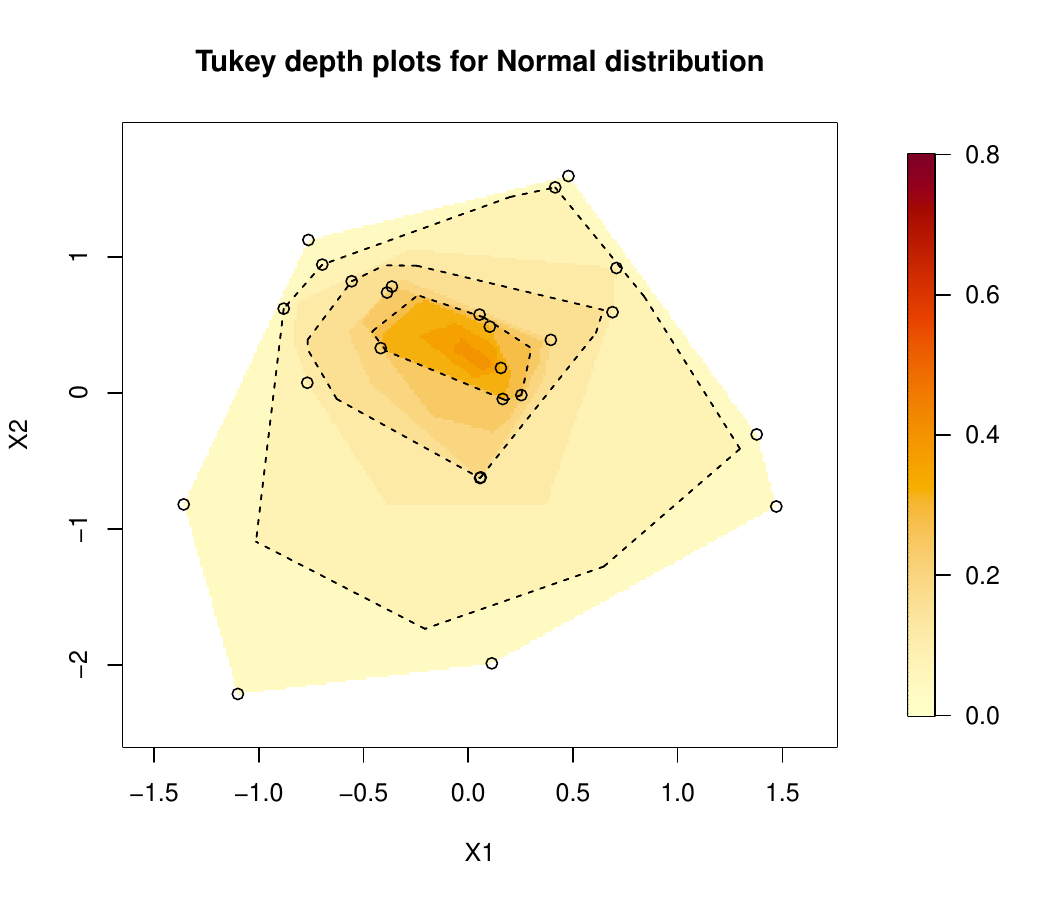} & 
\includegraphics[height=0.42\textwidth,trim=0.1cm 0 0.75cm 1.5cm,clip=true]{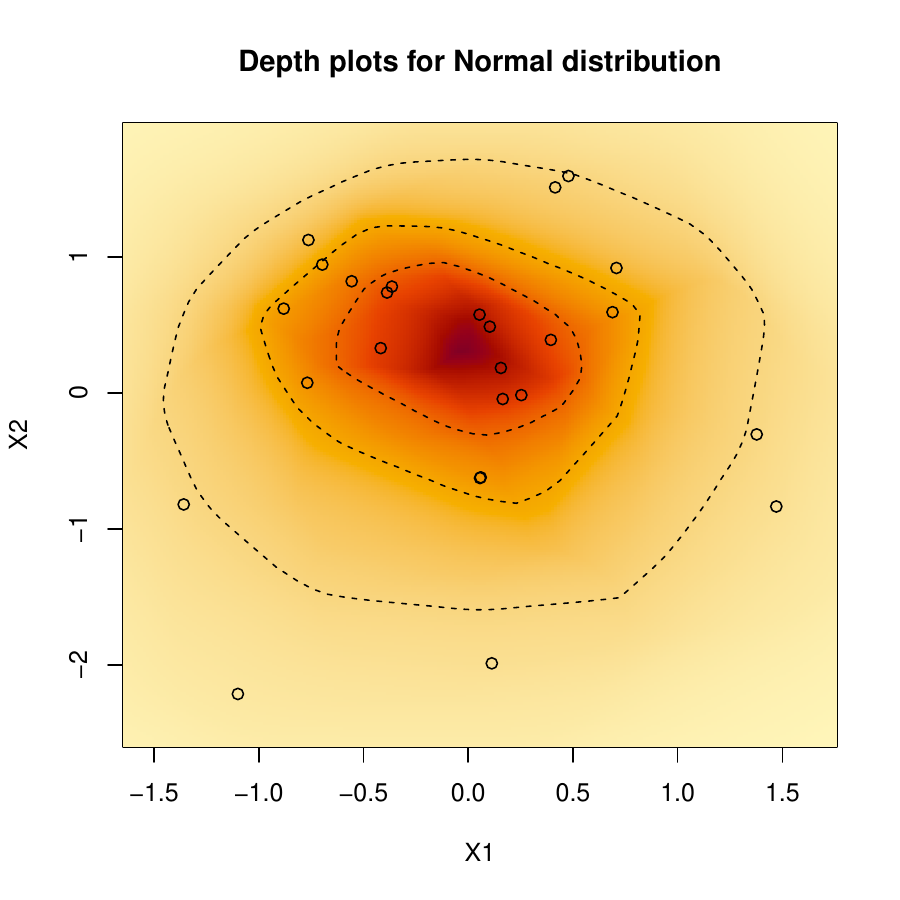}
\end{tabular}
\caption{Heat maps of the traditional Tukey depth $D_0(\cdot|\p_n)$ (left) and the proposed robust Tukey depth $D_{0.1}(\cdot|\p_n)$ (right) for a sample of $25$ points drawn from a bivariate standard  normal distribution. A few contours on the right plot (for depth levels $0.2$, $0.35$, and $0.5$) are presented to emphasize their convexity (later shown formally). With white color being for $0$, the values span from the smallest (shades of yellow) to the highest (shades of red).}
\label{fig:heatmap}
\end{figure}

In the sample version we take $\p_n$ instead of $\p$ and treat $\delta$ as the smoothing parameter, where $\delta\downarrow 0$ results in the empirical halfspace depth, see also Figure~\ref{fig:heatmap} for an illustration.
Unlike the classical halfspace depth the robustified depth is continuous and strictly positive beyond the convex hull of the given observations. %In fact, it is asymptotic to $\delta/\norm{\z}$ as $\norm{z}\to\infty$. 
The level sets $\{\z:D_\delta(\z|\p_n)=\a\}$ are the boundaries of convex sets for all positive $\a$ below the maximal depth.
Moreover, $D_{\delta}(\z|\p_n)$ approaches the true halfspace depth $D_0(\z|\p)$ as $n\to\infty$ and $\delta\downarrow 0$.
We note that in some other contexts such minimax problems can be indeed rewritten as regularized original problems, see~\cite{blanchet2016profile,bousquet2002stability,ghaoui_regularization,SVM}. 
%\pavlo{Maybe detail more the algorithm.}
%In the following we also show that our minimax problem is equivalent to the respective maximin problem, so that first the Tukey depth is found for all $\p'$ from the ambiguity ball and then the maximal depth is returned.
%The worst case probability law amounts to shifting certain mass closest to the given halfspace right to its 

\begin{figure}[h!]
\includegraphics[width=0.31\textwidth,trim=0 0 0 1.5cm,clip=true,page=1]{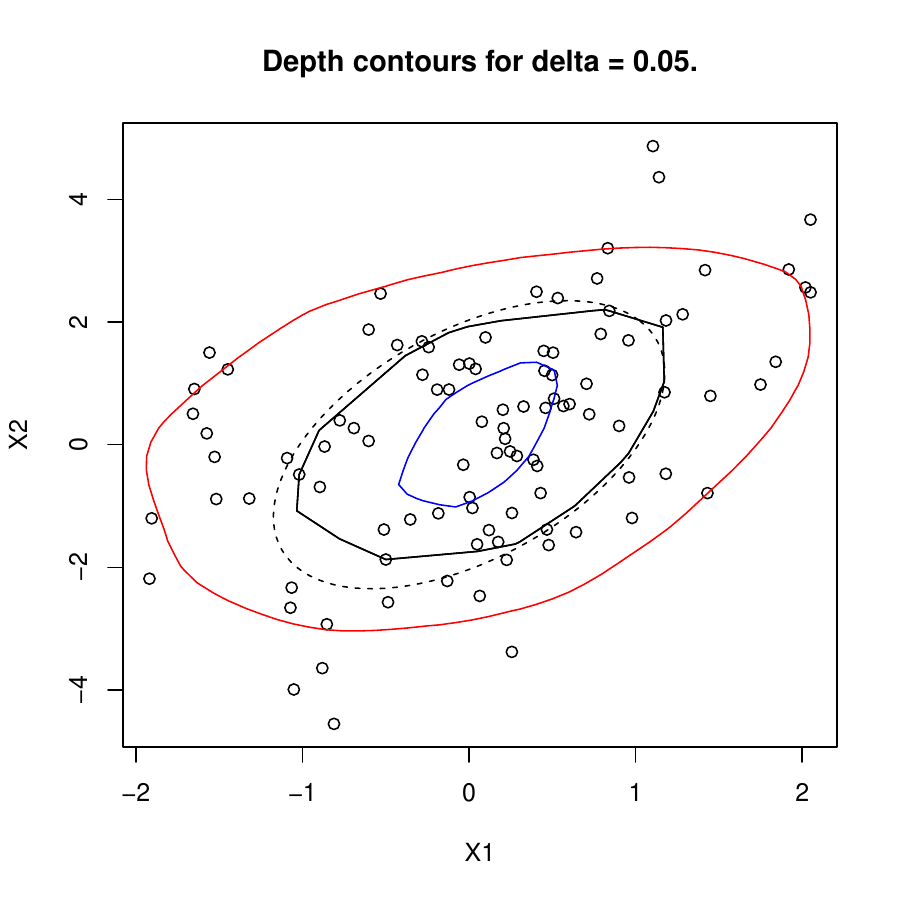}\quad
\includegraphics[width=0.31\textwidth,trim=0 0 0 1.5cm,clip=true,page=2]{picContours-n100-Normal.pdf}\quad
\includegraphics[width=0.31\textwidth,trim=0 0 0 1.5cm,clip=true,page=3]{picContours-n100-Normal.pdf}
\caption{Depth contours (upper-level sets of the depth function) to the level $\alpha=0.12$ for a sample of $100$ points drawn from a bivariate normal distribution centered at $\bs 0$ and with the covariance matrix $\bigl((1, 1)^\top, (1, 4)^\top\bigr)$ for the following depth functions: $D_0(\cdot|\p)$ (black dashed), $D_0(\cdot|\p_n)$ (black solid), $D_\delta(\cdot|\p_n)$ (red), and $\underline D_\delta(\cdot|\p_n)$ (blue). We use $\delta=0.05$ (left), $\delta=0.015$ (middle), and $\delta=0.005$ (right).}
\label{fig:bounds}
\end{figure}

In addition, we also briefly consider
\[\underline D_\delta(\z|\p)=\inf_{\u\in\R^d:\|\u\|=1}\left\{\inf_{\p':d_W(\p',\p)\leq \delta}\p'(\br{\u,\X-\z}> 0)\right\},\]
where the inner problem minimizes the mass in the open halfspace over the Wassertsein ball of measures around~$\p$.
Assuming $d_W(\p,\p_n)\leq \delta$ for some $\delta>0$ we get simple bounds on the true halfspace depth:
\[\underline D_\delta(\z|\p_n)\leq D_0(\z|\p)\leq D_\delta(\z|\p_n),\]
see Figure~\ref{fig:bounds}. This can be used to construct asymptotic confidence intervals and finite sample guarantees~\citep{esfahani_kuhn}.
Finding the right threshold $\delta>0$, however, is problematic. One can use modern measure concentration results~\citep{fournier2015rate}, but thus obtained bounds are normally overly conservative. We do not pursue the idea of sandwiching the true halfspace depth in the following, but rather focus on the properties of the robustified depth.

\subsection{Outline}
The main focus of this paper is on the \emph{robustified halfspace depth} defined in~\eqref{eq:minimax}.
The inner problem is studied in Section~\ref{sec:sol}, where we also establish its various useful properties and provide Algorithm~\ref{alg:sample} for the sample version. 
In Section~\ref{sec:props} we prove the above mentioned properties of the robustified halfspace depth for a general law~$\p$, establish its decay at infinity, and show that the corresponding maximin problem does not necessarily coincide with  our minimax  formulation in~\eqref{eq:minimax}. %, and also discuss its relation to other depths known in the literature~\citep{chernozhukov,einmahl2015,hlubinka,nagy_illumination}.
%%%Finite sample version of the robustified depth is further studied in Section~\ref{sec:sample}, where we prove asymptotic consistency and show that 
%%%the breakdown point of an upper level set is similar to the case of the classical halfspace depth which, nevertheless, results in a more robust model.
Finite sample version of the robustified depth is further studied in Section~\ref{sec:sample}, where we prove asymptotic consistency, suggest an affine-invariant version, and show that the breakdown point of an upper level set can be higher than that of the of the classical halfspace depth.
Furthermore, the median region may have the asymptotic breakdown point at $1/2$.
Section~\ref{sec:parchoice} explores the choice of the parameter $\delta$.
Simple representations of the outer level sets and the median region are presented in Section~\ref{sec:representations} including Algorithm~\ref{alg:median} for the construction of the latter.
Section~\ref{sec:relations} analyses relation of our proposal to other existing relevant depth notions, while 
some advantages of using the robustified depth are exemplified numerically in Section~\ref{sec:numerics}. All proofs are deterred to the Supplementary Materials.

\section{The inner problem}\label{sec:sol}
This section is devoted to showing that the inner problem in~\eqref{eq:minimax} has a simple explicit form, which readily follows from the general theory in~\cite{blanchet_distributional}.
Furthermore, we establish various useful properties of the solution and specialize to the case where the center of the ambiguity ball is given by the empirical law. It must be mentioned that the latter case can also be treated using the results from~\cite{esfahani_kuhn}.
\subsection{The general solution and the optimal transport plan}
For any $\z\in\R^d$ and any direction $\u\in\R^d,\norm{\u}=1$ we consider a random variable $Y\in\R$, which is the projection of $\z-\X$ onto the direction $\u$:
\[Y=\br{\u,\z-\X},\qquad\text{so that}\qquad\{\br{\u,\X-\z}\geq 0\}=\{Y\leq 0\}.\]
Importantly, the random variable $Y$ will be used to express the solution to the inner problem in~\eqref{eq:minimax}.
It will be shown that for a given direction $\u$ this $d$-dimensional optimization problem reduces to a one-dimensional problem, where the ambiguity is specified for the random variable~$Y$.
In this regard, we write $\p_Y$ for the law of $Y$ on $\R$ and $d_W^1$ for the Wasserstein distance between the laws on the real line with the cost function $c(x,y)=|x-y|$.

Consider the truncated expectation of $Y$ and its left-inverse:
\begin{equation}\label{eq:h}h(y)=\e (Y\ind{Y\in(0,y]}),\qquad h^{-1}(\delta)=\inf\{y\geq 0:h(y)\geq \delta\},\end{equation} %,\qquad \underline h(y)=-\e (Y\ind{Y\in[-y,0)})
where $y,\delta\geq 0$. Note that $h$ is a non-decreasing right-continuous function, and $h^{-1}$ is a non-decreasing left-continuous function with values in~$[0,\infty]$.
Note that $h(0+)=0$ and thus $h^{-1}(\delta)>0$ for $\delta>0$.
 %, and define $\underline h^{-1}(\delta)$ analogously.
%We will restrict to $\delta>0$ and then $h^{-1}(\delta)>0$. 
Finally, $a^+=\max(a,0)$ and $a^-=\max(-a,0)$ denote the positive and the negative parts of $a$, respectively.

%The following is essentially the result of~\cite{blanchet_distributional}. Robust bound on an non-exceedance probability can be written as the original non-exceedance probability for a certain larger threshold~$\lambda$.

\begin{proposition}[Inner problem]\label{prop:sup}
For $\delta\in(0,\infty)$ and $\z,\u\in\R^d$ with $\norm{\u}=1$ it holds that
\begin{align}
&\sup_{\p':d_W(\p',\p)\leq \delta}\p'(\br{\u,\X-\z}\geq 0)& &=\sup_{\p':d_W^1(\p',\p_Y)\leq \delta}\p'(Y\leq 0)\label{eq:opt}\\
&\qquad=\inf_{\lambda'>0}\{\delta/\lambda'+\e(1-Y^+/\lambda')^+\}& &=\delta/\lambda+\e(1-Y^+/\lambda)^+\label{eq:sol}\\
&\qquad=\p(Y\leq\lambda)-(h(\lambda)-\delta)/\lambda& &=\p(Y<\lambda)+(\delta-h(\lambda-))/\lambda,\label{eq:sol_alt}
\end{align}
where $\lambda=h^{-1}(\delta)\in(0,\infty]$ and~\eqref{eq:sol_alt} is not used when $\lambda=\infty$.
Furthermore,~\eqref{eq:sol} equals 1 iff $\e Y^+\leq \delta$. 

The minimization problem has the solution:
\begin{align}
&\inf_{\p':d_W(\p',\p)\leq \delta}\p'(\br{\u,\X-\z}> 0)& &=\inf_{\p':d_W^1(\p',\p_Y)\leq \delta}\p'(Y<0) \nonumber \\
& & &=\p(Y<-\underline\lambda)+(\underline h(\underline \lambda)-\delta)/\underline\lambda, \label{eq:sol2}
\end{align}
where $\underline\lambda=\underline h^{-1}(\delta)\in(0,\infty]$ and $\underline h(y)=\e (-Y\ind{-Y\in(0,y]})$ is the analogue of $h$ for the $-Y$ variable.
Furthermore, \eqref{eq:sol2} equals 0 iff $\e Y^-\leq \delta$.
\end{proposition}
%%%\begin{proof}
%%%The result in~\eqref{eq:sol} follows readily from~\cite[Thm.\ 1 and Thm.\ 2(a)]{blanchet_distributional} by noting that the distance from $\x\in\R^d$ to $A=\{\x:\br{\u,\x-\z}\geq 0\}$ is given by $y^+$ with $y=\br{\u,\z-\x}$.
%%%This also shows the equivalence of the two optimization problems in~\eqref{eq:opt}, which can also be proven directly by relating the ambiguity sets.
%%%The representation in~\eqref{eq:sol_alt} for $\lambda<\infty$ follows by rewriting the expectation of the positive part.
%%%Suppose that the value is 1 then it can not be that $\e Y^+>\delta$, because  then $\lambda<\infty$ and $\p(Y\leq \lambda)<1$ whereas $h(\lambda)\geq \delta$. 
%%%Assuming $\e Y^+\leq \delta$ we consider the two cases $\lambda=\infty$ and $\lambda>\infty$ to find that the value is~1.

%%%The minimization problem is analyzed by considering
%%%\[1-\sup_{\p':d_W(\p',\p)\leq \delta}\p'(\br{\u,\X-\z}\leq 0),\]
%%%and noting that the subset $A=\{\x:\br{\u,\x-\z}\leq 0\}$ is closed and the distance of $\x$ to $A$ is given by~$y^-$. 
%%%We conclude by noting that $1-(1-a)^+=\min(a,1)$ and simplifying the result $\e\min(Y^-/\underline\lambda,1)-\delta/\underline\lambda$.
%%%\end{proof}

With regard to the representation in~\eqref{eq:sol_alt} we note that for $\lambda<\infty$:
\[h(\lambda)\geq \delta\geq h(\lambda-)\]
and the solution reads simply $\p(Y\leq \lambda)$ when $Y$ has no mass at~$\lambda$. The latter is always true if $Y$ has no atoms. 
In this case we may view the solution as the original probability, where the hyperplane defining the halfspace is shifted appropriately in the direction opposite to~$\u$. 

Importantly, the optimal transport plans exist for both~\eqref{eq:opt} and~\eqref{eq:sol2}, see~\citet[Lem.\ 3]{blanchet_distributional}.
That is, there is a random variable $\X^*$ (not unique in general) on a possibly extended probability space $(\Omega,\mathcal F,\p)$ such that
\begin{equation}\label{eq:opt_transport}\sup_{\p':d_W(\p',\p)\leq \delta}\p'(\br{\u,\X-\z}\geq 0)=\p(\br{\u,\X^*-\z}\geq 0),\qquad \e\norm{\X^*-\X}\leq \delta.\end{equation}
Furthermore, the structure of $\X^*$ is very simple: it coincides with $\X$ when $\br{\u,\X-\z}=-Y\geq 0$ and it is the projection $\X+\u\br{\u,\z-\X}$ onto the hyperplane when $Y\in(0,\lambda)$ or when $Y=\lambda$ and $U\leq (\delta-h(\lambda-))/\lambda$ with some independent uniform~$U$. The latter uniform variable is needed to transport only some part of the mass corresponding to the atom of $Y$ at~$\lambda$, see also Figure~\ref{fig:illustration} below.
When the maximal value is~1 the above should be understood as transporting all the mass corresponding to $Y>0$ to the respective hyperplane, and this is the only case when $\e\norm{\X^*-\X}$ can be strictly below $\delta$; see Section~. 
For the infimum all the mass corresponding to $Y\in(-\underline\lambda,0)$ is moved to the hyperplane and additionally some mass corresponding to $Y=-\underline\lambda$ as well.
%If $h^{-1}(\delta)=\infty$ or $\underline h^{-1}(\delta)=\infty$ then all the mass corresponds to moving all the mass from $\{Y>0\}$ ($\{\}$) to~0.
Optimal transport plans of this kind are standard in various related problems, see~\cite{wozabal2014robustifying,Pflug_review} and references therein.

\subsection{The sample version}
Here we take $\p=\p_n=\frac{1}{n}\sum_{i\leq n}\delta_{\{\x_i\}}$ to be the empirical measure corresponding to the observations $\x_i\in\R^d$, and further simplify the expressions in Proposition~\ref{prop:sup}.
Let $y_i=\br{\u,\z-\x_i},i=1,\ldots,n$ be the projected observations, and denote by 
\[y^{(-\underline m)}\leq \cdots\leq y^{(-1)}\leq 0<y^{(1)}\leq y^{(2)}\leq \cdots\leq y^{(m)},\qquad s_i=\sum_{j=1}^i y^{(j)},\qquad \underline s_i=-\sum_{j=1}^i y^{(-j)}\] the sorted values of~$y_i$ and the partial sums of positive and non-positive projections, respectively, where $m,\underline m\geq 0$ and $\underline m+m= n$. For convenience, we let that $s_0=\underline s_0=0$, $s_i=s_{m}$ for $i>m$, and similarly $\underline s_i=\underline s_{\underline m}$ for $i>\underline m$.
In particular, we have 
\[p=\p_n(\br{\u,\X-\z}\geq 0)=(n-m)/n=\underline m/n.\]
It must be noted that the following result can also be obtained from~\citet[Cor.\ 5.3]{esfahani_kuhn}. 
\begin{corollary}[Sample version]\label{cor:sample} Assume that $\norm{\u}=1$.\\
If $1\leq k\leq m$ is such that $s_{k-1}<\delta n\leq s_{k}$ then 
\begin{equation}\label{eq:sol_emp}\sup_{\p':d_W(\p',\p_n)\leq \delta}\p'(\br{\u,\X-\z}\geq 0)=p+\frac{k-1}{n}+\frac{\delta-s_{k-1}/n}{y^{(k)}}\in\Big(p+\frac{k-1}{n},p+\frac{k}{n}\Big].\end{equation}
Furthermore, $s_m\leq  \delta n$ yields~1.

If $1\leq k\leq \underline m$ is such that $\underline s_{k-1}<\delta n\leq \underline s_{k}$ then 
\[\inf_{\p':d_W(\p',\p_n)\leq \delta}\p'(\br{\u,\X-\z}> 0)=p-\frac{k-1}{n}+\frac{\delta-\underline s_{k-1}/n}{y^{(-k)}}\in\Big[p-\frac{k}{n},p-\frac{k-1}{n}\Big).\]
Furthermore, $\underline s_{\underline m}\leq \delta n$ yields~0.
\end{corollary}
%%%\begin{proof}
%%%Note that $s_m< \delta n$ corresponds to $\hat\e Y^+<\delta$ and thus we get~1.
%%%Otherwise, choose $1\leq k\leq  m$ as stated. Then $\lambda=y^{(k)}$ and according to Proposition~\ref{prop:sup} the supremum is given by
%%%\[p+\frac{1}{n}\sum_{i=1}^{k-1}(1-y^{(i)}/y^{(k)})+\delta/y^{(k)},\]
%%%because the indices corresponding to $y^{(i)}=y^{(k)}$ are irrelevant. The first result now follows, where the upper bound stems from the inequality $\delta n-s_{k-1}\leq s_k-s_{k-1}=y^{(k)}$. 
%%%
%%%Similar argument for the second expression gives $\underline \lambda=-y^{(-k)}$ and then also
%%%\[\frac{\underline m-k+1}{n}-\frac{1}{n}\underline s_{k-1}/y^{(-k)}+\delta/y^{(-k)}\]
%%%and the result follows.
%%%\end{proof}
For better clarity and visibility  we provide the algorithm corresponding to Corollary~\ref{cor:sample}. It's complexity derives from the complexity of sorting, and so it is $O(n\log n)$.

\begin{algorithm}[H]
%\SetAlgoLined
 compute $y_i=\br{\u,\z-\x_i}$\;
 pick $y_i>0$ and sort them in increasing order\;
 let $m\in \{0,\ldots,n\}$ be their number\;
 compute cumulative sums $s_i=\sum_{j=1}^i y_j$ with $s_0=0$\;
 \If{$s_m\leq \delta n$}{  \KwResult{$1$}}
 find  the smallest $k\in\{1,\ldots,m\}$ such that $s_k\geq \delta n$\;
 \KwResult{$\frac{n-m}{n}+\frac{k-1}{n}+\frac{\delta-s_{k-1}/n}{y_k}$}

 \caption{Compute $\sup_{\p':d_W(\p',\p_n)\leq \delta}\p'(\br{\u,\X-\z}\geq 0)$.}\label{alg:sample}
\end{algorithm}
\medskip

%Even if $p_n=0$ and $k=1$ then the supremum is $\delta/y^{(1)}$, which is inversely proportional to the distance to the closest point on the other side of the hyperplane.
%In particular, the robust bound is always positive. 
Finally, we note that the optimal transport plan for the supremum consists of moving the points $\x_j$ corresponding to $y^{(i)}$,$i=1,\ldots, k-1$ onto the given hyperplane by means of projection.
Additionally, certain mass $\leq 1/n$ is shifted from the next closest point to the hyperplane so that $\sum(\text{mass}_i\times\text{distance}_i) =\delta$ unless $s_m/n<\delta$, see Figure~\ref{fig:illustration}.

\begin{figure}
\centering
\begin{tikzpicture}
\fill[gray!10] (4,0)--(4,4)--(5,4)--(5,0);
\filldraw (4.6949597, 1.2503322) circle(1pt);
\filldraw (3.7649484, 0.7963902) circle(1pt);
\filldraw (3.4675316, 3.5189625) circle(1pt);
\filldraw (3.1438376, 2.6455422) circle(1pt);
\filldraw (1.5477594, 3.8998995) circle(1pt);
\filldraw (0.5773439, 1.4530794) circle(1pt);
\filldraw (2.7536316 ,3.0210816) circle(1pt);
\filldraw (0.5980348, 2.6250300) circle(1pt);
\filldraw (2.7720839, 1.0943229) circle(1pt);
\filldraw (2.8840230, 2.2034970) circle(1pt);

\draw[red, ->] (4,2) -- (0,2);
\filldraw (4,2) circle(2pt) node[below right]{$\z$};
%\draw[->] (4,2)--(5,2);
%\draw (4.5,2) node[above]{$u$};
\draw (4,0) -- (4,4);
\draw[red,dotted] (3.7649484, 0.7963902) -- (3.7649484, 2);% node[red,below]{$y^{(1)}$};
\draw[red, dotted] (3.4675316, 3.5189625) -- (3.4675316, 2);% node[red,below]{$y^{(2)}$};
\draw[red, dotted] (3.1438376, 2.6455422) -- (3.1438376, 2);
\draw[red] (3.1438376, 2.1) -- (3.1438376, 1.9) node[red,below]{$y^{(3)}$};
\draw[blue] (3.7649484, 0.7963902) circle(2pt);
\draw[blue] (3.4675316, 3.5189625) circle(2pt);
\draw[blue,dashed] (3.1438376, 2.6455422) circle(2pt);
\draw[blue,->] (3.7649484, 0.7963902) -- (4, 0.7963902);
\draw[blue, ->] (3.4675316, 3.5189625) -- (4, 3.5189625);
\draw[blue, dashed,->] (3.1438376, 2.6455422) -- (4, 2.6455422);
\end{tikzpicture}
\caption{Illustration of the solution in~\eqref{eq:sol_emp}: $p=1/n$ and $k=3$ with the result in $(3/n,4/n]$. The optimal transport plan~\eqref{eq:opt_transport} is depicted by blue arrows, where the dashed arrow indicates that less than $1/n$ mass is moved.}
\label{fig:illustration}
\end{figure}

\subsection{Properties of the solution}
This section is devoted to various properties of the solution to the inner problem, which form the basis for the following study of the robustified depth $D_\delta(\z|\p)$.
For a fixed $\delta\geq 0$  and some random variable $Y_i\in\R$ we let 
\[v_i=v_i(\delta)\in[0,1]\text{ be the optimal value in~\eqref{eq:opt} with }Y\stackrel{d}{=}Y_i.\]
In the case $\delta=0$ we take $v_i=\p(Y_i\leq 0)$. 
Note that $v_i$ is a non-decreasing function of $\delta$.

For further use we first consider a deterministic $Y=c$. Proposition~\ref{prop:sup} readily implies that 
\begin{equation}\label{eq:det}
v=\begin{cases}\delta/c,&\text{if }c>\delta,\\1,&\text{otherwise.}\end{cases}
\end{equation}
Next, recall that $Y_1$ is stochastically smaller than $Y_2$ ($Y_1\prec Y_2$) if $\p(Y_1 \leq y)\geq \p(Y_2 \leq y)$ for all $y\in \R$. 
We start with a monotonicity result.
%%The proofs of all of these results are postponed to Appendix~\ref{sec:proofs_props}.
\begin{lemma}[Monotonicity]\label{lem:mon}
If $Y^+_1\prec Y^+_2$ then $v_1\geq v_2$ for any~$\delta\geq 0$.
\end{lemma}
Note that $Y_1\prec Y_2$ implies $Y_1^+\prec Y_2^+$.%, but we still have the equality $v_2=v_1$. 
%This follows immediately from the optimal transport plan or the explicit solution in Proposition~\ref{prop:sup}.

\begin{lemma}[Mixture]\label{lem:mix}
Let $Y$ be a mixture of $Y_1$ and $Y_2$ with probabilities $p\in(0,1)$ and $1-p$, respectively.
Then for $\delta\geq 0$ we have
\[v=\sup\{pv_1(\delta c_1)+(1-p)v_2(\delta c_2):c_1,c_2\geq 0, pc_1+(1-p)c_2\leq 1\}.\] 

In particular, there is the sandwich bound:
\[pv_1(\delta/p)+(1-p)v_2(0)\leq v\leq pv_1(\delta/p)+(1-p)v_2(\delta/(1-p)).\]
\end{lemma}
These bounds, for example, imply the following property which will be used later:
\begin{equation}\label{eq:inc_v}
pv(\delta/p) \quad\text{is non-decreasing in }p\in[0,1],
\end{equation}
where $p=0$ yields 0.
Take a trivial mixture with $Y_i\stackrel{d}{=}Y$ to find that $v(\delta)\geq pv(\delta/p)$.
Now the inequality for $0< p_1<p_2\leq 1$ follows from $v(\delta/p_2)\geq (p_1/p_2)v((\delta/p_2)/(p_1/p_2))$.
Furthermore, we observe that 
\begin{equation}\label{eq:lowerbound}
\p(Y\leq c)\geq p\text{ for some }c>0\qquad \text{ then }\qquad v(\delta)\geq \min(p,\delta/c).
\end{equation}
This follows by regarding $Y$ as a mixture of $Y|Y\leq c$ and $Y|Y>c$, using the lower bound and monotonicity in~\eqref{eq:inc_v}, and then applying Lemma~\ref{lem:mon} and~\eqref{eq:det}.

\begin{lemma}[Boundary values]\label{lem:infinite}
The following is true for $\delta\geq 0$:
\begin{itemize}
\item[(i)] If $\p(Y_t>c)\to 1$ for every $c>0$ then $v_t\to 0$,
\item[(ii)] If $\p(Y_t\leq 0)\to 1$ then $v_t\to 1$.
\end{itemize}

Let $Y$ be a mixture of $Y_1$ and $Y_t$ with probabilities $p\in(0,1)$ and $1-p$, respectively. Then 
$v\to pv_1(\delta/p)$ in case (i) and $v\to pv_1(\delta/p)+(1-p)$ in case (ii).
\end{lemma}

The next result is crucial for the property P6 of the depth below.
\begin{lemma}[Strict monotonicity]\label{lem:strict}
Let $Y_2=Y_1+\eta$ for $\eta>0$. Then $v_1=v_2$ for $\delta>0$ implies that both $v_i$ are~1.
\end{lemma}
We conclude with a continuity result.
\begin{lemma}[Continuity]\label{lem:cont}
Consider $\delta_n\to \delta>0$ and $Y_n\convd Y$. Then $v_n(\delta_n)\to v(\delta)$.
\end{lemma}
Note that the boundary case $\delta=0$ is excluded from the above result, but see Proposition~\ref{prop:approx} below.

\section{Properties of the robustified halfspace depth}\label{sec:props}
In this section we show that the robustified halfspace depth defined in~\eqref{eq:minimax} satisfies the desirable properties of a depth function~\citep{zuo_serfling} with a relaxation of affine invariance, and that it has various further nice features. In addition, we study its decay at infinity and show that the corresponding maximin problem may give a strictly smaller value.
\subsection{The basic properties}
Throughout this section we assume that the law $\p$ of the random vector $\X\in\R^d$ is fixed.
%function $D_\delta(\z)=D_\delta(\z|\p)$ defined in~\eqref{eq:minimax} for a fixed threshold $\delta>0$. 
We use $\p_{\X'}$ to denote the law of the random variable $\X'\in \R^d$.
The first result shows that the robustified halfspace depth $D_\delta$ inherits all the main properties of the classical halfspace depth $D_0$, whereas the second result establishes various additional nice properties.
We show that our depth function is continuous and strictly positive with level sets being the boundaries of the convex upper level sets. This should be compared to the `staircase' form of the empirical halfspace depth.
Finally, we show that $D_\delta$ decreases to the halfspace depth as $\delta\downarrow 0$.
Thus we regard $D_\delta$ as a smoothed/regularized version of the classical depth function.

It is noted that some further properties readily follow from the ones stated below as, for example, monotonicity relative to the deepest point and maximality at the center for a centrally symmetric distribution.
We would like to emphasize that no assumptions are made on $\p$ and hence all the properties hold also for the empirical law $\p_n$.

\begin{proposition}[Standard properties]\label{prop:stdprops}
For any $\delta>0$ the robustified halfspace depth $D_\delta(\cdot|\p)$ satisfies the following properties:
\begin{itemize}
\item[P1] Isometric invariance: $D_\delta(A\z+\bs b|\p_{A\X+\bs b})=D_\delta(\z|\p_{\X})$ for any $\z,\bs b\in \R^d$ and any matrix $A\in \R^{d\times d}$ with $AA^\top=I$;
\item[P2] Null at infinity: $D_\delta(\z|\p)\to 0$ as $\norm{\z}\to\infty$. 
\item[P3] Convex upper level sets: the sets $\{\z\in\R^d: D_\delta(\z|\p)\geq \a\}$ are convex for every $\a\in(0,1]$.
%\item[P4] Maximality at the center: if $\X\stackrel{d}{=}-\X$ then $D_\delta(\bs 0|\p)=\sup_{\z}D_\delta(\z|\p)$.
\end{itemize}
\end{proposition}
%%%\begin{proof}
%%%P1. 
%%%Since every vector can be represented as $A\x+\bs b$ we have
%%%\[D_\delta(A\z+\bs b|\p_{A\X+\bs b})=\inf_{\norm{u}=1}\sup_{\p':d_W(\p'_{A\X+\bs b},\p_{A\X+\bs b})\leq \delta}\p'(\br{\u,(A\X+\bs b)-(A\z+\bs b)}\geq 0).\]
%%%From the definition of the Wasserstein distance and $\norm{(A\x-\bs b)-(A\y-\bs b)}=\norm{\x-\y}$ %and the fact that $A\X+\bs b\stackrel{d}{=}A\X'+\bs b$ iff $\X\stackrel{d}{=}\X'$ 
%%%we readily deduce that 
%%%\[d_W(\p'_{A\X+\bs b},\p_{A\X+\bs b})=d_W(\p'_{\X},\p_{\X}).\] 
%%%It is left to note that $\br{\u,(A\X+\bs b)-(A\z+\bs b)}=\br{A^\top\u,\X-\z}$, where the vector $A^\top \u$ runs over all unit vectors. The result now follows.
%%%
%%%P2. For every $\z\neq 0$ we pick $\u=\z/\norm{\z}$ and observe that 
%%%\[Y=\br{\u,\z-\X}=\norm{\z}-\br{\z,\X}/\norm{\z}\geq \norm{\z}-\norm{\X}.\] 
%%%According to Lemma~\ref{lem:mon} the depth is upper-bounded by the value $\tilde v$ corresponding to $\tilde Y= \norm{\z}-\norm{\X}\to \infty$ a.s.
%%% Conclude by Lemma~\ref{lem:infinite}(i).
%%%
%%%
%%%P3. It is sufficient to show that for any $\z_1,\z_2\in \R^d$ and any $t\in [0,1]$ we have
%%%\[D_\delta(t\z_1+(1-t)\z_2|\p)\geq \min(D_\delta(\z_1|\p),D_\delta(\z_2|\p)).\]
%%%Note that $\br{\u,\z_2-\z_1}\geq 0$ implies $\br{\u,\X-\z_2}\leq\br{\u,\X-(t\z_1+(1-t)\z_2)}$, whereas $\br{\u,\z_2-\z_1}\leq 0$ implies  $\br{\u,\X-\z_1}\leq\br{\u,\X-(t\z_1+(1-t)\z_2)}$.
%%%The result now follows from~\eqref{eq:minimax}.
%%%\end{proof}

\begin{proposition}[Additional properties]\label{prop:addprops}
For any $\delta>0$ the function $D_\delta(\cdot|\p)$ satisfies the following additional properties:
\begin{itemize}
\item[P4] Positivity: $D_{\delta}(\z|\p)>0$ for all $\z\in\R^d$;
\item[P5] Joint continuity: $D_{\delta}(\z|\p)$ is continuous in~$(\z,\delta)\in\R^d\times(0,\infty)$;
\item[P6] Level sets: for any $\a\in(0,1)$ the set $\{\z\in\R^d: D_\delta(\z|\p)=\a\}$ is the boundary of the respective upper level set $\{\z\in\R^d: D_\delta(\z|\p)\geq \a\}$, or both are empty.
\end{itemize}
\end{proposition}

From the definition~\eqref{eq:minimax} it is clear that the depth $D_\delta(z|\p)$ is non-decreasing in $\delta\geq 0$.
Finally, we show that our robustified depth $D_\delta(\z|\p)$ approaches the halfspace depth as $\delta\downarrow 0$. Note that uniform convergence fails in general, since  $D_\delta(\z|\p)$ is continuous in $\z$ for every $\delta>0$ according to P5, whereas the limit function $D_0(\z|\p)$ is not necessarily continuous.
\begin{proposition}[Tukey depth approximation]\label{prop:approx}
For all $\z\in \R^d$ it holds that 
\[D_\delta(\z|\p)\downarrow D_0(\z|\p)\qquad \text{ as }\quad\delta\downarrow 0.\]
\end{proposition}
%%%\begin{proof}
%%%It is sufficient to show for any direction $\u$ that 
%%%\[\sup_{d_W(\p',\p)\leq \delta}\p'(\br{\u,\X-\z}\geq 0)-\p(\br{\u,\X-\z}\geq 0)\downarrow 0\]
%%%as $\delta\downarrow 0$.
%%%We may assume that $\p(Y>0)>0$ since otherwise both terms are equal to~1. For $\delta>0$ small enough we have $\e Y^+>\delta$ and so $\lambda_\delta=h^{-1}(\delta)<\infty$.
%%%Now the solution in~\eqref{eq:sol_alt} gives the representation of the above difference:
%%%\begin{equation}\label{eq:approx}\p(Y\in(0,\lambda_\delta])-(h(\lambda_\delta)-\delta)/\lambda_\delta,\qquad (h(\lambda_\delta)-\delta)/\lambda_\delta\in[0,\p(Y=\lambda_\delta)).\end{equation}
%%%This indeed converges to~0 if $\lambda_\delta\to 0$.
%%%Assume that $\lambda_\delta\to\lambda>0$ and note that necessarily $\p(Y\in(0,\lambda))=0$.
%%%If $Y$ has no mass at $\lambda$ then we are done, and if such mass is positive then it must be that $(h(\lambda_\delta)-\delta)/\lambda_\delta$ converges to this mass.
%%%Hence the expression in~\eqref{eq:approx} goes to 0 and the proof is  complete.
%%%\end{proof}

For any $\delta>0$ consider the maximal depth
\[\oa(\delta)=\max_{\z\in\R^d} D_\delta(\z|\p)\in(0,1],\]
which is achieved according to P2 and P5. 
\begin{lemma}\label{lem:oa}
The maximal depth $\oa(\delta)\in(0,1]$ is non-decreasing and continuous for $\delta>0$. 
\end{lemma}
%%%\begin{proof}
%%%Our proof relies on the joint continuity of the depth function, see P5.
%%%Consider $\delta_n\uparrow \delta$ and let $\z$ be a point achieving the maximal depth for~$\delta$.
%%%Now $D_{\delta_n}(\z|\p)\to D_\delta(\z|\p)$ and the latter can not be exceeded for ambiguity radius~$\delta_n$, showing left-continuity of~$\oa$. 
%%%Let $\delta_n\downarrow \delta$ with $\z_n$ being points with maximal depth. These $z_n$ must belong to a compact set, see P2 and use monotonicity of depth in~$\delta_n$. 
%%%Thus $(z_n)$ contains a convergent subsequence establishing right-continuity of $\oa$.
%%%\end{proof}

The median
\begin{equation}\label{eq:median}M_\delta=\{\z\in\R^d:D_\delta(\z|\p)=\oa(\delta)\}\end{equation} is a convex region 
which is, in fact, an interval (commonly a single point) if $\oa<1$, see P6. %example: vertecis of a rectangle
The case $\oa=1$ is special as it allows for median regions with a non-empty interior, and it also leads to a simple characterization:
\begin{equation}\label{eq:median1}M_\delta=\{\z\in\R^d:\sup_{\norm{u}=1}\e \br{\u,\z-\X}^+\leq \delta\},\qquad \text{ when }\oa(\delta)=1,\end{equation}
see Proposition~\ref{prop:sup}.
In words, such points $\z$ have the property that the expected distance of $\X$ to any halfspace anchored at $\z$ is at most $\delta$.
The median region for the empirical distribution is further studied in \S\ref{sec:representations}. 

Let us note that one should be careful when approximating the median regions. That is for a fixed $\delta>0$ the weak convergence $\X_n\convd \X$ does not in general imply that the median regions $M_{n\delta}$ well approximate~$M_\delta$.
In particular, $M_\delta$ may be a set with a non-empty interior, whereas all $M_{n\delta}$ can be single points or intervals which happens, for example, when $\br{\u,\X_n}^+$ are non-integrable for all $\u\neq \bs 0$ and all~$n$. 
It seems plausible that the upper level sets for $\a\in(0,\oa)$ do converge with respect to the Hausdorf distance, say. Such approximation results are beyond the scope of this paper.
%Its proofThis is beyond the scope of the present paper.

%\subsection{Maximin formulation}
%Here we show that our robustified depth $D_\delta(z)$ defined in~\eqref{eq:minimax} can be alternatively written as the maximin problem.
%\begin{proposition}
%For every $\delta\geq 0$ and $z\in\R^d$ there is the identity
%\[D_\delta(z)=\sup_{\p':d(\p',\p)\leq \delta}\inf_{u\in\R^d:\|u\|=1}\left\{\p'(\br{u,X-z}\geq 0)\right\}.\]
%\end{proposition}
%\begin{proof}
%It is only required to show the inequality
%\[\sup_{\p':d(\p',\p)\leq \delta}\inf_{u\in\R^d:\|u\|=1}\left\{\p'(\br{u,X-z}\geq 0)\right\}\geq D_\delta(z),\]
%because the reversed inequality is standard, see the argument in~\cite[Lem.\ 36.1]{rockafellar}.
%Let $u^*$ be an optimal direction (its existence is proven in ...) and let $X^*$ be the random such that $(X,X^*)$ provides the optimal coupling for the direction~$u^*$.
%In other words, the depth is given by $\p(\br{u^*,X^*-z}\geq 0)$. Furthermore, the maximin problem is lower bounded by
%\[\inf_{u\in\R^d:\|u\|=1}\left\{\p(\br{u,X^*-z}\geq 0)\right\}\] and so it is left to show that all such probabilities are lower bounded by the one with $u=u^*$.
%This is not true!
%\end{proof}
Let us now analyze the case of the standard normal vector $\X$ in more detail.
\begin{example}[Standard normal]\label{ex:normal}\rm
%Due to translation invariance of the depth function, see P1, we 
Assume that $\X$ is a standard $d$-dimensional normal vector. %with covariance matrix $\Sigma$ and zero mean. Additionally, we assume that $\Sigma$
%That is, we may write $\X=A\Z$ for a $d$-dimensional standard normal $\Z$ and $d\times d$ matrix $A$ satisfying $AA^\top=\Sigma$.
Consider $Y=\br{\u,\z-\X}\sim\mathcal N(\br{\u,\z},1)$ when $\norm{\u}=1$. Assume for the moment that $\z\neq \bs 0$ and note that $\br{\u,\z}$ is maximal when $\u=\z/\norm{\z}$, and thus this is an optimal direction according to Lemma~\ref{lem:mon}. For $\z=\bs 0$ any direction is optimal. 
It is left to analyze the random variable $Y=\norm{\z}+Z$ for a standard normal $Z$. 

Let $\Phi,\varphi$ be the normal cumulative distribution function and the density, respectively.  
According to Proposition~\ref{prop:sup} the depth is given by 
\[D_\delta(\z|\p)=\Phi(\lambda-\norm{\z}),\qquad \text{where }\int_0^{\lambda}x\varphi(x-\norm{\z})\D x=\delta\]
with $\lambda=\infty$ (leading to 1) when the solution does not exist; recall that the standard Tukey depth corresponds to $\lambda=0$.
In particular, the upper level sets are balls centered at the origin. The maximal depth attains 1 iff $\delta\geq \e Z^+=1/\sqrt{2\pi}\approx 0.4$.

One can also calculate the maximal depth $\oa(\delta)<1$ for $\delta<1/\sqrt{2\pi}$, which must be attained by $\z=\bs 0$. Thus we need to solve 
$\varphi(0)-\varphi(\lambda)=\delta$ according to the formula for the expectation of the truncated normal.
This readily gives
\[\oa(\delta)=\Phi\Big(\sqrt{-2\log(1-\delta\sqrt{2\pi})}\Big)\in(1/2,1),\qquad \delta\in(0,1/\sqrt{2\pi}),\]
see Figure~\ref{fig:astar} below for an illustration.

\end{example} 

\subsection{Decay at infinity}
The property P2 states that the robustified depth becomes 0 at infinity. Here we study its rate of decay and compare it to the halfspace depth.
First, we provide an upper bound.
\begin{lemma}\label{lem:decay_bound}
If $\e\norm{\X}^p<\infty$ for some $p\in(0,1]$ then for $\delta>0$ there exists a constant $c>0$ such that 
\begin{align*}
&D_\delta(\z|\p)\leq c\norm{\z}^{-p}
\end{align*}
for all $\z\in\R^d\backslash\{\bs 0\}$.
\end{lemma}
%%%\begin{proof}
%%%Take $\u=\z/\norm{\z}$ assuming $\z\neq \bs 0$ and reconsider the optimal transport plan $\X^*$ satisfying~\eqref{eq:opt_transport}.
%%%Note that $\br{\u,\X^*-\z}\leq \norm{\X^*}-\norm{\z}$, and so we have the bound
%%%\[D_\delta(\z|\p)\leq \p(\norm{\X^*}-\norm{\z}\geq 0)\leq \frac{\e\norm{\X^*}^p}{\norm{\z}^p},\]
%%%where $\e\norm{\X^*-\X}\leq \delta$; here we have also used Markov's inequality.
%%%Using $(a+b)^p\leq a^p+b^p$ for $p\leq 1$ we find that 
%%%\[\e \norm{\X^*}^p\leq \e(\norm{\X^*-\X}+\norm{\X})^p\leq \e\norm{\X^*-\X}^p+\e\norm{\X}^p.\]
%%%It is left to observe that $\e\norm{\X^*-\X}^p\leq \e(1+\norm{\X^*-\X})\leq (1+\delta)$.
%%%\end{proof}

Next, we show that the decay can never be faster than linear.
\begin{lemma}\label{lem:decay_sublin}
For $\delta>0$ it holds that
\[\liminf_{\norm{\z}\to\infty}\norm{\z}D_\delta(\z|\p)\geq \delta.\]
\end{lemma}
%%%\begin{proof}
%%%Choose a constant $r>0$ large enough so that $\p(\norm{\X}\leq r)\geq 1/2$, and let $B_r$ be the ball of radius $r$ centered at the origin.
%%%For any point $\z$ and any direction $\u$ we must have 
%%%\[\br{\u,\z-\x}\leq \norm{\z}+r,\qquad \forall \x\in B_r.\]
%%%Thus $Y=\br{\u,\z-\X}\leq \norm{\z}+r$ with probability $p\geq 1/2$, and according to~\eqref{eq:lowerbound}
%%%we find
%%%\[D_\delta(\z|\p)\geq \min(1/2,\delta/(\norm{\z}+r)).\]
%%%%where $\hat v$ corresponds to $\hat Y=\norm{\z}+r$, which does not depend on the direction. But $\hat v(\delta)=\frac{\delta}{\norm{\z}+r}\wedge 1$ which follows easily from~\eqref{eq:sol_alt}.
%%%This lower bound behaving asymptotically as $\delta/\norm{\z}$ and the result follows.
%%%\end{proof}

Thus we find that for an integrable $\X$ the depth decays linearly. It will be shown in Proposition~\ref{prop:outer} that the empirical depth (which is necessarily compactly supported) has the exact rate $\delta/\norm{\z}$.
\begin{corollary}[Linear decay]\label{cor:lindec}
Assume $\e \norm{\X}<\infty$. Then for $\delta>0$ there exist two constants $0<c_1<c_2$ such that 
\[c_1\norm{\z}^{-1}\leq D_\delta(\z|\p)\leq c_2\norm{\z}^{-1}\] 
when $\norm{\z}$ is large enough.
\end{corollary}
%%%\begin{proof}
%%%Combine Lemma~\ref{lem:decay_bound} and Lemma~\ref{lem:decay_sublin}.
%%%\end{proof}

In the following we discuss a simple class of examples where the decay is slower than linear. Necessarily, such must have infinite first moment $\e \norm{\X}=\infty$.
Here we rely on some basic theory of regular variation~\citep{BGT}.
%For simplicity we assume that $X$ is regularly varying with index $-\a$, where $\a\in(0,1)$. According to~\cite{basrak} this assumption (standard in extreme value theory) can be stated as 
%\begin{equation}\label{eq:assumption_MRV}\lim_{t\to\infty}\frac{\p(\br{u,X}>t)}{t^{-\a}\ell(t)}\to w(u)\end{equation}
%for all $u\in\mathbb S^{d-1}$, where $w(u)>0$ for some $u\neq 0$ and $\ell$ is a slowly varying function. It is noted that such $X$ satisfies $\e\norm{X}^p<\infty$ for $p<\a$, and  $\e\norm{X}^p=\infty$ for $p>\a$.
%In particular, Lemma~\ref{lem:decay_bound} gives 
%\[\lim_{\norm{z}\to\infty}\norm{z}^pD_\delta(z)=0\] for $p<\a$.

\begin{lemma}\label{lem:spheredist}
Assume that $\X$ has a spherical distribution and $\p(\br{\u,\X}> t)=t^{-\a}\ell(t)$ with $\a\in (0,1)$ and some slowly varying (at infinity) function~$\ell$.
Then $D_\delta(\z|\p)\sim\norm{\z}^{-\a}\ell(\norm{\z})$ as $\norm{\z}\to\infty$ for any $\delta\geq 0$.
In particular, 
\[\lim_{\norm{\z}\to\infty}\norm{\z}^pD_\delta(\z|\p)=\begin{cases}0, &p<\a,\\ \infty, &p>\a.\end{cases}\]
\end{lemma}
%%%\begin{proof}
%%%Consider $Y=\br{\u,\z-\X}\stackrel{d}{=}\br{\u,\z}-Z$, where $Z$ has the distribution of $\br{\u,\X}$; the same for all $\u$ by assumption.
%%%The largest $\br{\u,\z}$ equals $\norm{\z}$, and according to Lemma~\ref{lem:mon} it gives the smallest value, and hence the depth. According to
%%%\eqref{eq:sol_alt} we then have
%%%\[D_\delta(\z|\p)\geq \p(\norm{\z}-Z< \lambda)=\p(Z>\norm{\z}-\lambda)=(\norm{\z}-\lambda)^{-\a}\ell(\norm{\z}-\lambda),\]
%%%where $\lambda=\lambda(\norm{\z})$ is defined using the random variable~$\norm{\z}-Z$.
%%%For the lower bound it is thus left to show that $\lambda/\norm{\z}\to 0$, see the uniform convergence theorem~\cite[Thm.\ 1.2.1]{BGT}.
%%%For the upper bound we use $\p(\norm{\z}-Z\leq \lambda)$ resulting in the same asymptotic behavior.
%%%
%%%On the contrary, assume that for some $\epsilon>0$ we have $\lambda>4\epsilon\norm{\z}$ for some $\z$ with arbitrarily large norm. For simplicity we write $\eta=\norm{\z}$.
%%%From the definition of $\lambda$ we then must have
%%%\[\delta\geq \e\big((\eta-Z)\ind{\eta-Z\in (0,2\epsilon \eta]}\big)\geq \epsilon\eta\p(Z\in [(1-2\epsilon) \eta,(1-\epsilon)\eta)]).\]
%%%The latter probability is asymptotic to $((1-2\epsilon)^{-\a}-(1-\epsilon)^{-\a})\eta^{-\a}\ell(\eta)$, and so for $\a\in(0,1)$ the right hand side in the display increases to $\infty$ as $\eta\to\infty$, which is a contradiction.
%%%\end{proof}

In conclusion, the main difference from the halfspace depth arises in the integrable case, where the traditional depth may decay at an arbitrarily fast rate.
One particular example concerns empirical distributions, where the halfspace depth is 0 for large enough $\norm{\z}$, whereas the robustified depth decays as $\delta/\norm{\z}$. 
In the non-integrable case one may possibly distinguish between the directions along which $\z$ escapes to infinity, since potentially the decay may be drastically different for different directions.

\subsection{Maximin formulation}\label{sec:maxmin}
One may be interested in swapping minimax in~\eqref{eq:minimax} for the maximin formulation.
By the standard max-min inequality we have
\begin{equation}\label{eq:maxmin}\sup_{\p':d_W(\p',\p)\leq \delta}\left\{\inf_{\u\in\R^d:\|\u\|=1}\p'(\br{\u,\X-\z}\geq 0)\right\}\leq D_\delta(\z|\p),\end{equation}
and the obvious question is if the two formulations coincide.
The answer is no in general, and we provide a simple example with the strict inequality below.

We take $d=2,\z=\bs 0$ and consider 
\begin{equation}\label{eq:maxmin_p}\p=\frac{1}{2}\delta_{\{(1,1)\}}+\frac{1}{2}\delta_{\{(1,-1)\}}\end{equation} putting half of the mass at the point $(1,1)$ and the other half at $(1,-1)$.
For simplicity we assume that $\delta\in(0,1/2)$.
It is not difficult to see from Corollary~\ref{cor:sample} that the optimal direction is $\u=(-1,0)$ resulting in $D_\delta(\bs 0|\p)=\delta$. 
%The optimal law is not unique, and one such puts mass $\delta/2$ at $(0,\pm 1)$ and $(1-\delta)/2$ at $(1,\pm 1)$.

\begin{lemma}\label{lem:maxmin}
For $\delta\in(0,1/2)$ and $\p$ in \eqref{eq:maxmin_p} it holds that 
\[\sup_{\p':d_W(\p',\p)\leq \delta}\left\{\inf_{\u\in\R^d:\|\u\|=1}\p'(\br{\u,\X}\geq 0)\right\}=\frac{1}{\sqrt 2}\delta<\delta=D_\delta(\bs 0|\p).\]
\end{lemma}

The main problem in using maximin formulation is that it does not seem to yield an explicit solution in general unlike the robustified depth suggested in this paper.
Furthermore, our depth concept is in accordance with the modern approach of dealing with the optimizer's curse in stochastic programs~\citep{esfahani_kuhn}.

\section{Finite sample version}\label{sec:sample}
Throughout this section we assume that $\p=\p_n=\frac{1}{n}\sum_{i=1}^n\delta_{\{\x_i\}}$ is the empirical measure corresponding to the $n$ iid observations of~$\X$.
\subsection{Consistency}

We have a consistency result, uniform in both the point $\z$ and the ambiguity radius $\delta$. %for any sequence $\delta_n\geq 0$.
\begin{proposition}\label{prop:consistency}
For $\p_n$ obtained by independently sampling from $\p$ it holds that
\[\sup_{\z\in\R^d,\delta\geq 0}|D_{\delta}(\z|\p_n)-D_{\delta}(\z|\p)|\to 0\quad\text{a.s.}\]
as $n\to\infty$.
\end{proposition}
%%%\begin{proof}
%%%Firstly, the case $\delta=0$ is classical~\citep{donoho_gasko}.
%%%According to $|\inf f-\inf g|\leq \sup|f-g|$ and the representation~\eqref{eq:sol} of the depth for $\delta>0$, it is sufficient to show that 
%%%\[\sup_{\z,\u\in\R^d}|\e_n(1-\br{\u,\z-\X}^+)^+-\e(1-\br{\u,\z-\X}^+)^+|\to 0,\]
%%%where the division by $\lambda'$ is incorporated into~$\u$. 
%%%
%%%Now the result follows from empirical process theory~\cite[Thm.\ II.24 and Lem.\ II.25]{pollard}. It is only needed to observe that the graph of the function $f_{\u,\z}(\x)=(1-\br{\u,\z-\x}^+)^+\in[0,1]$ can be constructed from three half-spaces in $\R^{d+1}$ using union and intersection operations, and so the respective class has polynomial discrimination~\cite[Lem.\ II.15]{pollard}.
%%%\end{proof}

Upon recalling the continuity of depth in the ambiguity radius we readily obtain the basic consistency result.
\begin{corollary}\label{cor:convdelta}
For any $\z\in\R^d$ and any sequence $\delta_n\downarrow 0$ it holds that
\[D_{\delta_n}(\z|\p_n)\to D_0(\z|\p)\quad \text{a.s.}\]
as $n\to \infty$.
\end{corollary}
%%%\begin{proof}
%%%Consider the bound
%%%\[|D_{\delta_n}(\z|\p_n)- D_0(\z|\p)|\leq |D_{\delta_n}(\z|\p_n)- D_{\delta_n}(\z|\p)|+|D_{\delta_n}(\z|\p)- D_0(\z|\p)|\]
%%%and observe that both terms converge to 0 according to Proposition~\ref{prop:consistency} and Proposition~\ref{prop:approx}.
%%%\end{proof}

It must be noted that a similar result follows from the general theory in \cite[Thm.\ 3.6]{esfahani_kuhn} based on the measure concentration result in~\cite{fournier2015rate}.
This result, however, relies on a somewhat restrictive assumption $\e e^{\norm{\X}^a}<\infty$ for some $a>1$, which is slightly stronger than saying that $\norm{\X}$ is light-tailed.
Furthermore, it is required there that $\delta_n$ decays sufficiently slowly.
Additionally, \cite{esfahani_kuhn} provide a finite sample guarantee.

\subsection{Breakdown point}
\begin{figure}
\centering
\includegraphics[width=0.4\textwidth]{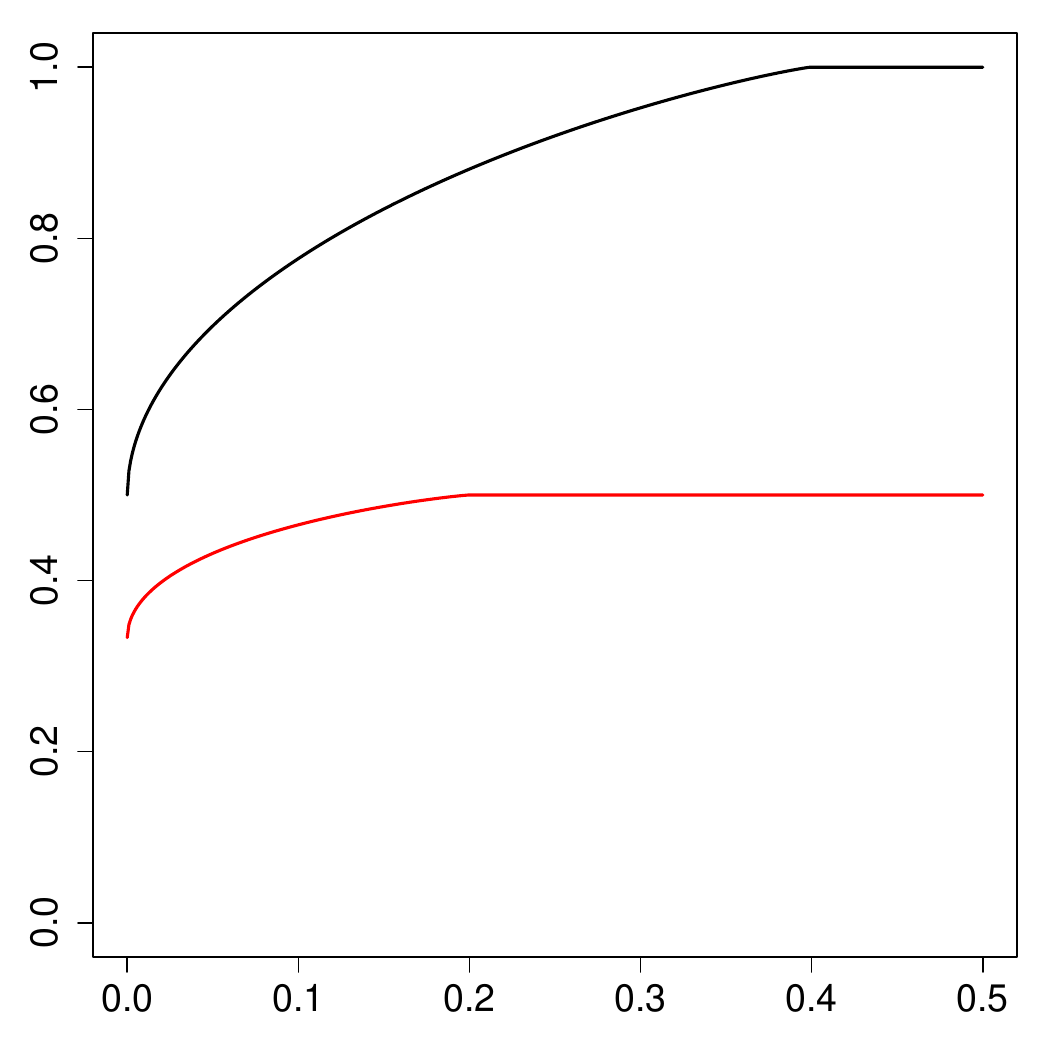}
\caption{$\oa(\delta)$ and $\a^*(\delta)$ (red) for the standard normal in an arbitrary dimension.}
\label{fig:astar}
\end{figure}

We start by defining an important quantity $\a^*=\a^*(\delta)\in(0,1/2]$ using the equation
\begin{equation}\label{eq:astar}\frac{\a^*}{1-\a^*}=\oa\Big(\frac{\delta}{1-\a^*}\Big),\end{equation}
where $\oa$ denotes the maximal depth.
Figure~\ref{fig:astar} illustrates the maximal depth $\oa$ and also $\a^*$ for the standard normal distribution in an arbitrary dimension~$d$, see Example~\ref{ex:normal}.
Let us summarize the basic properties of~$\a^*$. 
\begin{lemma}\label{lem:astar}
Let $\delta\geq 0$.
The equation~\eqref{eq:astar} has a unique solution $\a^*$ in $(0,1)$. Moreover,  $\a^*\in(0,1/2]$ with $\a^*=1/2$ iff $\oa(2\delta)=1$.
For any $\a\in[0,1)$ it holds that
\[\frac{\a}{1-\a}<\oa\Big(\frac{\delta}{1-a}\Big)\qquad\text{iff}\qquad \a<\a^*.\]
Furthermore, $\a^*(\delta)\leq \oa(\delta)$ and $\a^*(\delta)$ is non-decreasing in $\delta\geq 0$.
\end{lemma}
%%%\begin{proof}
%%%According to~\eqref{eq:inc_v} the function $p\oa(\delta/p)$ is non-decreasing in $p\in[0,1]$, which follows by considering all points $\z$ and all directions~$\bs u$; it is assumed to be 0 at~0.
%%%Hence $p\oa(\delta/p)-(1-p)$ is strictly increasing with values $-1$ and $\oa(\delta)>0$ at the end points.
%%%This function inherits continuity from $\oa$, see Lemma~\ref{lem:oa}.
%%%Hence there is a unique root, proving the first claim. The equivalence of inequalities follows from the strict monotonicity.
%%%Finally, $\a^*=1/[1+1/\oa(\delta/(1-\a^*))]$ which must be in $(0,1/2]$ since $\oa\leq 1$, and the value $1/2$ corresponds to $\oa(\delta/(1-1/2)=1$.
%%%The fact that $\a^*(\delta)$ is non-decreasing is inherited from the same property of $\oa(\delta)$.
%%%Note also that $(1-\oa)\oa(\delta/(1-\oa))\leq \oa$ and hence $\a^*\leq \oa$.
%%%\end{proof}

Consider the upper level set 
\[U_\delta(\a)=U_\delta(\a|\p_n)=\{\z\in\R^d:D_\delta(\z|\p_n)\geq \a\}\] 
and let $d_H(A,B)$ be the Hausdorff distance between two compact sets (the maximal distance from a point in one set to the closest point in the other set).
We let $d_H(A,B)=\infty$ if $A$ or $B$ is empty.
The breakdown point of $U_\delta(\a)$ for a given $\p_n,\delta,\a$ is defined~\citep{donoho_gasko,nagy_survey} by
\[BP(U_\delta(\a),\p_n)=\min_{m\in\mathbb N}\left\{\frac{m}{m+n}:\sup_{\y_1,\ldots\y_m\in\R^d}d_H(U_\delta(\a|\p_n),U_\delta(\a|\p_{n+m}))=\infty\right\},\]
where $\p_{n+m}=\frac{1}{n+m}(\sum_{i=1}^n\delta_{\{\x_i\}}+\sum_{i=1}^m\delta_{\{\y_i\}})$.
In words, it is the minimal proportion of new (contaminating) observations which can make the upper level set arbitrarily different from the given one. 
In the case of the classical halfspace depth $\delta=0$ the breakdown point is given by~\cite{nagy_illumination,donoho_gasko}, and here we get a similar result. 
%The difference is in a slightly relaxed upper bound on~$\a$, see also the comments following Proposition~\ref{prop:breakdown}.
Note, however, that the robust depth is at least as large as the traditional Tukey depth and so the effective $\a$ in our case (resulting in a similar upper level set) is larger, making our depth more robust to contamination.
In this regard, we also mention that our upper bound on~$\a$ is larger.
We write $\a^*(\delta,\p_n)$ for the solution of~\eqref{eq:astar} with respect to the measure~$\p_n$ and radius~$\delta$.
\begin{proposition}\label{prop:bpastar}
For any $\delta\geq 0$, $\p_n$ and $0<\a\leq\a^*(\delta,\p_n)$ it holds that
\[BP(U_\delta(\a),\p_n)=\frac{\lceil n\a/(1-\a)\rceil}{n+\lceil n\a/(1-\a)\rceil}.\]

Moreover, for $0<\a<\a^*(\delta,\p)$ we also have
\[BP(U_\delta(\a),\p_n)\to \a \]
as $n\to\infty$ on a set of measure one.
\end{proposition}

\begin{figure}[h!]
\begin{tabular}{ccc}
& Traditional Tukey depth $D_0(\cdot|\p_{n+m})$ & Robust Tukey depth $D_{0.1}(\cdot|\p_{n+m})$ \\
\begin{sideways}\qquad\qquad Contamination $m=30$\end{sideways} &\includegraphics[width=0.47\textwidth,trim=0 0 0 1.5cm,clip=true,page=2]{picMaps-robMed-n120.pdf} & 
\includegraphics[width=0.47\textwidth,trim=0 0 0 1.5cm,clip=true,page=1]{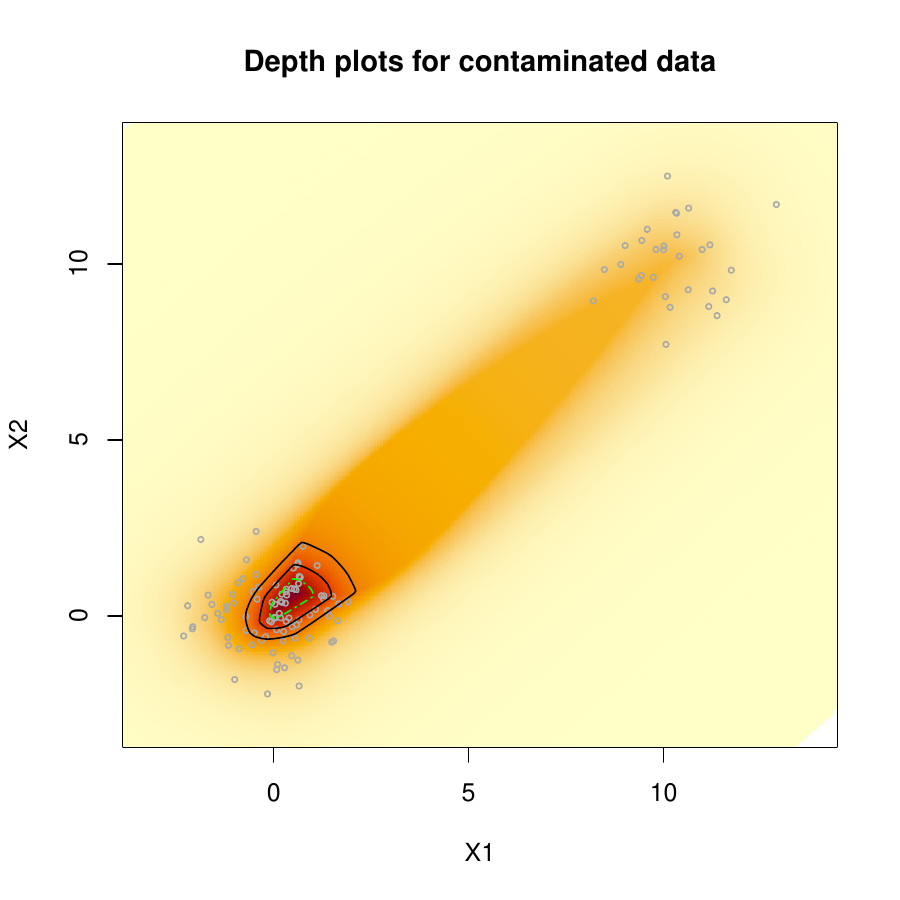} \\
\begin{sideways}\qquad\qquad Contamination $m=50$\end{sideways} &\includegraphics[width=0.47\textwidth,trim=0 0 0 1.5cm,clip=true,page=4]{picMaps-robMed-n120.pdf} & 
\includegraphics[width=0.47\textwidth,trim=0 0 0 1.5cm,clip=true,page=3]{picMaps-robMed-n120.pdf}
\end{tabular}
\caption{Heat maps of the traditional Tukey depth $D_0(\cdot|\p_{n+m})$ (left) and the proposed robust Tukey depth $D_{0.1}(\cdot|\p_{n+m})$ (right) for a sample of $n+m=120$ points drawn from a bivariate standard  normal distribution except contaminating $m=30$ (top) and $m=50$ (bottom). Depth contours are plotted for depth values $(\frac{4}{7},\frac{5}{7},\frac{6}{7},1)\times\oa(\cdot,\p_{n+m})$. The contour of depth $\frac{6}{7}$ is plotted with dashed green line.}
%%Robust Tukey depth contours are plotted for depth values   $0.35,0.45,0.5$ with $\delta=0.1$; traditional Tukey depth contours are plotted for depth values $0.25,0.35,0.4\,(\text{for }m=50\text{ only})$.}
\label{fig:breakdown}
\end{figure}

Note that $\a^*(0)=\oa(0)/(1+\oa(0))$ and we recover the upper bound on $\a$ in the classical setting.
The only difference in the breakdown point of the upper level set for the robustified and classical halfspace depths is in the allowed range of the level~$\a$, 
but see the above comment on the effective~$\a$ implying that our depth is more robust; see also Figure~\ref{fig:breakdown} for an illustration.
This result can be reformulated to show that for large $\a$ (exceeding the above bound) the breakdown occurs at $\a^*(\delta,\p)$ due to the empty set problem. % and then again our setting is more robust than the classical one.
We note that the level set for $\a=\oa$ breaks down with a single additional point in the sense that such level is not achieved in the new data set.
A more relevant statistic is the median $M_\delta$ itself, defined in~\eqref{eq:median}, and we study it below.
%It will be shown in the following that the median does not break down at less than $\oan$ additional observations. The exact breakdown depends on the data set~$\p_n$.

\begin{proposition}\label{prop:breakdown}
For any $\delta\geq 0$ and $\a^*_n=\a^*(\delta,\p_n)$ it holds that
\[BP(M_\delta,\p_n)\geq \frac{\lceil n\a^*_n/(1-\a^*_n)\rceil}{n+\lceil n\a^*_n/(1-\a^*_n)\rceil}.\]
%where $p^*_n\in(0,1/2]$ is the unique solution of $(1-p^*_n)/p^*_n=q^*(\delta/p^*_n,\p_n)$.
If $\oa(2\delta,\p_n)=1$ then $BP(M_\delta,\p_n)= 1/2$.

Furthermore, 
\[\liminf_{n\to\infty} BP(M_\delta,\p_n)\geq \a^*(\delta,\p)\]
%where $p^*$ is the unique solution of $(1-p^*)/p^*=q^*(\delta/p^*,\p)$. 
on a set of measure one, where the limit is $1/2$ if $\oa(2\delta,\p)=1$.
\end{proposition}
%%%\begin{proof}
%%%Suppose there exists $\z$ outside of the convex hull of $\x_i$, and such that it belongs to the new median region.
%%%Its depth is no more than $p+\delta/d$ with $p=m/(n+m)$, where $d$ is the distance from the convex hull, see Lemma~\ref{lem:mix} as well as Lemma~\ref{lem:mon} with~\eqref{eq:det}; alternatively one may use Lemma~\ref{lem:infinite}.
%%%The points in the old median region have new depth at least $(1-p)\oa(\delta/(1-p),\p_n)$, and hence we must have 
%%%\[p\geq (1-p)\oa(\delta/(1-p),\p_n).\]
%%%According to Lemma~\ref{lem:astar} we must have $p\geq \a^*_n$ and then $m\geq n\a^*_n/(1-\a^*_n)$.
%%%
%%%%The equivalence of $p^*=1/2$ and $q^*(2\delta)=1$ follows readily from the definition of~$p^*$.
%%%Let us now show that $m=n$ new points are sufficient for the breakdown, so that the breakdown point is exactly $1/2$ when $\a_n^*=1/2$.
%%%Suppose breakdown does not occur, and so the new median must be contained in some inflation of the given convex hull. Now we copy the given $n$ points and shift them sufficiently far away, so that the two inflations do not intersect.
%%%But we may regard the new set of points as the original data set, which readily leads to a contradiction.
%%%The final statement follows from the convergence $\a^*_n\to \a^*(\delta,\p)$ a.s.
%%%\end{proof}

Importantly, $\a^*(\delta)$ is non-decreasing  and so the breakdown bound grows with~$\delta$. 
Furthermore, any median region with non-empty interior has a breakdown point of~$1/2$, which is considerably better than in the case of the classical halfspace depth in general. 
%\ji{References? I see an example with $1/4$: just take 3 points}
In the case of a multivariate standard normal $\oa(2\delta)=1$ occurs when $\delta\geq 2^{-3/2}\pi^{-1/2}\approx 0.2$, see Figure~\ref{fig:astar}.
In this case the median has a breakdown point of $1/2$, even though the maximal depth can be as small as~0.881. 

%We provide a numerical illustration of breakdown for a contaminated bivariate normal distribution in Figure~\ref{fig:breakdown}.

%\subsection{Can the median breakdown point exceed the lower bound?}
%An example where the breakdown point exceeds the given lower bound. Dimension $d=1$ is easy but such an example may seem degenerate. 
%For $d=2$ we consider an equilateral triangle and points masses $1/3$ at its vertices. Find $\delta$ such that $p^*=3/5$ and so the bound becomes $m=2$.
%Let us show that breakdown does not happen.
%Suppose $z$ arbitrarily far away belongs to the median. Then the maximal distance of the points on the triangle can not exceed~$2/3$ (maybe less if the depth at $z$ is less).
%Suppose all the new mass is put at~$z$. Furthermore, define $z'$ as the vertex in the convex hull of $z$ and the triangle if such exists, or the point on the edge at its intersection with the line passing through $z$ and the other vertex. It is not difficult to see that the depth of $z'$ is larger than the depth at far removed $z$.
%\ji{What happens if the depth at $z$ is smaller than $2/3$? This is quite a bit messy. Maybe finish here or give an simulated example without claiming that the bound is not tight?}
%\ji{There should be a similar problem for the Tukey median, no?}

\subsection{Towards affine-invariant depth}

Different to the traditional halfspace depth, the robustified version described until now is isometrically invariant only. In other words, the proposed depth is preserved under translation, rotation and reflection. Absence of general affine invariance comes as the price for distributional robustness in a natural way, due to weaker invariance properties of the Wasserstein distance. If affine invariance is a necessary requirement, provided that the first two moments $\bs{\mu}_{\X}$ and $\Sigma_{\X}^{-1}=\Lambda_{\X}\Lambda_{\X}^\top$ are known or (which is more realistic) can be reliably estimated, affine invariance can be implemented according to transformation/re-transformation procedure~\citep{ChakrabortyC96}: consider the standardized $\tilde \X=\Lambda_{\X}(\X-\bs{\mu}_{\X})$ and the respective $\tilde \z=\Lambda_{\X}(\z-\bs{\mu}_{\X})$ and calculate the depth of $\tilde\z$ w.r.t.\ the law of~$\tilde\X$.
% $D_\delta(\tilde\z|\p_{\tilde\X})=D_\delta(A\z+\bs b|\p_{A\X+\bs b})$ for any $\z,\bs b\in \R^d$ and any $A\in \R^d\times\R^d$ with $AA^\top=I$.

To avoid unnecessary transformations and maintain practitioner's convenience, we formulate the affine-invariant version of the robustified halfspace depth by incorporating affine invariance in the distribution dissimilarity used in the inner problem. This yields the following formulation:

\begin{equation}\label{eq:minimaxaffinv}
	D_{\delta(\Sigma)}(\z|\p)=\inf_{\u\in\R^d:\|\u\|=1}\left\{\sup_{\p':d_W(\p',\p)\leq \delta}\p'(\br{\Lambda\u,\X-\z}\geq 0)\right\}\,.
\end{equation}

\begin{proposition}\label{prop:affinv}
	Let $\p$ possess second moment and positive definite covariance matrix $\Sigma_X$, then $D_{\delta(\Sigma)}$ defined in equation~\eqref{eq:minimaxaffinv} is affine invariant, \textit{i.e.} satisfies
	\begin{itemize}
		\item[P1A] $D_{\delta(\Sigma)}(A\z+\bs b|\p_{A\X+\bs b}) = D_{\delta(\Sigma)}(\z|\p_{\X})$ for any $\z,\bs b\in \R^d$ and any non-singular matrix $A\in \R^{d\times d}$.
	\end{itemize}
\end{proposition}
%%%\begin{proof}
%%%	The proof is straightforward noting $\V=A\X + \bs b$ and noticing that, for each $\u\in\R^d$, $\|u\|=1$, the following holds
%%%	\begin{equation*}
%%%		\br{\Lambda_{V}\u,\V} = (\Lambda_{V}\u)^\top\V = \bigl((A^\top)^{-1}\Lambda_X\u\bigr)^\top(A\X) = \u^\top\Lambda_X^\top A^{-1}AX = \br{\Lambda_X\u,X},
%%%	\end{equation*}
%%%	and the same for $\z$; the second equality holds because $(\Sigma_V)^{-1}=(A^\top)^{-1}\Lambda_X\Lambda_X^\top A^{-1}$.
%%%\end{proof}

From first sight, this affine invariance comes at cost of robustness, but this is the case only if non-robust (\textit{e.g.}, moment) estimate is used for $\Sigma$. To avoid breakdown, robust estimates for $\Sigma$ should be used instead, maintaining breakdown level not lower than the required one, such as, \textit{e.g.}, minimum covariance determinant estimator~\citep{RousseeuwL87,LopuhaaR91} \citep[see also][for an efficient (stochastic) algorithm]{RousseeuwD99}. Though such estimators are algorithmically NP-complete in exact version~\citep{Bernholt06}, this does alter (much) the general complexity of the (exact) algorithm, which is NP-complete already; the reader is referred to~\cite{JohnsonP78} or~\cite{DyckerhoffM16} for the exact algorithm for the halfspace depth. Approximate algorithms for both estimators have polynomial (linear in $n$) complexity, see, \textit{e.g.},~\cite{DyckerhoffMN21} for the halfspace depth.

%??? Discuss that anyway it's the same complexity, give references (paper on complexity of MCD, my computational taxonomy, ...).

%% TODO: Discussion on comensurability.

We illustrate the use of affine invariance on the following Example~\ref{exa:anomdet} of depth application for outlier detection.

\begin{example}[Outlier detection]\label{exa:anomdet}

To better illustrate practical importance of the affine-invariance property, we construct the following synthetic empirical example. To $100$ points drawn from the bivariate Gaussian distribution located at $(15, 55)^\top$ with variances $4$ and $16$ (for abscissa and ordinate, respectively) having correlation $0.625$, we add 5 outliers on Mahalanobis distance 5, and five clustered outliers at $(21, 35)^\top$ (stemming from normal distribution with variance $0.25$). Further, we calculate depths of all $110$ points and plot them in anomaly-wise ({\it i.e.}, first normal observations then outliers, in a supervised manner) and depth-increasing order, see Figure~\ref{fig:anomdet}. Please note, that the three figures below are regarded by a practitioner without the knowledge about anomalies ({\it i.e.} point's symbols and coloring). Figure~\ref{fig:anomdet} illustrates that, in the considered setting, while outliers are rather masked from the traditional halfspace depth, they can be distinguished (provided the threshold is properly chosen, which can be done, {\it e.g.}, using the receiver-operating characteristics) using the simple depth values with the robust halfspace depth, and can be literally separated using its affine-invariant version; a pattern that is recognizable even in absence of knowledge about outliers.

\begin{figure}[!b]
	\begin{center}
		\includegraphics[width=0.5\textwidth,trim=0cm 0.5cm 0 2cm,clip=true]{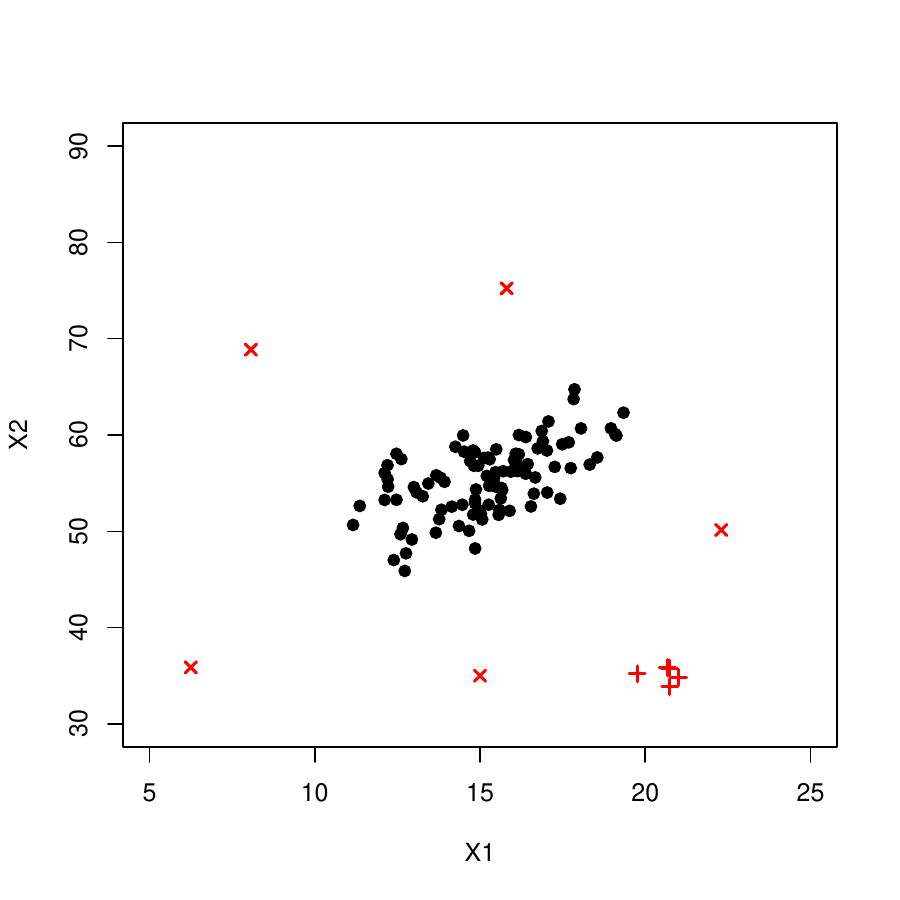}
	\end{center}
	\caption{Bi-variate normal sample ($100$ points) with two groups of outliers ($5$ distributed and $5$ clustered.}\label{fig:anomdetPoints}
\end{figure}
\begin{figure}
	\begin{tabular}{ccc}
		`Traditional' Tukey $D_0$ & Robust halfspace $D_{0.1}$ & AI robust halfspace $D_{0.1(\hat\Sigma_X)}$ \\
		\includegraphics[width=0.3\textwidth,trim=0 0.5cm 1cm 1cm,clip=true,page=1]{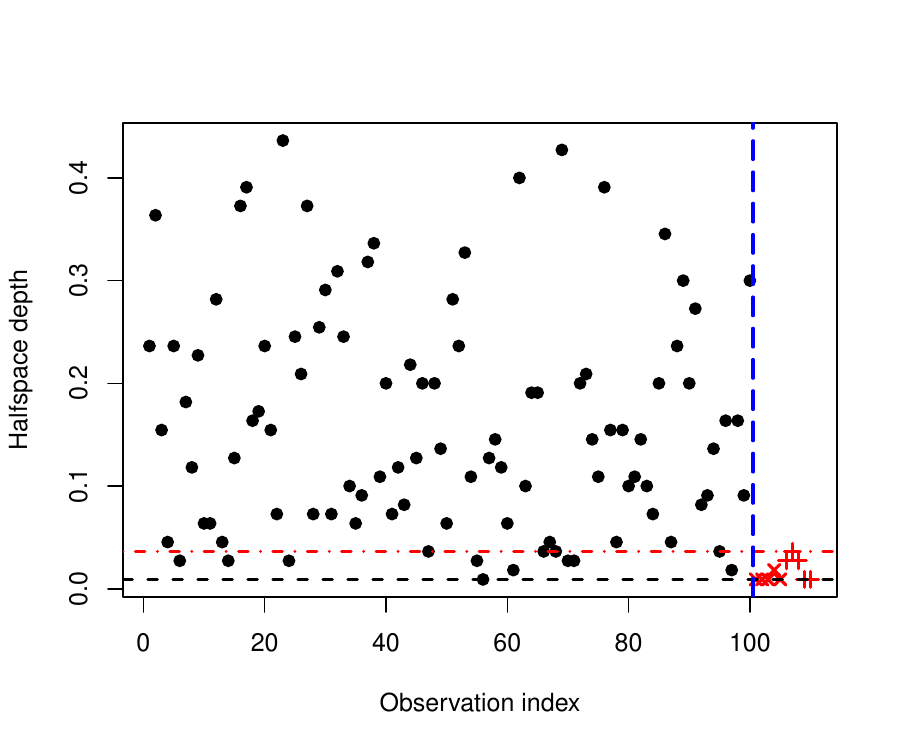} & \includegraphics[width=0.3\textwidth,trim=0 0.5cm 1cm 1cm,clip=true,page=4]{genAnom.pdf} & \includegraphics[width=0.3\textwidth,trim=0 0.5cm 1cm 1cm,clip=true,page=7]{genAnom.pdf} \\
		\includegraphics[width=0.3\textwidth,trim=0 0.5cm 1cm 1cm,clip=true,page=3]{genAnom.pdf} & \includegraphics[width=0.3\textwidth,trim=0 0.5cm 1cm 1cm,clip=true,page=6]{genAnom.pdf} & \includegraphics[width=0.3\textwidth,trim=0 0.5cm 1cm 1cm,clip=true,page=9]{genAnom.pdf}
	\end{tabular}
	\caption{Illustration of the depth-induced order for the standard (left), robust (middle), and affine-invariant robust (right) halfspace depth for the data from Figure~\ref{fig:anomdetPoints} (with outliers being assembled at the last $10$ indices). Three upper plots maintain natural order of the observation in the data set, while those below order them by depth value.}\label{fig:anomdet}
\end{figure}

\end{example}

\subsection{Algorithmic aspects}

All the above plots are generated using Algorithm~\ref{alg:sample} for $1000$ equally spaced directions.
In the following we show that robustified depth also allows for gradient-based optimization techniques, unlike the traditional depth. %It is expected that the number of local minima deminishes as $\delta$ grows.

Choose some smooth parameterization of the angles $\u=\u(\theta)$, where $\theta$ is $d-1$-dimensional.
Thus we have
\[\partial y_j/\partial\theta_i=\sum_k(z_k-x_{jk})\partial u_k/\partial \theta_i\]
at least for $\theta$ not on the boundary of the domain.
Consider the representation of the value $v(\bs u;\delta)$ in Corollary~\ref{cor:sample} and assume that $\u$ is such that 
\begin{equation}
y^{(-1)}\neq 0,\qquad y^{(k-1)}< y^{(k)}<y^{(k+1)},\qquad s_k\neq \delta n,
\end{equation}
where $1\leq k\leq m$ and $y^{(0)}=0,y^{(m+1)}=\infty$.
By continuity of $y_i$ in the angle $\u$ it must be that $p,k$ are constant for small perturbations of $\u$ and also that $y^{(1)},\ldots,y^{(k)}$ correspond to the same $\x_i$  up to the order of the first $k-1$.
Hence we readily find that 
\begin{equation}\label{eq:derivatives}
\frac{\partial v(\u;\delta)}{\partial \theta_i}=\frac{-1}{y^{(k)}}\left(\frac{1}{n}\sum_{j=1}^{k-1} \partial y^{(j)}/\partial \theta_i+\Big(v(\u;\delta)-\frac{n-m+k-1}{n}\Big)\partial y^{(k)}/\partial \theta_i\right),
\end{equation}
where $\partial y^{(j)}/\partial\theta_i$ is given above for the specific index.
It is noted that the above assumption is satisfied for Lebesgue-almost all directions in the $d-1$-dimensional simplex.

\begin{figure}[t!]
\centering
\includegraphics[width=0.48\textwidth]{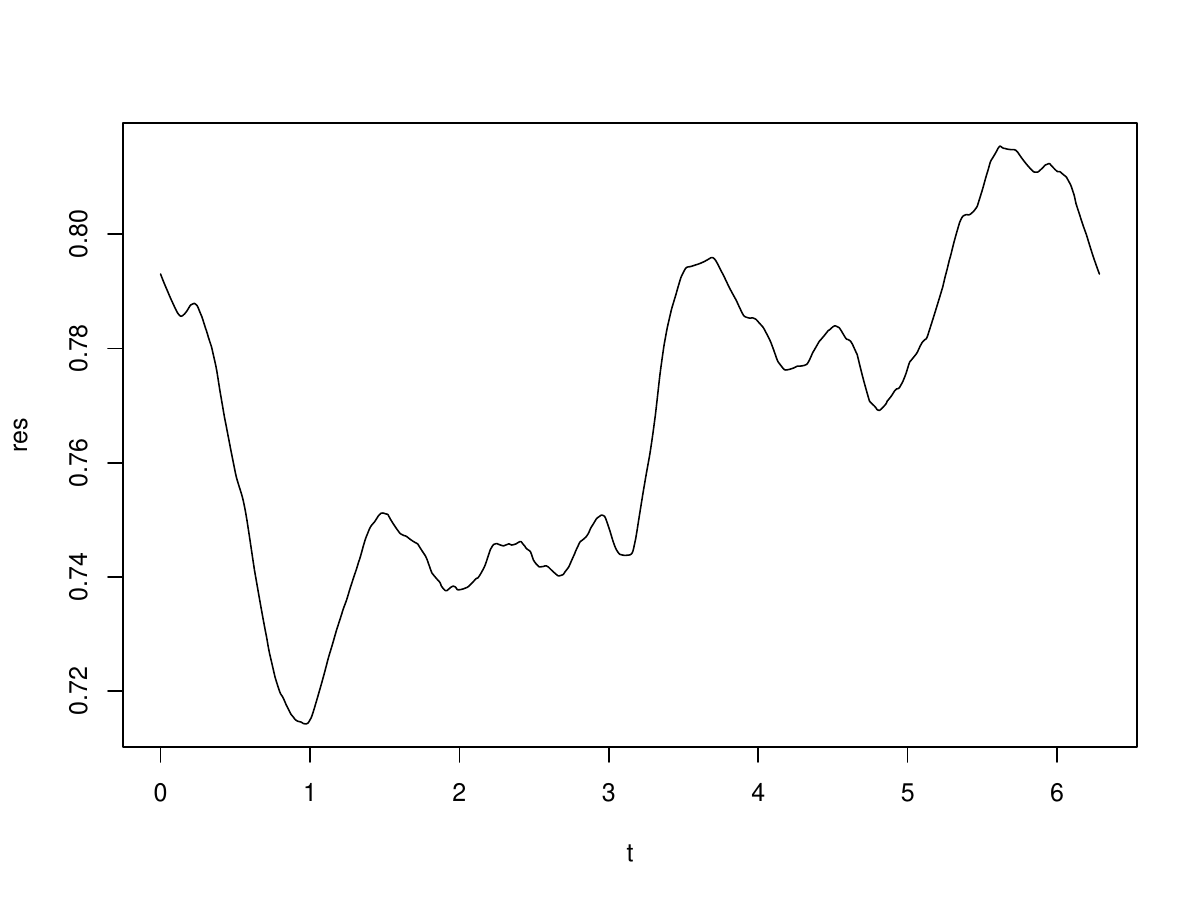}
\includegraphics[width=0.48\textwidth]{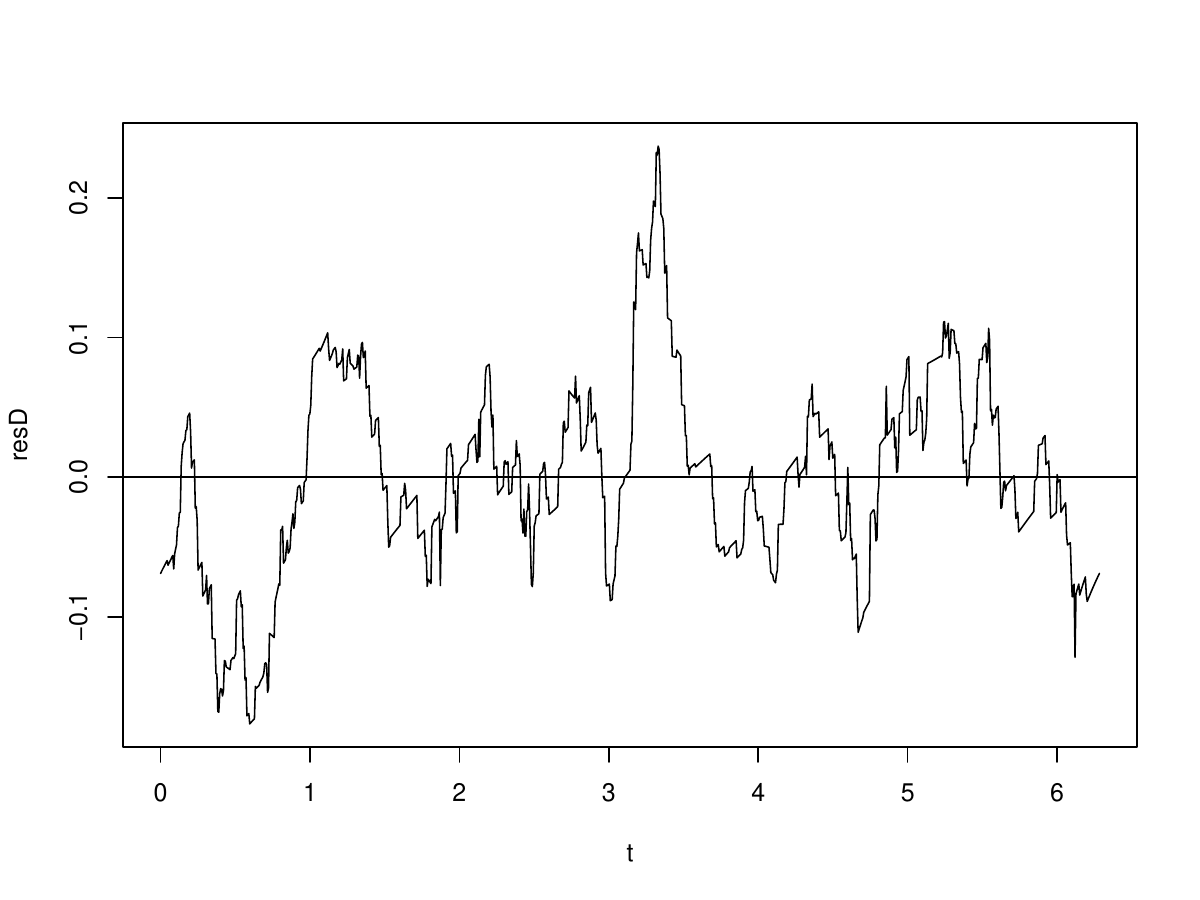}
\caption{The value function $v(\theta)$ and its derivative $v'(\theta)$ for $d=2,\delta=0.1$ and a sample with $n=100$ calculated at 1000 equally spaced $\theta\in[0,2\pi]$.}
\label{fig:derivative}
\end{figure}

In dimension $d=2$ we take $\u=(\cos\theta,\sin\theta),\theta\in[0,2\pi]$ and write $v(\theta)=v(\u(\theta);\delta)$.
An illustration of the value function and its derivative~\eqref{eq:derivatives} is given in Figure~\ref{fig:derivative} in the case of a bivariate standard normal sample with $n=100$ and $\delta=0.1$.
The minimum is achieved at $\theta=1.81$.

%Note also that for $u_d$ small we get numerical instability due to division by~$u_d$ in calculation of the derivative. One may either use some other parameterization of the angles, or to fix the index corresponding to the largest $u$
%There is a problem with the above parameterization for small $u_d$.
%What to use in $d\geq 3$?
%We may take 
%\[\u=(\cos\theta_1,\sin\theta_1\cos\theta_2,\ldots,\sin\theta_1\cdots\sin\theta_{d-1}),\]
%where $\theta_i\in[0,2\pi]$ and $\theta_{d-1}\in[-\pi/2,\pi/2]$?
%Another option is $\theta=(u_1,\ldots,u_{d-1},\pm)$, but that leads to numerical instability when $u_d$ is small (division by small number).
%One way around is to fix the index corresponding to the largest $|u_i|$.
%
%\pavlo{Having the derivative allows for application of gradient-based optimization techniques, which are known for fast speed of convergence and widely used available implementations. In this, it is superior to the traditional Tukey depth which calculation narrows down to optimization of a step-wise constant function where the derivative equals zero almost everywhere. The problem of non-(quasi-)convexity is inherited from the Tukey depth though. This is usually deal with by such techniques as, e.g., simulated annealing or multiple initialization, while in each uni-modal sub-problem gradient descent can be  employed directly.}
%

\section{On the choice of parameter $\delta$}\label{sec:parchoice}

First, let us consider recovering the value of $\delta$ that exactly compensates for the optimizer's curse, \textit{i.e.}, debiases the empirical depth. Closer inspection of \eqref{eq:h}--\eqref{eq:sol_alt} reveals that, clearly, in this case $\delta$ shall be dependent on the sample size $n$, but also on the true (\textit{i.e.}, population) depth $D_0(\z|\p)$. This oracle-based way though is not applicable directly in practice, will provide us with certain evidence about the behavior of optimal $\delta$. We implement the following procedure.

Let $\mathbb{R}\ni\lambda>0$ be such that:
\begin{equation}\label{eq:emp:lambda}
	\e\bigl( \inf_{\u\in\R^d:\|\u\|=1} \frac{1}{n}\sum_{i=1}^n\ind{\br{\u,\X_i-\z}\geq -\lambda} \bigr)  = D_0(\z|\p)\,.
\end{equation}

Reminding the reader that $Y=\br{\u,\z-X}$, the optimal value of $\delta$ can be estimated as
\begin{equation}\label{eq:emp:delta}
	\delta = \inf_{\u\in\R^d:\|\u\|=1} \e (Y\ind{Y\in(0,\lambda]})\,,
\end{equation}
where the expectation is independent from observations and thus can be expressed as an integral over projected probability measure.

\begin{figure}[t!]
\includegraphics[width=0.3\textwidth,trim=0 0.5cm 1cm 1.5cm,clip=true]{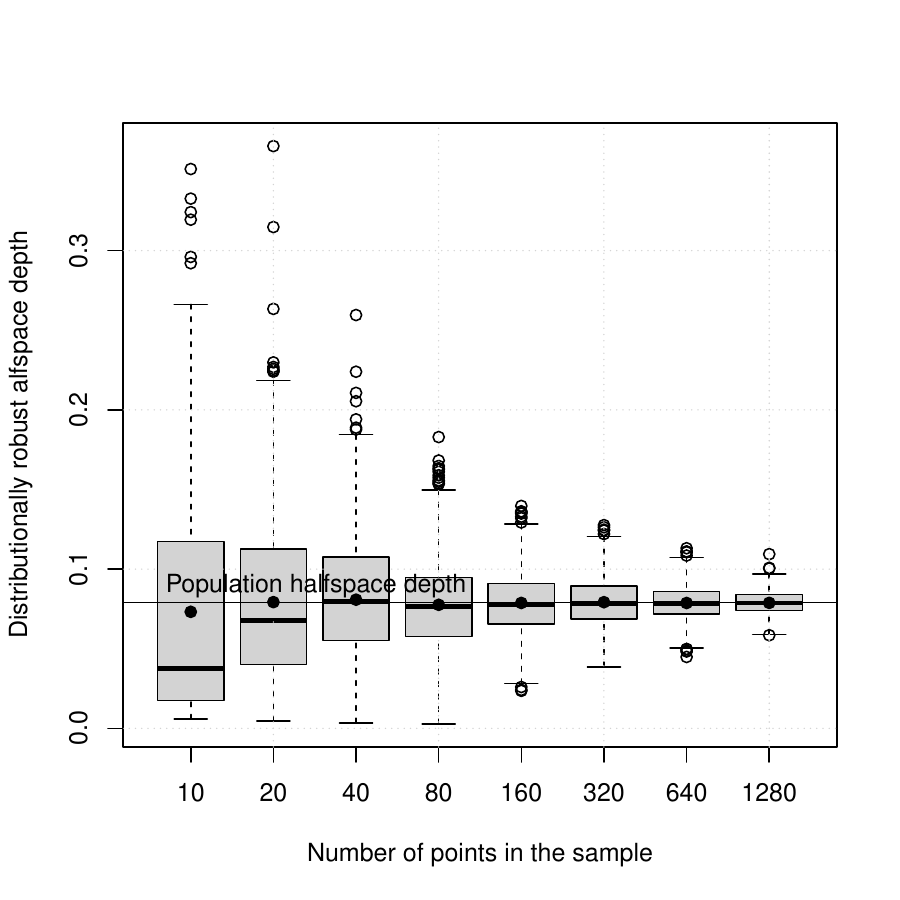} \quad \includegraphics[width=0.3\textwidth,trim=0 0.5cm 1cm 1.5cm,clip=true]{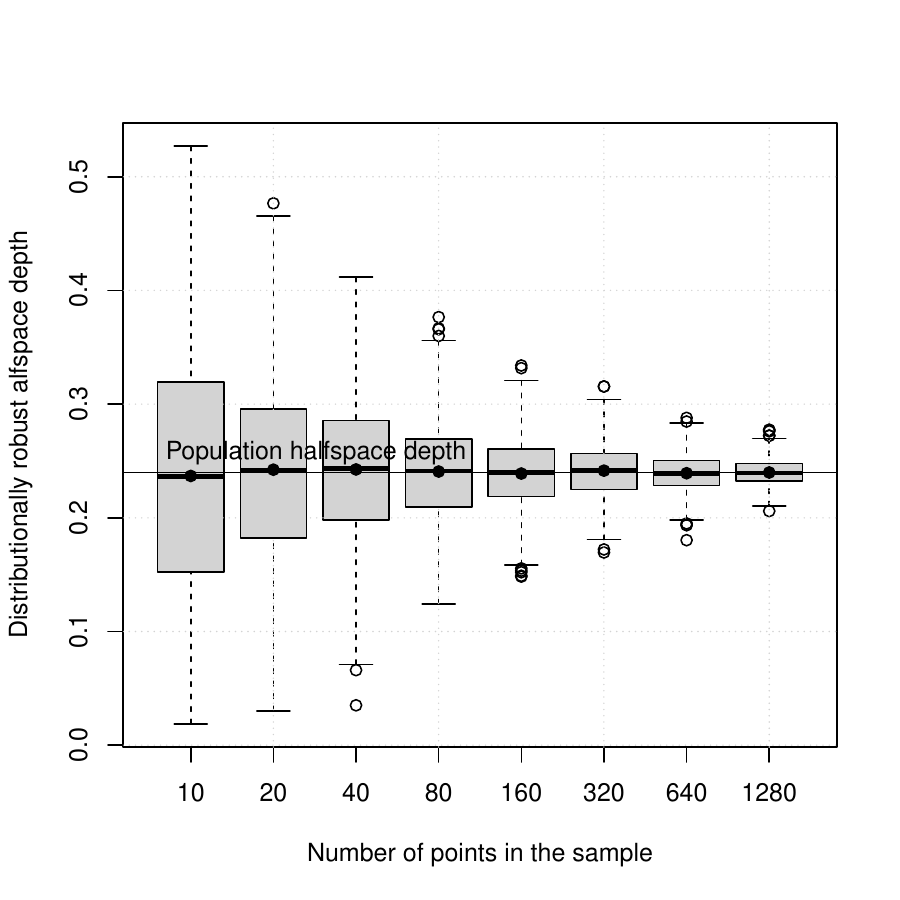} \quad \includegraphics[width=0.3\textwidth,trim=0 0.5cm 1cm 1.5cm,clip=true]{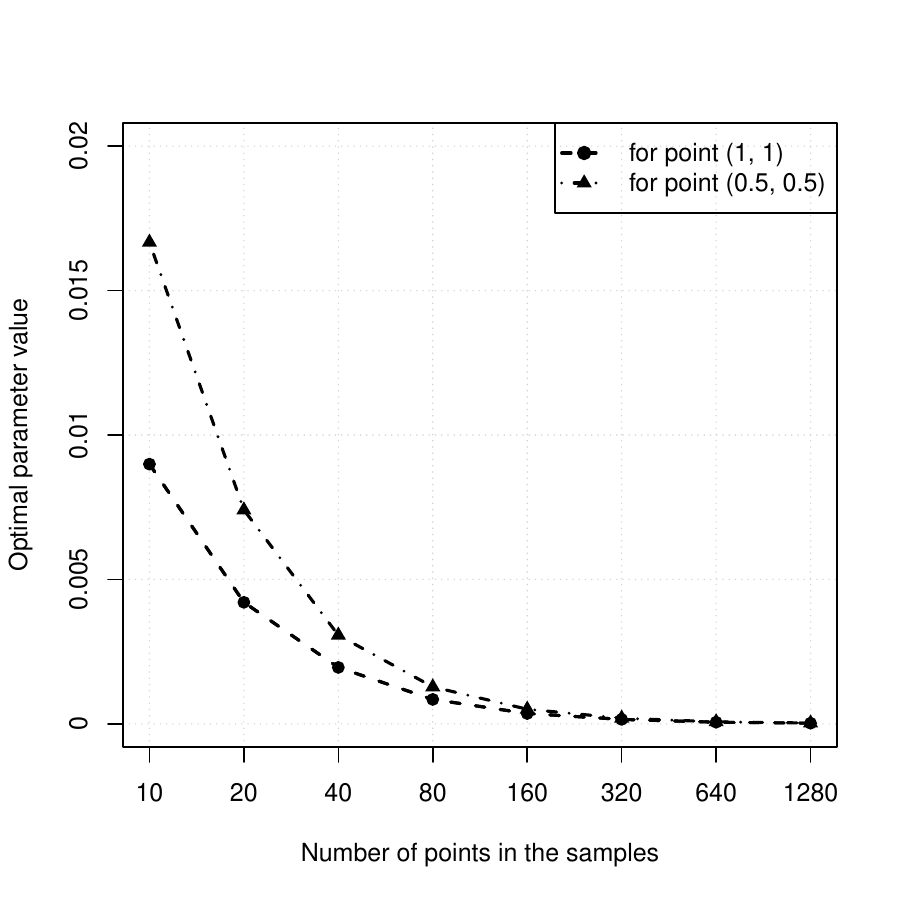}
\caption{Boxplots (over $1000$ random samples) of the empirical robustified hoalfspace depth of points $(1, 1)^\top$ (left) and $(0.5, 0.5)^\top$ (middle) w.r.t. a sample from bi-variate (correlated) normal distribution, for different sample sizes; chosen values of the parameter $\delta$ in each setting (right).}\label{fig:choice_delta_curse}
\end{figure}

While from equation \eqref{eq:emp:delta} one derives that $\delta$ depends on $\br{\u,X}$, equation \eqref{eq:emp:lambda} additionally testifies dependence on $\br{\u,\z}$. To nivellate at least dependence on $\br{\u,X}$, we consider bivariate standard normal distribution $X\sim \mathcal{N}\bigl((0, 0)^\top,\begin{pmatrix} 1 & 0 \\ 0 & 1 \end{pmatrix}\bigr)$, and apply the above described procedure to construct unbiased depth estimator for two points: $(1, 1)^\top$ and $(0.5, 0.5)^\top$; see Figure~\ref{fig:choice_delta_curse} indicating different values of $D_0(\z|\p_n)$ varying $n$.

%%First: is it possible to choose "optimal" (in a certain sense) value of $\delta$?

%%Let us ... to compensate for optimizer's curse.

%%Employing numerical evaluation of the quantities in equations (...) - (...), 

Nevertheless, one should not underestimate importance of such approach in practice. First, if the task is indeed to avoid exactly the optimizer's curse, a simulation study---based on available prior information and data set at hand---can provide a useful grid of values for $\delta$ to be later repeatedly used in practice. Furthermore, choosing reasonable $\delta$ using a generic procedure can provide mentioned above 
%%(see also the discussion in Section~\ref{sec:conclusions}) 
advantages independent of data; we illustrate this on the example of supervised classification (using class labels to determine the value of $\delta$) in Section~\ref{ssec:supclass}.

%%Think about cross-validation tuning of this parameter in DD-plot (even if diagonal) classification... Compare with other Tukey-depth versions.

%%\section{Outer level sets and the median region in the sample version}\label{sec:representations}
\section{Depth-trimmed regions}\label{sec:representations}

Importantly, in the sample version quite a bit more can be said about the shape of the median region and the outer level sets, that is, those corresponding to a small level~$\a$.
Characterization of such sets is addressed in this section.
\subsection{Outer level sets and depth decay}
%In the following we consider the empirical law $\p_n=\frac{1}{n}\sum\delta_{\{x_i\}}$.
The convex hull of the observations $\x_1,\ldots,\x_n$ is denoted by~$H$.
Let $d(\z,H)=\inf\{\norm{\z-\y}:\y\in H\}$ be the distance from $\z$ to the convex hull.
The next result shows that for $\a\leq 1/n$ the upper level set is $\delta/\a$-thickened convex hull~$H$, see Figure~\ref{fig:outer}. 
\begin{proposition}[Outer level sets]\label{prop:outer}
The following is true for any $\delta>0$ and $\z\in \R^d$:
\begin{itemize}
\item If $d(\z,H)\geq \delta n$ then $D_\delta(\z|\p_n)=\delta/d(\z,H)\leq 1/n$.
\item  For $\a\in(0,1/n]$ the level set is given by $\{\z:D_\delta(\z|\p_n)=\a\}=\{\z:d(\z,H)= \delta/\a\}$.
\item $D_\delta(\z|\p_n)\sim \delta/\norm{\z}$ as $\norm{\z}\to\infty$.
\end{itemize}
\end{proposition}
%%%\begin{proof}
%%%Let $\z^*$ be the unique point in $H$ such that $\norm{\z-\z^*}=d(\z,H)\geq \delta n$. 
%%%Consider the direction $\u=(\z-\z^*)/\norm{\z-\z^*}$ and observe that this direction maximizes $\min_i \br{\u,\z-\x_i}$.
%%%This can be easily seen by examining the vertices of the face containing~$\z^*$.
%%%For such $\u$ consider the solution to the inner problem as described in Corollary~\ref{cor:sample}.
%%%Note that $m=n$ and $y^{(1)}=d(\z,H)\geq \delta n$. Hence we get the solution~$\delta/d(\z,H)\leq 1/n$. Any other direction yields a larger value, since either $y^{(1)}$ is smaller or $m<n$. This proves the first statement.
%%%
%%%Since $\delta/\a\geq \delta n$ we see from the first statement that $\{\z\in\R^d:d(\z,H)= \delta/\a\}$ must have depth~$\a$.
%%%%The points $\z$ on the outside of this closed curve must have a smaller depth. 
%%%The second result now follows from the property P6 of the depth.
%%%
%%%For the third result we let $c$ be the maximal distance of the points in $H$ from the origin. Then $\norm{\z}-c\leq d(\z,H)\leq \norm{\z}+c$ showing that $d(\z,H)\sim \norm{\z}$, and the result follows easily from the first statement.
%%%\end{proof}

\begin{figure}[h!]
\centering
\begin{tikzpicture}[scale=0.6]
\path[draw=gray!10,fill=gray!10,line width=2.4cm,line cap=round,line join=round](0.5773439, 1.4530794)--(0.5980348, 2.6250300)--(1.5477594, 3.8998995)--(3.4675316, 3.5189625)--(4.6949597, 1.2503322)--(3.7649484, 0.7963902)--(0.5773439, 1.4530794);
\path[draw=gray!20,fill=gray!20,line width=1.4cm,line cap=round,line join=round](0.5773439, 1.4530794)--(0.5980348, 2.6250300)--(1.5477594, 3.8998995)--(3.4675316, 3.5189625)--(4.6949597, 1.2503322)--(3.7649484, 0.7963902)--(0.5773439, 1.4530794);

\draw (1.3,1.7) node{$H$};

\filldraw (4.6949597, 1.2503322) circle(1pt);
\filldraw (3.7649484, 0.7963902) circle(1pt);
\filldraw (3.4675316, 3.5189625) circle(1pt);
\filldraw (3.1438376, 2.6455422) circle(1pt);
\filldraw (1.5477594, 3.8998995) circle(1pt);
\filldraw (0.5773439, 1.4530794) circle(1pt);
\filldraw (2.7536316 ,3.0210816) circle(1pt);
\filldraw (0.5980348, 2.6250300) circle(1pt);
\filldraw (2.7720839, 1.0943229) circle(1pt);
\filldraw (2.8840230, 2.2034970) circle(1pt);

\draw (0.5773439, 1.4530794)--(0.5980348, 2.6250300)--(1.5477594, 3.8998995)--(3.4675316, 3.5189625)--(4.6949597, 1.2503322)--(3.7649484, 0.7963902)--(0.5773439, 1.4530794);
\end{tikzpicture}
\caption{Illustration of the outer upper level sets for $\a\leq 1/n$.}
\label{fig:outer}
\end{figure}

\subsection{The median region}

%The purpose of this section is to investigate the shape of the median region (those $z$ with $D_\delta(z)=1$), which has breakdown point at $1/2$. ...
Here our focus is on the median region $M_\delta$ defined in~\eqref{eq:median1} for the empirical law~$\p_n$ and $\delta>0$ so large that $\oa(\delta)=1$.
Define the class of sets
\begin{equation}\label{eq:separation}
\mathcal C=\Big\{I\subseteq \{1,\ldots,n\}:\{\x_i:i\in I\} \text{ and } \{\x_i:i\notin I\}\text{ can be separated by a hyperplane}\Big\},
\end{equation} 
which will be crucial for the following. Importantly, the cardinality of $\mathcal C$ is $O(n^d)$, which should be compared to $2^n$ possible subsets. 

\begin{proposition}\label{prop:median}
For $\delta>0$ such that $\oa(\delta)=1$ we have
\[M_\delta=\bigcap_{I\in\mathcal C,I\neq \emptyset} B_I=\bigcap_{I\neq \emptyset} B_I,\] where $B_I$ is the ball centered at $\sum_{i\in I} \x_i/|I|$ with radius $\delta n/|I|$. 
%$J$ runs over all non-empty subsets of $\{1,\ldots,n\}$ and $$. % runs over the non-empty subsets such that $\{\x_i:i\in I\}$ can be separated from the rest by a hyperplane.
Moreover, for $\delta$ large enough $M_\delta$  is a ball of radius $\delta$ centered at $\sum_i \x_i/n$.
\end{proposition}
We illustrate this result in Figure~\ref{fig:depth1}.
\begin{figure}[h!]
\centering
\begin{tikzpicture}[scale=0.8]
\def\r{3}
\filldraw (0,0) circle (2pt) node[above]{$\x_1$};
\filldraw (4,0) circle (2pt)node[below]{$\x_2$};
\filldraw (1,2) circle (2pt)node[left]{$\x_3$};
\def\a{(0,0) circle (\r)}
\def\b{(4,0) circle (\r)}
\def\c{(1,2) circle (\r)}
\def\ab{(2,0) circle (\r/2)}
\def\ac{(0.5,1) circle (\r/2)}
\def\bc{(2.5,1) circle (\r/2)}
\def\abc{(5/3,2/3) circle (\r/3)}
\draw[dashed]\a;
\draw\b;
\draw[dashed]\c;
\draw\ab;
\draw\ac;
\draw\bc;
\draw[dashed]\abc;
\begin{scope}
        \clip \b;
        \clip \ab;
        \clip \ac;
        \fill[red] \bc;
    \end{scope}
 \draw (2,0) circle (1pt);
\draw (2.5,1) circle (1pt);
\draw (0.5,1) circle (1pt);
\draw (5/3,2/3) circle (1pt);
\draw[<->](4.1,0)--(4+\r-0.1,0);
\draw (4+\r/2,0) node[above] {$3\delta$};
\end{tikzpicture}
\caption{The sample median region $M_\delta$ for $n=3$ and $\x_1,\x_2,\x_3$ depicted by black dots.}
\label{fig:depth1}
\end{figure}
Our proof relies on the following basic geometric optimization problem, concerning the maximal total Euclidean distance of given points to an arbitrary half-space.
Such a result  does not seem to appear in the standard books such as~\cite{boyd}, and so we provide a proof in the Supplementary Materials.
%%Appendix~\ref{sec:proofs_props}.

\begin{lemma}\label{lem:geometry}
 Consider non-zero vectors $\vv_i\in\R^d$ for $i=1,\ldots,n$. Then 
\begin{equation}\label{eq:reg1}{\rm argmax}_{\u\in\R^d:\norm{\u}=1}\left\{\sum_{i=1}^n\br{\u,\vv_i}^+\right\}\end{equation}
is non-empty and all of its elements have the form $\sum_{i\in I}\vv_i/\norm{\sum_{i\in I}\vv_i}$, where the set $I\neq\emptyset$ satisfies
\begin{equation}\label{eq:I}I=\{j: \sum_{i\in I}\br{\vv_i,\vv_j}>0\}.\end{equation}

The optimal value is given by 
$\max_I\{\norm{\sum_{i\in I}\vv_i}\}=\max_J\{\norm{\sum_{i\in J}\vv_i}\}$, where $I$ satisfies the above property and $J$ runs over all possible subsets of $\{1,\ldots,n\}$. 
\end{lemma}

This result will be used with $\bs v_i=\bs x_i-\z$, in which case the sets $I$ in~\eqref{eq:I} necessarily belong to~$\mathcal C$.
For example, the hyperplane passing through $\z+\epsilon\bs\nu$ with the normal $\bs\nu=\sum_{i\in I}\vv_i$ is such for a sufficiently small $\epsilon>0$.
Importantly, the class $\mathcal C$ does not depend on the chosen~$\z$. This allows us to prove Proposition~\ref{prop:median}, see the Supplementary Materials.

Let us discuss the description of $M_\delta$ in Proposition~\ref{prop:median}.
In practice, we may find the list of pairs $(\sum_{i\in I}\x_i,|I|)$ for all $I$ corresponding to the separable points in $O(n^d\log n)$ time.
Then for each $\z$ we can verify if it belongs to the median in $O(N)$ time, where $N=O(n^d)$ is the length of the above list. 
Some algorithms concerning the intersection of balls can be found in~\cite{balls_algs}.
For a fixed $\z$ it is sufficient to separate $\vv_i$ by a hyperplane passing through the origin, and the corresponding algorithm in dimension $d=2$ is given in the following.
%An easy way around this problem is to discretize a certain central region and to check for all points therein if they belong to~$M_\delta$. 
%We provide a fast and exact algorithm in the case of dimension $d=2$ in the following.

%It seems that exact checking if $\z\in M_\delta$ in $d=2$ can be done in $n\log n$ time (mostly arising from sorting the directions of $\vv_i$).

\subsection{An algorithm in $\R^2$}
Here we assume that $n\geq 1$ points $\vv_i=\x_i-\z\in\R^2$ are given. 
We present an algorithm computing the maximal Wasserstein cost of shifting the empirical distribution into an arbitrary half-space having $z$ on its boundary.
In other words, it finds the minimal $\delta>0$ resulting in $D_\delta(\z)=1$.
We note that the complexity of this algorithm is $O(n\log n)$ which comes from sorting the points $\vv_i$ according to their angles. 
We say that index $n$ is followed by index 1.
%\ji{This should be a very small modification of the algorithm for Tukey depth. Should we drop the following description and refer to such algorithms?}

\begin{algorithm}[H]\label{alg:median}
%\SetAlgoLined
 sort points $\vv_i$ so that their angle in $[0,2\pi)$ is non-decreasing\;
 $\delta=0$\;
 \For{$m=1,\ldots 2n$}{
  \eIf{$m=1$}{
    \tcp{initialization}
    $i:=1$\;
    find largest $j$ such that ${\rm det}[\vv_1,\vv_j]\geq 0$\;
    compute $\vv:=\vv_1+\cdots +\vv_j$\;
   }{
      take $k:=(j+1)\ind{j<n}+\ind{j=n}$ following $j$\;
      compute $d:={\rm det}[\vv_i,\vv_k]$\;
      \If{$d\geq 0$ and $i\neq k$}{ 
      \tcp{add $\vv_k$}
         $\vv:=\vv+\vv_k$\;
         $j:=k$\;
      }
      \If{$d\leq 0$}{ 
      \tcp{exclude $\vv_i$}
         \eIf{$i\neq j$}{
         $\vv:=\vv-\vv_i$\;
         $i:=(i+1)\ind{i<n}+\ind{i=n}$\;
         }{ 
				$\vv:=\vv_k$\;
				$i,j:=k$\;         
         }
      }
  }
  $\delta:=\max\{\delta,\norm{\vv}/n\}$\;
 }
 \KwResult{$\delta$}
 \caption{Compute the minimal $\delta$ resulting in $D_\delta(\z)=1$. Assumption: $d=2$.}
\end{algorithm}
\medskip

The basic observation behind this algorithm is that the $2n$ proposed subsets necessarily contain all subsets $I$ such that $\{\vv_i,i\in I\}$ can be separated from the rest by a line passing through~$\bs 0$; additionally, there are some non-separable subsets when some points $\vv_i$ lay on the same line passing through the origin.
These, in turn, include all the subsets $I$ satisfying the property~\eqref{eq:I} as mentioned above. 
In order to understand the number of the latter subsets we perform the following experiment:
generate $100$ points from a bivariate normal with correlation a) $0$ and b) 0.7, calculate the number of such sets and replicate 1000 times.
We find that the average number of sets $I$ satisfying~\eqref{eq:I} is approximately $37$ in a) and $9$ in b) out of $2n=200$.

Next, we use Algorithm~\ref{alg:median} to construct the median regions for various $\delta$ and a certain empirical distribution, see Figure~\ref{fig:depth1}.
We discretize a certain central region and compute the minimal $\delta$ resulting in $D_\delta(\z|\p_n)=1$ for all the given points~$\z$.
Then we select those points exceeding a given threshold.
\begin{figure}
\centering
\includegraphics[width=0.475\textwidth,trim=0 0 0 1.5cm,clip=true,page=1]{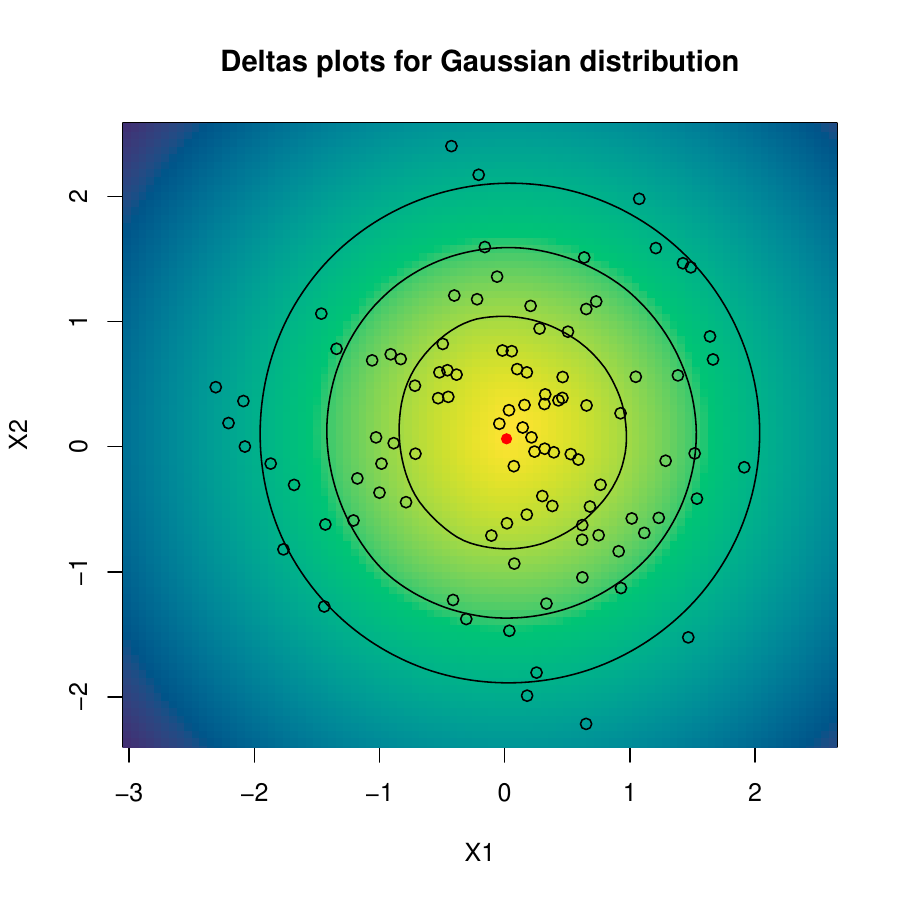} \quad \includegraphics[width=0.475\textwidth,trim=0 0 0 1.5cm,clip=true,page=2]{picMedian-n100.pdf}
%%\caption{100 samples from a bivariate normal with correlation 0.7 and the median regions corresponding to $\delta=0.6,0.8,1,2$ (smaller deltas correspond to smaller regions). 
%%The black dot at $(0,0.09)$ is the median at $\delta=0.483$.
%%Interpretation: for any line passing through $\z$ the data set can be moved to any side of the line with less than $\delta$ cost. \ji{Pasha, can you generate this and maybe without a big circle?}
%%}
\caption{Minimal $\delta$s delivering (maximal) depth $1$, depicted in color (darker color means higher values). The three depicted contours correspond (from inside to outside) to $\delta=1, 1.5, 2$. Left: $100$ bivariate observations drawn from a standard normal distribution. Right: $100$ bivariate observations drawn from a normal distribution with correlation $0.5$. The filled red dot is the median for the smallest $\delta$ $=$ $0.398$ (left) and $0.497$ (right) resulting in depth~$1$.}
\label{fig:depth1}
\end{figure}

\bigskip
Let us briefly discuss a possible application.
From the basic theory we see that 
\[D_\delta(\z|\p)=1\qquad\text{iff}\qquad \forall \u\neq\bs 0\quad\exists \p^*:d_W(\p^*,\p)\leq \delta,\, \p^*(\br{\u,\X-\z}>0)=0,\]
where we have used $1-\p^*(\br{\u,\X-\z}\geq 0)=\p^*(\br{-\u,\X-\z}>0)$ and then swapped $-\u$ for $\u$.
Suppose that some law $\p$ is given together with the bound $d_W(\p,\p_0)\leq\delta$ for a known $\delta>0$.
We would like to determine the points $\z$ such that for any direction $\u$ it is possible that $\br{\u,\X-\z}>0$ has no mass under~$\p_0$.
That is, we cannot choose a direction guaranteeing to see some mass.
Given the available information such points correspond to $D_\delta(\z|\p)=1$, i.e., to the median region of $\p$ for the ambiguity~$\delta$.

\section{Relation to other depth notions}\label{sec:relations}

As a consequence of its popularity, many modifications of the halfspace depth are available in the literature; we elaborate on relation to those right below.

Introduced by \cite{RamsayDL19} integrated rank-weighted halfspace depth allows for depth-sensitivity beyond the empirical convex hull, a property useful in applications, \textit{e.g.}, anomaly detection. \cite{ClemenconMS23} suggest its affine-invariant version along with finite-sample analysis. Since (affine-invariant) integrated rank-weighted halfspace depth relies on integration and not minimization, proposed here mechanism is not readily applicable.

\cite{chernozhukov} and~\cite{HallinDBCAM21} proposed Monge-Kantorovich depth based on the optimal transport theory, which yields spatial ranks and signs, and---different to other existing depth notions---directly indexes depth contours by their probability mass. Their basic construction consists of transforming $\p$ into the spherical uniform distribution and applying the classical halfspace depth (based on this uniform distribution) to the transformed points. 
We employ optimal transport theory in a different way. Our new distribution is not some fixed well-chosen reference distribution---it is the worst-case distribution for the given direction $\u$ and ambiguity radius~$\delta>0$. 

\cite{hlubinka} (see also~\cite{KotikH17}) proposed generalization of the halfspace depth using a weight function. That framework allows for the following formulation
\[\inf_{\u}\frac{\e w(\br{\u,\X-\z})}{\e w(-\br{\u,\X-\z})},\]
where $w(y)=0$ for $y<0$. In our case the problem~\eqref{eq:minimax} can be rewritten as $\inf_{\norm{\u}=1}$ $f\big(\mathcal L(\br{\u,\X-\z})\big)$, where the function $f$ is applied to the law of the projection on~$\u$, see~\eqref{eq:sol} below, and this seems to be the closest representation.

\cite{nagy_illumination} introduced the illumination depth based on a fixed upper level set $R$ corresponding to he classical halfspace depth.
More precisely, the illumination depth of a point $\z$ outside of $R$ is defined as 
the volume of the convex hull of $\{\z\}\cup R$ divided by the volume of~$R$. 
Note, however, that this depth depends only on the shape of $R$ and the relative position of~$\z$. 
For the empirical law $\p_n$ this leads to a non-trivial ordering beyond the convex hull of the set of observations.
The similarity to our depth comes from the fact that both depths are asymptotically inversely proportional to the distance from $\z$ to the convex hull of observations as $\z$ escapes to infinity.

\cite{einmahl2015} use extreme value theory to adjust the empirical halfspace depth when far away from the center. 
Assuming multivariate regular variation of $\p$, they extrapolate from moderately remote regions into the remote regions with few or no observations.
In the case of robustified halfspace depth the extrapolation for the empirical law is simplistic---the depth decays as $\delta/\norm{\z}$ as will be shown in the following.

\begin{figure}[t!]
\includegraphics[width=0.47\textwidth,page=1]{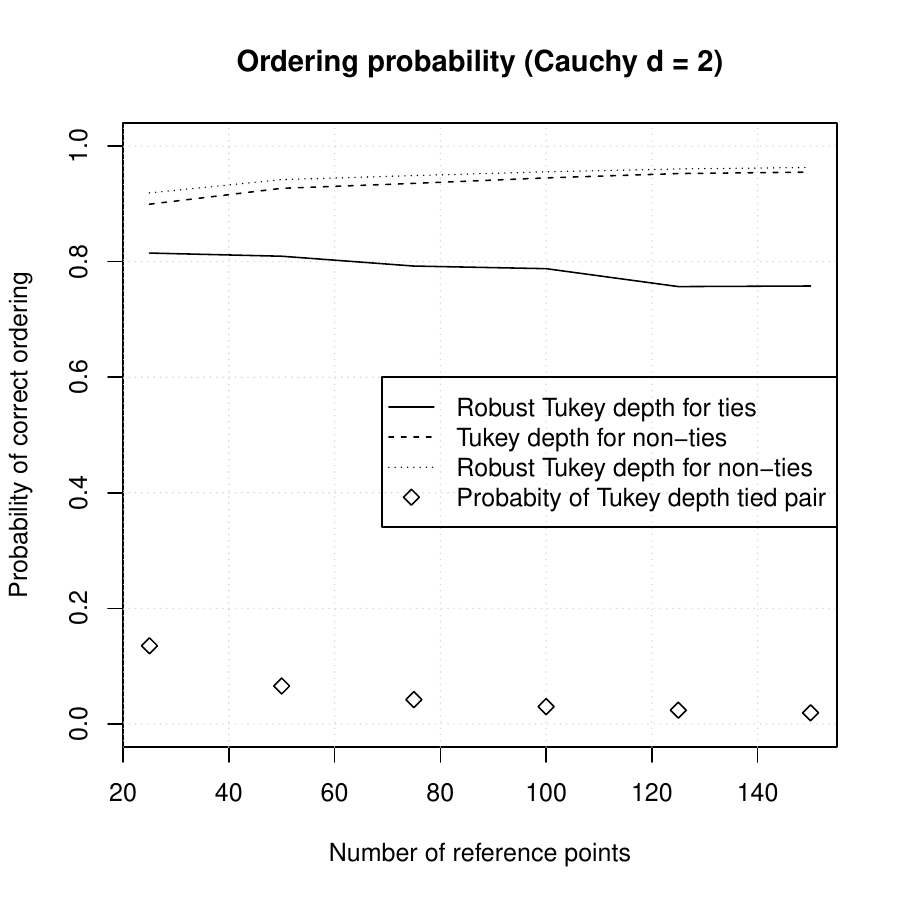}\quad
\includegraphics[width=0.47\textwidth,page=2]{pic-probOrder.pdf}
\caption{Conditional probability of the correct ordering of a pair of points stemming from an elliptical Cauchy distribution for dimension $d=2$ (left) and $d=10$ (right). The three lines correspond to (a) robust Tukey depth (dotted)  and (b) Tukey depth (dashed) where the latter is not tied, and (c) the robust Tukey depth where the Tukey depth is tied (solid). The rhombuses show a probability of drawing a pair tied by the empirical Tukey depth.}
\label{fig:depthorder}
\end{figure}

\section{Numerics}\label{sec:numerics}

\subsection{Ordering experiment}

Here we illustrate the orderings provided by the robust Tukey depth and by the traditional Tukey depth in the case of a relatively small sample. We use an elliptical Cauchy distribution in dimensions $d=2$ and $d=10$, where the scatter matrix is given by $(2^{-|i-j|})_{i,j=1}^d$.
Population depth is used to define a reference ordering, which in this case coincides with the density-based ordering~\citep{LiuS93} and so it is found  using the Mahalanobis distance from the center.

Firstly, we focus on the probability that a pair of points (sampled independently from the underlying Cauchy distribution) is correctly ordered by (a) the empirical robust depth with $\delta=0.1$ and (b) the empirical traditional Tukey depth.
For a fair comparison, this pair of points, in fact, is drawn conditional on having distinct empirical Tukey depth (not tied). We simply iterate until getting such a pair. 
Secondly, we estimate the probability of correct ordering by (c) the empirical robust depth conditional on a pair of points having the same empirical Tukey depth (tied). For each sample size $n$ we generate $1999$ empirical laws and  one pair of points with the required conditional law. 

As we can observe from Figure~\ref{fig:depthorder}, not only the robust Tukey depth provides ranking very similar to the traditional one, but also allows for (a rather consistent and not available) ordering of the Tukey depth's ties.

%\begin{figure}
%\includegraphics[width=0.48\textwidth]{uniform_true.pdf}
%\includegraphics[width=0.48\textwidth]{uniform_25_tukey.pdf}\\
%\includegraphics[width=0.48\textwidth]{uniform_25_002.pdf}
%\includegraphics[width=0.48\textwidth]{uniform_25_01.pdf}
%\caption{Uniform distribution on the unit square and a sample of 25 points. True Tukey depth, empirical Tukey depth, robustified depth with $\delta=0.02$ and with $\delta=0.1$}
%\end{figure}

\begin{figure}[!t]
	\begin{tabular}{ccc}
		\quad Normal-shift alternative & \quad Adaptive choice of $\delta$ & \quad 200 + 200 obs. (100 rep.) \\
		\includegraphics[width=0.28\textwidth,trim=0 0.5cm 1cm 1cm,clip=true]{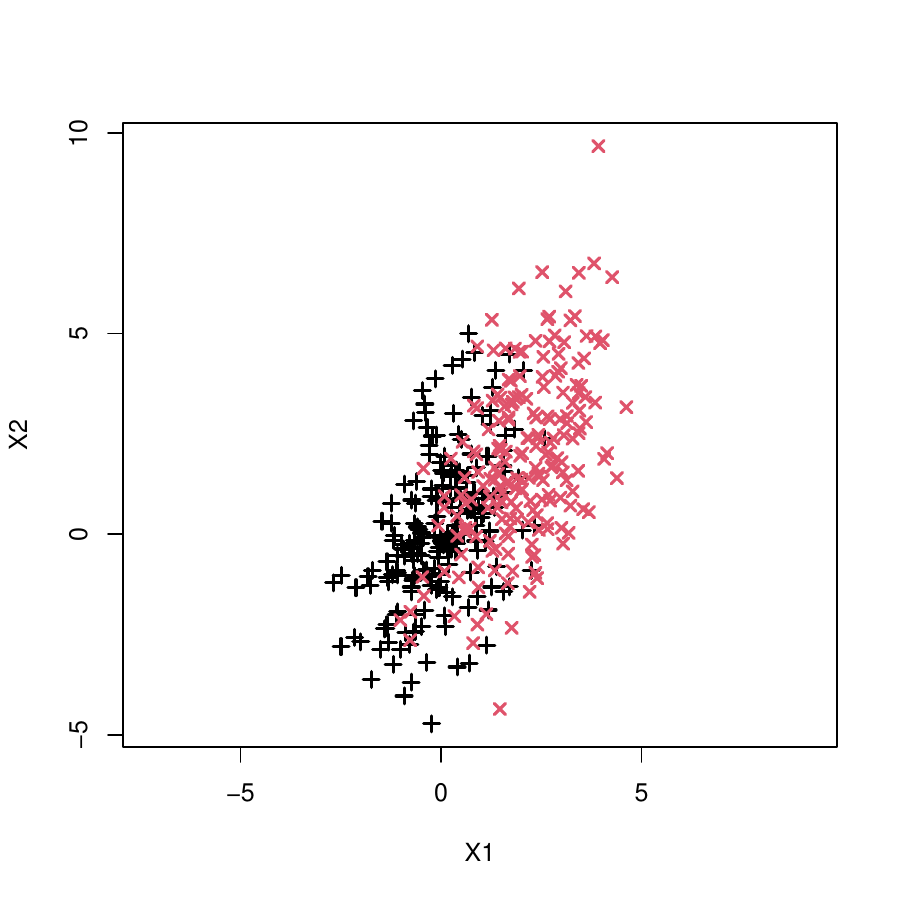} & \includegraphics[width=0.28\textwidth,trim=0 0.5cm 1cm 1cm,clip=true]{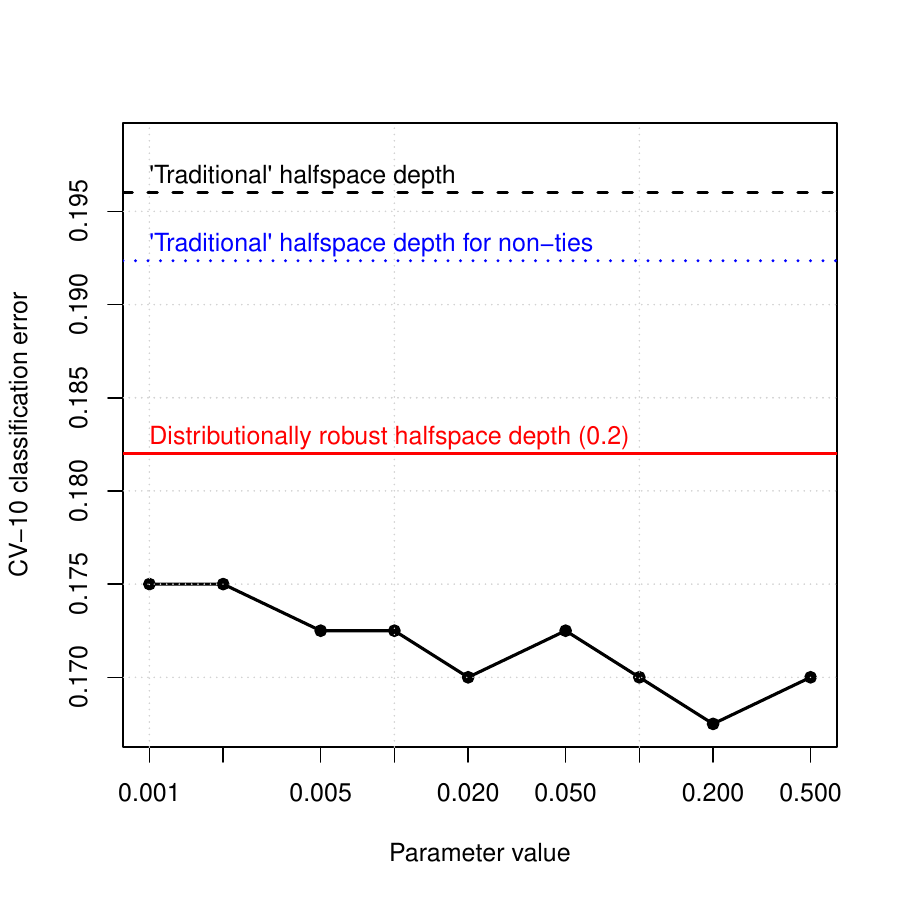} & \includegraphics[width=0.28\textwidth,trim=0 0.5cm 1cm 1cm,clip=true]{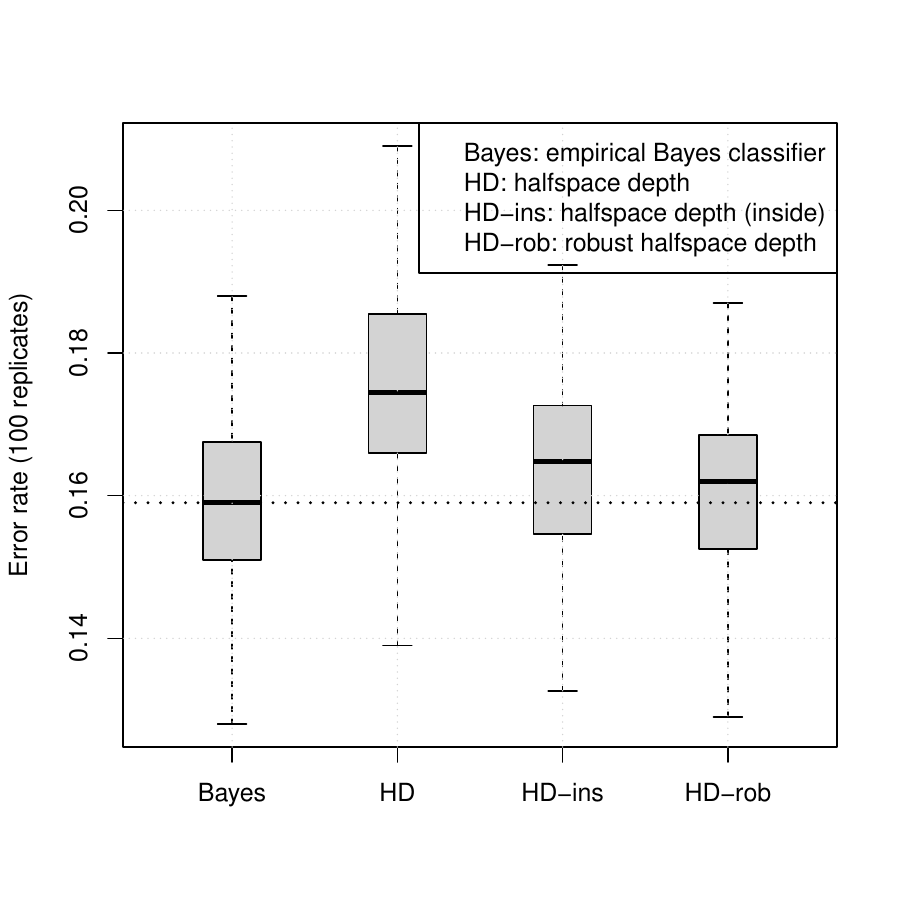} \\
		100 + 100 obs. (100 rep.) & \quad 50 + 50 obs. (100 rep.) & \quad 20 + 20 obs. (100 rep.) \\
		\quad \includegraphics[width=0.28\textwidth,trim=0 0.5cm 1cm 1cm,clip=true]{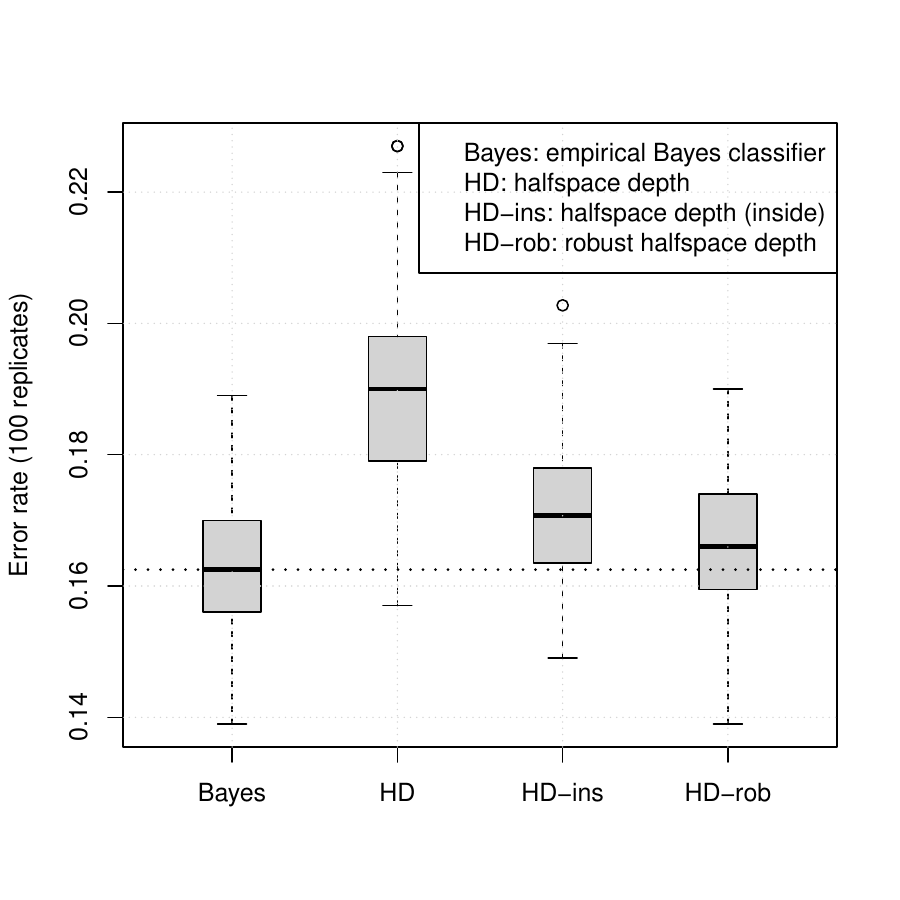} & \includegraphics[width=0.28\textwidth,trim=0 0.5cm 1cm 1cm,clip=true]{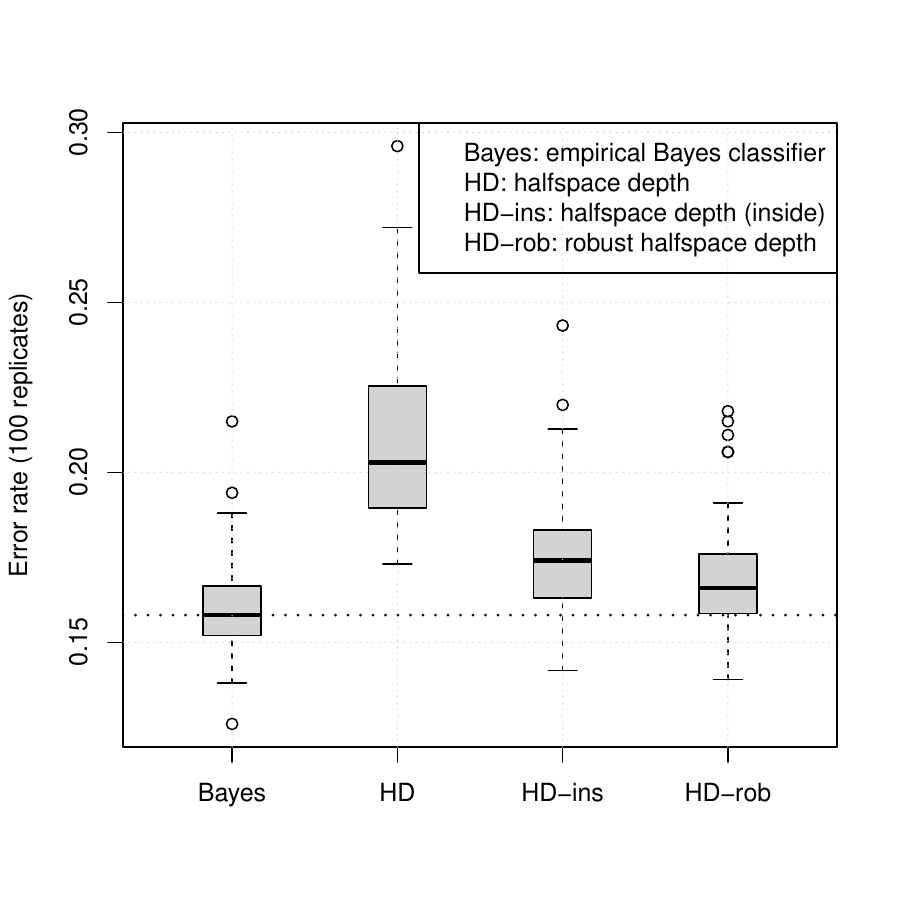} & \includegraphics[width=0.28\textwidth,trim=0 0.5cm 1cm 1cm,clip=true]{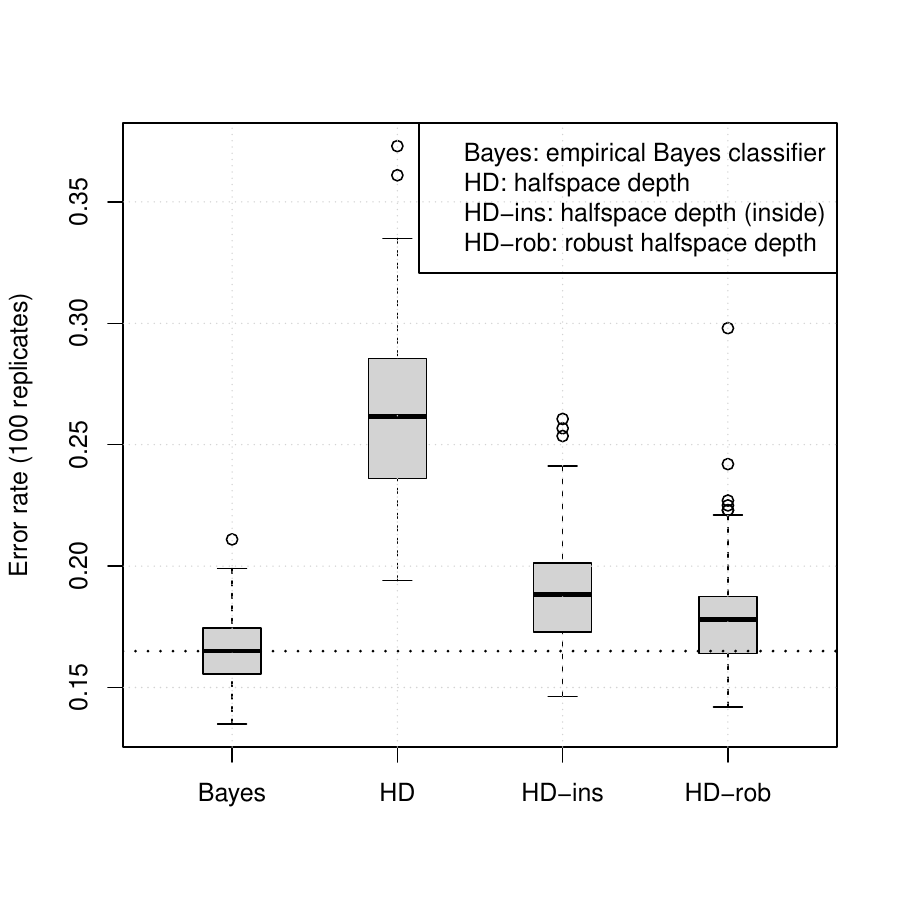} \\
	\end{tabular}
	\caption{Tuning parameter $\delta$ for supervised classification by the maximum depth classifier for the normal-shift alternative. Top, left: training sample having $200+200$ observations stemming from bi-variate correlated normal distributions differing in location shift only. Top, middle: For a training sample, $10$-fold cross-validation classification error for different choices of $\delta$ (as well evaluated on $1000$ test data error for traditional halfspace depth (also excluding ties) and robust halfspace depth with optimal choice of $\delta=0.2$). The $4$ remaining plots indicate bosxplots of classification errors on test sample over $100$ draws for varying size of training set.}\label{fig:maxdepth}
\end{figure}

\subsection{Supervised classification}\label{ssec:supclass}

As announced in Section~\ref{sec:parchoice}, below we suggest an example of using available in data feedback to tune parameter $\delta$, for the supervised classification task. For simplicity of exposition, we stick to a simple example of the normal-shift alternative, which is optimally handled by the maximum depth classifier~\cite{GhoshC05} (see also~\cite{LangeMM14} for a more elaborate procedure). More precisely, with two classes having equal prior probabilities represented by two elliptically symmetric (here normal) random vectors  in $\mathbb{R}^d$ (here $d=2$) differing in location only, say $X_0$ and $X_1$, the optimal classifier for any $\z\in\mathbb{R}^d$ is given by:
\begin{equation}\label{eq:maxdepth}
	g(\z|X_0,X_1) = \ind{D(\z|X_1) > D(\z|X_0)}\,
\end{equation}
with ties broken at random. On the $DD$-plot~\cite{LiCAL12} with points $\bigl(D(\z|X_0), D(\z|X_1)\bigr)^T$, rule~\eqref{eq:maxdepth} corresponds to the diagonal line.

Using classes' labels as feedback, \textit{e.g.}, by means of cross-validation, the value of parameter $\delta$ can be chosen in an adaptive way. For two normally distributed classes, in Figure~\ref{fig:maxdepth} we illustrate this process by plotting the training data set, cross-validation error, as well as simulation-estimated (over $100$ repetitions, using a large separate test sample) errors of this adaptive procedure along with those for traditional halfspace depth (both breaking ties at random or ignoring them) and empirical Bayes classfier (=linear discriminant analysis).

%%%%%%%%%%%%%%%%%%%%%%%%%%%%%%%%%%%%%%%%%%%%%%
%% Single Appendix:                         %%
%%%%%%%%%%%%%%%%%%%%%%%%%%%%%%%%%%%%%%%%%%%%%%
%\begin{appendix}
%\section*{???}%% if no title is needed, leave empty \section*{}.
%\end{appendix}
%%%%%%%%%%%%%%%%%%%%%%%%%%%%%%%%%%%%%%%%%%%%%%
%% Multiple Appendixes:                     %%
%%%%%%%%%%%%%%%%%%%%%%%%%%%%%%%%%%%%%%%%%%%%%%
%\begin{appendix}
%\section{???}
%
%\section{???}
%
%\end{appendix}

%%%%%%%%%%%%%%%%%%%%%%%%%%%%%%%%%%%%%%%%%%%%%%
%% Support information, if any,             %%
%% should be provided in the                %%
%% Acknowledgements section.                %%
%%%%%%%%%%%%%%%%%%%%%%%%%%%%%%%%%%%%%%%%%%%%%%
%\begin{acks}[Acknowledgments]
% The authors would like to thank ...
%\end{acks}
%%%%%%%%%%%%%%%%%%%%%%%%%%%%%%%%%%%%%%%%%%%%%%
%% Funding information, if any,             %%
%% should be provided in the                %%
%% funding section.                         %%
%%%%%%%%%%%%%%%%%%%%%%%%%%%%%%%%%%%%%%%%%%%%%%
\section*{Acknowledgements}
J.\ Ivanovs gratefully acknowledges financial support of Sapere Aude Starting Grant 8049-00021B ``Distributional Robustness in Assessment of Extreme Risk''. P.\ Mozharovskyi gratefully acknowledges the support of the Young Researcher Grant of the French National Agency for Research (ANR JCJC 2021) in category Artificial Intelligence registered under the number ANR-21-CE23-0029-01.

%%\bibliographystyle{Chicago}
%%\bibliography{model_risk}

\begin{thebibliography}{}

\bibitem[\protect\citeauthoryear{Aurenhammer}{Aurenhammer}{1988}]{balls_algs}
Aurenhammer, F. (1988).
\newblock Improved algorithms for discs and balls using power diagrams.
\newblock {\em Journal of Algorithms\/}~{\em 9\/}(2), 151--161.

\bibitem[\protect\citeauthoryear{Bernholt}{Bernholt}{2006}]{Bernholt06}
Bernholt, T. (2006).
\newblock Robust estimators are hard to compute.
\newblock Technical report, University of Dortmund.

\bibitem[\protect\citeauthoryear{Bingham, Goldie, and Teugels}{Bingham
  et~al.}{1989}]{BGT}
Bingham, N.~H., C.~M. Goldie, and J.~L. Teugels (1989).
\newblock {\em Regular Variation}, Volume~27 of {\em Encyclopedia of
  Mathematics and its Applications}.
\newblock Cambridge University Press, Cambridge.

\bibitem[\protect\citeauthoryear{Blanchet, Kang, and Murthy}{Blanchet
  et~al.}{2019}]{blanchet2016profile}
Blanchet, J., Y.~Kang, and K.~Murthy (2019).
\newblock Robust {W}asserstein profile inference and applications to machine
  learning.
\newblock {\em Journal of Applied Probability\/}~{\em 56\/}(3), 830--857.

\bibitem[\protect\citeauthoryear{Blanchet and Murthy}{Blanchet and
  Murthy}{2019}]{blanchet_distributional}
Blanchet, J. and K.~Murthy (2019).
\newblock Quantifying distributional model risk via optimal transport.
\newblock {\em Mathematics of Operations Research\/}~{\em 44\/}(2), 565--600.

\bibitem[\protect\citeauthoryear{Bousquet and Elisseeff}{Bousquet and
  Elisseeff}{2002}]{bousquet2002stability}
Bousquet, O. and A.~Elisseeff (2002).
\newblock Stability and generalization.
\newblock {\em Journal of Machine Learning Research\/}~{\em 2}, 499--526.

\bibitem[\protect\citeauthoryear{Boyd and Vandenberghe}{Boyd and
  Vandenberghe}{2004}]{boyd}
Boyd, S. and L.~Vandenberghe (2004).
\newblock {\em Convex Optimization}.
\newblock Cambridge University Press, Cambridge.

\bibitem[\protect\citeauthoryear{Chakraborty and Chaudhuri}{Chakraborty and
  Chaudhuri}{1996}]{ChakrabortyC96}
Chakraborty, B. and P.~Chaudhuri (1996).
\newblock On a transformation and re-transformation technique for constructing
  an affine equivariant multivariate median.
\newblock {\em Proceedings of the American Mathematical Society\/}~{\em 124},
  2539--2547.

\bibitem[\protect\citeauthoryear{Chernozhukov, Galichon, Hallin, and
  Henry}{Chernozhukov et~al.}{2017}]{chernozhukov}
Chernozhukov, V., A.~Galichon, M.~Hallin, and M.~Henry (2017).
\newblock Monge-{K}antorovich depth, quantiles, ranks and signs.
\newblock {\em The Annals of Statistics\/}~{\em 45\/}(1), 223--256.

\bibitem[\protect\citeauthoryear{Cl{\'e}men{\c{c}}on, Mozharovskyi, and
  Staerman}{Cl{\'e}men{\c{c}}on et~al.}{2023}]{ClemenconMS23}
Cl{\'e}men{\c{c}}on, S., P.~Mozharovskyi, and G.~Staerman (2023).
\newblock {Affine invariant integrated rank-weighted statistical depth:
  properties and finite sample analysis}.
\newblock {\em Electronic Journal of Statistics\/}~{\em 17\/}(2), 3854--3892.

\bibitem[\protect\citeauthoryear{Cuesta-Albertos and
  Nieto-Reyes}{Cuesta-Albertos and Nieto-Reyes}{2008}]{CuestaAlbertosNR08}
Cuesta-Albertos, J. and A.~Nieto-Reyes (2008).
\newblock The random {Tukey} depth.
\newblock {\em Computational Statistics and Data Analysis\/}~{\em 52},
  4974--4988.

\bibitem[\protect\citeauthoryear{Cuturi and Doucet}{Cuturi and
  Doucet}{2014}]{ML_earth3}
Cuturi, M. and A.~Doucet (2014, 22--24 Jun).
\newblock Fast computation of wasserstein barycenters.
\newblock In E.~P. Xing and T.~Jebara (Eds.), {\em Proceedings of the 31st
  International Conference on Machine Learning}, Volume~32 of {\em Proceedings
  of Machine Learning Research}, Bejing, China, pp.\  685--693. PMLR.

\bibitem[\protect\citeauthoryear{Donoho and Gasko}{Donoho and
  Gasko}{1992}]{donoho_gasko}
Donoho, D.~L. and M.~Gasko (1992).
\newblock Breakdown properties of location estimates based on halfspace depth
  and projected outlyingness.
\newblock {\em The Annals of Statistics\/}~{\em 20\/}(4), 1803--1827.

\bibitem[\protect\citeauthoryear{Dyckerhoff}{Dyckerhoff}{2004}]{Dyckerhoff04}
Dyckerhoff, R. (2004).
\newblock Data depths satisfying the projection property.
\newblock {\em Allgemeines Statistisches Archiv\/}~{\em 88}, 163--190.

\bibitem[\protect\citeauthoryear{Dyckerhoff and Mozharovskyi}{Dyckerhoff and
  Mozharovskyi}{2016}]{DyckerhoffM16}
Dyckerhoff, R. and P.~Mozharovskyi (2016).
\newblock Exact computation of the halfspace depth.
\newblock {\em Computational Statistics and Data Analysis\/}~{\em 98}, 19--30.

\bibitem[\protect\citeauthoryear{Dyckerhoff, Mozharovskyi, and Nagy}{Dyckerhoff
  et~al.}{2021}]{DyckerhoffMN21}
Dyckerhoff, R., P.~Mozharovskyi, and S.~Nagy (2021).
\newblock Approximate computation of projection depths.
\newblock {\em Computational Statistics and Data Analysis\/}~{\em 157}, 107166.

\bibitem[\protect\citeauthoryear{Einmahl, Li, and Liu}{Einmahl
  et~al.}{2015}]{einmahl2015}
Einmahl, J. H.~J., J.~Li, and R.~Y. Liu (2015, 12).
\newblock Bridging centrality and extremity: Refining empirical data depth
  using extreme value statistics.
\newblock {\em Annals of Statistics\/}~{\em 43\/}(6), 2738--2765.

\bibitem[\protect\citeauthoryear{Esfahani and Kuhn}{Esfahani and
  Kuhn}{2018}]{esfahani_kuhn}
Esfahani, P.~M. and D.~Kuhn (2018).
\newblock Data-driven distributionally robust optimization using the
  wasserstein metric: Performance guarantees and tractable reformulations.
\newblock {\em Mathematical Programming\/}~{\em 171\/}(1-2), 115--166.

\bibitem[\protect\citeauthoryear{Fournier and Guillin}{Fournier and
  Guillin}{2015}]{fournier2015rate}
Fournier, N. and A.~Guillin (2015).
\newblock On the rate of convergence in wasserstein distance of the empirical
  measure.
\newblock {\em Probability Theory and Related Fields\/}~{\em 162\/}(3-4),
  707--738.

\bibitem[\protect\citeauthoryear{Ghaoui and Lebret}{Ghaoui and
  Lebret}{1997}]{ghaoui_regularization}
Ghaoui, L.~E. and H.~Lebret (1997).
\newblock Robust solutions to least-squares problems with uncertain data.
\newblock {\em SIAM Journal on Matrix Analysis and Applications\/}~{\em
  18\/}(4), 1035--1064.

\bibitem[\protect\citeauthoryear{Ghosh and Chaudhuri}{Ghosh and
  Chaudhuri}{2005}]{GhoshC05}
Ghosh, A.~K. and P.~Chaudhuri (2005).
\newblock On maximum depth and related classifiers.
\newblock {\em Scandinavian Journal of Statistics\/}~{\em 32\/}(2), 327--350.

\bibitem[\protect\citeauthoryear{Hallin, del Barrio, Cuesta-Albertos, and
  Matr{\'a}n}{Hallin et~al.}{2021}]{HallinDBCAM21}
Hallin, M., E.~del Barrio, J.~Cuesta-Albertos, and C.~Matr{\'a}n (2021).
\newblock {Distribution and quantile functions, ranks and signs in dimension d:
  A measure transportation approach}.
\newblock {\em The Annals of Statistics\/}~{\em 49\/}(2), 1139--1165.

\bibitem[\protect\citeauthoryear{Hlubinka, Kot\'{\i}k, and
  Venc\'{a}lek}{Hlubinka et~al.}{2010}]{hlubinka}
Hlubinka, D., L.~Kot\'{\i}k, and O.~Venc\'{a}lek (2010).
\newblock Weighted halfspace depth.
\newblock {\em Kybernetika\/}~{\em 46\/}(1), 125--148.

\bibitem[\protect\citeauthoryear{Johnson and Preparata}{Johnson and
  Preparata}{1978}]{JohnsonP78}
Johnson, D.~S. and F.~P. Preparata (1978).
\newblock The densest hemisphere problem.
\newblock {\em Theoretical Computer Science\/}~{\em 6}, 93--107.

\bibitem[\protect\citeauthoryear{Jun~Li and Liu}{Jun~Li and
  Liu}{2012}]{LiCAL12}
Jun~Li, J. A. C.-A. and R.~Y. Liu (2012).
\newblock Dd-classifier: Nonparametric classification procedure based on
  dd-plot.
\newblock {\em Journal of the American Statistical Association\/}~{\em
  107\/}(498), 737--753.

\bibitem[\protect\citeauthoryear{Kot\'{i}k and Hlubinka}{Kot\'{i}k and
  Hlubinka}{2017}]{KotikH17}
Kot\'{i}k, L. and D.~Hlubinka (2017).
\newblock A weighted localization of halfspace depth and its properties.
\newblock {\em Journal of Multivariate Analysis\/}~{\em 157}, 53--69.

\bibitem[\protect\citeauthoryear{Kuhn, Esfahani, Nguyen, and
  Shafieezadeh-Abadeh}{Kuhn et~al.}{2019}]{KuhnENSA19}
Kuhn, D., P.~M. Esfahani, V.~A. Nguyen, and S.~Shafieezadeh-Abadeh (2019).
\newblock Wasserstein distributionally robust optimization: Theory and
  applications in machine learning.
\newblock arXiv:1908.08729.

\bibitem[\protect\citeauthoryear{Lange, Mosler, and Mozharovskyi}{Lange
  et~al.}{2014}]{LangeMM14}
Lange, T., K.~Mosler, and P.~Mozharovskyi (2014).
\newblock Fast nonparametric classification based on data depth.
\newblock {\em Statistical Papers\/}~{\em 55}, 49--69.

\bibitem[\protect\citeauthoryear{Liu and Singh}{Liu and Singh}{1993}]{LiuS93}
Liu, R.~Y. and K.~Singh (1993).
\newblock A quality index based on data depth and multivariate rank tests.
\newblock {\em Journal of the American Statistical Association\/}~{\em
  88\/}(421), 252--260.

\bibitem[\protect\citeauthoryear{Lopuhaa and Rousseeuw}{Lopuhaa and
  Rousseeuw}{1991}]{LopuhaaR91}
Lopuhaa, H.~P. and P.~J. Rousseeuw (1991).
\newblock Breakdown points of affine equivariant estimators of multivariate
  location and covariance matrices.
\newblock {\em The Annals of Statistics\/}~{\em 19\/}(1), 229--248.

\bibitem[\protect\citeauthoryear{Mosler and Mozharovskyi}{Mosler and
  Mozharovskyi}{2022}]{MoslerM22}
Mosler, K. and P.~Mozharovskyi (2022).
\newblock Choosing among notions of multivariate depth statistics.
\newblock {\em Statistical Science\/}~{\em 37}, 348--368.

\bibitem[\protect\citeauthoryear{Nagy and Dvo\v{r}\'{a}k}{Nagy and
  Dvo\v{r}\'{a}k}{2020}]{nagy_illumination}
Nagy, S. and J.~Dvo\v{r}\'{a}k (2020).
\newblock Illumination depth.
\newblock {\em Journal of Computational and Graphical Statistics\/}~{\em
  30\/}(1), 78--90.

\bibitem[\protect\citeauthoryear{Nagy, Sch\"{u}tt, and Werner}{Nagy
  et~al.}{2019}]{nagy_survey}
Nagy, S., C.~Sch\"{u}tt, and E.~M. Werner (2019).
\newblock Halfspace depth and floating body.
\newblock {\em Statistics Surveys\/}~{\em 13}, 52--118.

\bibitem[\protect\citeauthoryear{Pele and Werman}{Pele and
  Werman}{2009}]{ML_earth2}
Pele, O. and M.~Werman (2009).
\newblock Fast and robust earth mover's distances.
\newblock In {\em 2009 IEEE 12th International Conference on Computer Vision},
  pp.\  460--467. IEEE.

\bibitem[\protect\citeauthoryear{Pflug and Pohl}{Pflug and
  Pohl}{2017}]{Pflug_review}
Pflug, G.~C. and M.~Pohl (2017).
\newblock A review on ambiguity in stochastic portfolio optimization.
\newblock {\em Set-Valued and Variational Analysis\/}.

\bibitem[\protect\citeauthoryear{Pollard}{Pollard}{1984}]{pollard}
Pollard, D. (1984).
\newblock {\em Convergence of Stochastic Processes}.
\newblock Springer Series in Statistics. Springer-Verlag, New York.

\bibitem[\protect\citeauthoryear{Rachev and R\"uschendorf}{Rachev and
  R\"uschendorf}{1998}]{rachev_ruschendorf}
Rachev, S.~T. and L.~R\"uschendorf (1998).
\newblock {\em Mass transportation problems. {V}ol. {I}}.
\newblock Probability and its Applications. Springer-Verlag, New York.

\bibitem[\protect\citeauthoryear{Ramsay, Durocher, and Leblanc}{Ramsay
  et~al.}{2019}]{RamsayDL19}
Ramsay, K., S.~Durocher, and A.~Leblanc (2019).
\newblock Integrated rank-weighted depth.
\newblock {\em Journal of Multivariate Analysis\/}~{\em 173}, 51--69.

\bibitem[\protect\citeauthoryear{Resnick}{Resnick}{2008}]{resnick_extremes}
Resnick, S.~I. (2008).
\newblock {\em Extreme values, regular variation and point processes}.
\newblock Springer Series in Operations Research and Financial Engineering.
  Springer, New York.
\newblock Reprint of the 1987 original.

\bibitem[\protect\citeauthoryear{Rockafellar}{Rockafellar}{1970}]{rockafellar}
Rockafellar, R.~T. (1970).
\newblock {\em Convex Analysis}.
\newblock Princeton Mathematical Series, No. 28. Princeton University Press,
  Princeton, N.J.

\bibitem[\protect\citeauthoryear{Rousseeuw and Leroy}{Rousseeuw and
  Leroy}{1987}]{RousseeuwL87}
Rousseeuw, P.~J. and A.~M. Leroy (1987).
\newblock {\em Robust Regression and Outlier Detection}.
\newblock New York: John Wiley \& Sons.

\bibitem[\protect\citeauthoryear{Rousseeuw and Van~Driessen}{Rousseeuw and
  Van~Driessen}{1999}]{RousseeuwD99}
Rousseeuw, P.~J. and K.~Van~Driessen (1999).
\newblock A fast algorithm for the minimum covariance determinant estimator.
\newblock {\em Technometrics\/}~{\em 41}, 212--223.

\bibitem[\protect\citeauthoryear{Rubner, Tomasi, and Guibas}{Rubner
  et~al.}{2000}]{ML_earth1}
Rubner, Y., C.~Tomasi, and L.~J. Guibas (2000).
\newblock The earth mover's distance as a metric for image retrieval.
\newblock {\em International Journal of Computer Vision\/}~{\em 40\/}(2),
  99--121.

\bibitem[\protect\citeauthoryear{Scarf}{Scarf}{1958}]{scarf1958min}
Scarf, H. (1958).
\newblock A min max solution of an inventory problem.
\newblock {\em Studies in the Mathematical Theory of Inventory and
  Production\/}.

\bibitem[\protect\citeauthoryear{Tukey}{Tukey}{1975}]{tukey}
Tukey, J.~W. (1975).
\newblock Mathematics and the picturing of data.
\newblock In {\em Proceedings of the International Congress of Mathematicians,
  Vancouver, 1975}, Volume~2, pp.\  523--531.

\bibitem[\protect\citeauthoryear{Wozabal}{Wozabal}{2012}]{wozabal2012framework}
Wozabal, D. (2012).
\newblock A framework for optimization under ambiguity.
\newblock {\em Annals of Operations Research\/}~{\em 193\/}(1), 21--47.

\bibitem[\protect\citeauthoryear{Wozabal}{Wozabal}{2014}]{wozabal2014robustifying}
Wozabal, D. (2014).
\newblock Robustifying convex risk measures for linear portfolios: {A}
  nonparametric approach.
\newblock {\em Operations Research\/}~{\em 62\/}(6), 1302--1315.

\bibitem[\protect\citeauthoryear{Xu, Caramanis, and Mannor}{Xu
  et~al.}{2009}]{SVM}
Xu, H., C.~Caramanis, and S.~Mannor (2009).
\newblock Robustness and regularization of support vector machines.
\newblock {\em Journal of Machine Learning Research\/}~{\em 10}, 1485--1510.

\bibitem[\protect\citeauthoryear{Zuo and Serfling}{Zuo and
  Serfling}{2000}]{zuo_serfling}
Zuo, Y. and R.~Serfling (2000).
\newblock General notions of statistical depth function.
\newblock {\em The Annals of Statistics\/}~{\em 28\/}(2), 461--482.

\end{thebibliography}

\clearpage

\setcounter{section}{0}

\begin{center}
	{\Large Supplementary Materials to \\ \indent\\ Distributionally robust halfspace depth} \\
	
	\indent\\
	
	by\\
	
	\indent\\
	
	{\large Jevgenijs Ivanovs and Pavlo Mozharovskyi}
\end{center}

\section{Proofs of properties of the inner problem}

\begin{proof}[Proof of Proposition~\ref{prop:sup}]
The result in~\eqref{eq:sol} follows readily from~\cite[Thm.\ 1 and Thm.\ 2(a)]{blanchet_distributional} by noting that the distance from $\x\in\R^d$ to $A=\{\x:\br{\u,\x-\z}\geq 0\}$ is given by $y^+$ with $y=\br{\u,\z-\x}$.
This also shows the equivalence of the two optimization problems in~\eqref{eq:opt}, which can also be proven directly by relating the ambiguity sets.
The representation in~\eqref{eq:sol_alt} for $\lambda<\infty$ follows by rewriting the expectation of the positive part.
Suppose that the value is 1 then it can not be that $\e Y^+>\delta$, because  then $\lambda<\infty$ and $\p(Y\leq \lambda)<1$ whereas $h(\lambda)\geq \delta$. 
Assuming $\e Y^+\leq \delta$ we consider the two cases $\lambda=\infty$ and $\lambda>\infty$ to find that the value is~1.

The minimization problem is analyzed by considering
\[1-\sup_{\p':d_W(\p',\p)\leq \delta}\p'(\br{\u,\X-\z}\leq 0),\]
and noting that the subset $A=\{\x:\br{\u,\x-\z}\leq 0\}$ is closed and the distance of $\x$ to $A$ is given by~$y^-$. 
We conclude by noting that $1-(1-a)^+=\min(a,1)$ and simplifying the result $\e\min(Y^-/\underline\lambda,1)-\delta/\underline\lambda$.
\end{proof}

\begin{proof}[Proof of Corollary~\ref{cor:sample}]
Note that $s_m< \delta n$ corresponds to $\hat\e Y^+<\delta$ and thus we get~1.
Otherwise, choose $1\leq k\leq  m$ as stated. Then $\lambda=y^{(k)}$ and according to Proposition~\ref{prop:sup} the supremum is given by
\[p+\frac{1}{n}\sum_{i=1}^{k-1}(1-y^{(i)}/y^{(k)})+\delta/y^{(k)},\]
because the indices corresponding to $y^{(i)}=y^{(k)}$ are irrelevant. The first result now follows, where the upper bound stems from the inequality $\delta n-s_{k-1}\leq s_k-s_{k-1}=y^{(k)}$. 

Similar argument for the second expression gives $\underline \lambda=-y^{(-k)}$ and then also
\[\frac{\underline m-k+1}{n}-\frac{1}{n}\underline s_{k-1}/y^{(-k)}+\delta/y^{(-k)}\]
and the result follows.
\end{proof}

\begin{proof}[Proof of Lemma~\ref{lem:mon}]
Observe from~\eqref{eq:sol} that $Y$ and $Y^+$ yield the same value. Hence we may focus on the non-negative random variables and assume that $Y_1\prec Y_2$.
The case $\delta=0$ follows immediately from the definition of stochastic ordering, and so we may assume $\delta>0$.
Our proof relies on the direct use of couplings instead of the representations in Proposition~\ref{prop:sup}, which are difficult to handle in this context.
Consider $Y_2$ and the corresponding $Y_2^*\in\R$ attaining the optimal value $v_2$ on a possibly extended probability space $(\Omega,\mathcal F,\p)$ allowing for a further independent uniform random variable, see~\eqref{eq:opt_transport}.
Stochastic ordering implies by a standard argument that there is a random variable $Y_1\leq Y_2$ on this latter probability space with the given law.
Furthermore, define $Y^*_1=Y^*_2-Y_2+Y_1$ and note that $\e|Y^*_1-Y_1|=\e|Y^*_2-Y_2|\leq \delta$.
Finally, we observe that 
\[v_2=\p(Y^*_2\leq 0)\leq \p(Y^*_2\leq Y_2-Y_1)=\p(Y^*_1\leq 0)\leq v_1,\]
because $Y^*_1$ belongs to the ambiguity ball around the law of~$Y$. 
\end{proof} 

\begin{proof}[Proof of Lemma~\ref{lem:mix}]
The result is obvious for $\delta=0$ and so we assume $\delta>0$.
Choose any $c_1,c_2\geq 0$ such that $pc_1+(1-p)c_2\leq 1$.
Consider $(Y_1^*,Y_1)$ and $(Y_2^*,Y_2)$ on the same probability space, where 
the first pair represents the optimal transport plan for $Y_1$ with ambiguity radius $\delta c_1$, and the second for $Y_2$ with radius $\delta c_2$, see~\eqref{eq:opt_transport}.
If $c_i=0$ we take $Y_i^*=Y_i$. 
Define $(Y^*,Y)$ as their mixture with probabilities $p$ and $1-p$, respectively.
Note that 
\[\e|Y^*-Y|=p\e|Y_1^*-Y_1|+(1-p)\e|Y_2^*-Y_2|\leq p\delta c_1+(1-p)\delta c_2\leq \delta,\]
and also 
\[v\geq \p(Y^*\leq 0)=p\p(Y_1^*\leq 0)+(1-p)\p(Y_2^*\leq 0)=pv_1(\delta c_1)+(1-p)v_2(\delta c_2).\]
Thus we have established the lower bound on~$v$ in terms of the supremum over the legal $c_i$.

Now consider a probability space and random variables $Y$ and $I$, where  $I$ is 1 with probability $p$ and 2 with probability $(1-p)$, and $Y$ given $I=i$ has the law of $Y_i$.
We may extend this probability space so that there is a random variable $Y^*$ with $(Y^*,Y)$ being the optimal coupling. 
Observe that 
\[v=\p(Y^*\leq 0)=p\p(Y^*\leq 0|I=1)+(1-p)\p(Y^*\leq 0|I=2).\]
Letting $c_i=\e(|Y^*-Y|\,|I=i)/\delta$ we find that necessarily $pc_1+(1-p)c_2\leq 1$.
Note that $\p(Y^*\leq 0|I=i)\leq v_i(\delta c_i)$, because the law of $Y^*| I=i$ is in the ambiguity ball centered at $Y_i$ with radius $c_i\delta$.
Hence $v\leq pv_1(\delta c_1)+(1-p)v_2(\delta c_2)$ for the given above $c_1,c_2$.
\end{proof}

\begin{proof}[Proof of Lemma~\ref{lem:infinite}]
Assume (i). Then for any $\epsilon\in(0,1/2)$ small we have $\p(Y_t<c)\leq \epsilon$ for large enough~$t$.
Regard $Y_t$ as the mixture corresponding to $Y_t<c$ and $Y_t\geq c$. The sandwich bounds in Lemma~\ref{lem:mix} imply that 
\[v_t\leq \epsilon+v_{2t}(2\delta)\]
for large enough~$t$, where $v_{2t}$ corresponds to the second mixing component. According to Lemma~\ref{lem:mon} such must be smaller than the value corresponding to the deterministic variable at~$c$. But the latter deminishes to 0 as $c\to\infty$ according to~\eqref{eq:det}. For (ii) we take the mixture corresponding to $Y_t\leq 0$ and $Y_t>0$, and use the lower bound.

With regard to the mixture, the statement follows readily from the sandwich bounds.
%
%Lemma~\ref{lem:mix} gives the sandwich bounds:
%\[pv_1(\delta/p)\leq v\leq pv_1(\delta/p)+(1-p)v_2(\delta/(1-p)),\]
%because $c_1\leq 1/p,c_2\leq \delta/(1-p)$ and $v_i$ are increasing in the radius.
%It is left to show that $v_2\to 0$ if $t\to\infty$ and $v_2\to 1$ if $t\to -\infty$.
%The latter follows immediately from its lower bound $\p(Y_2'+t\leq 0)\to 1$.
%
%Consider an arbitrary $T>0$ and view $Y_2$ as a mixture of $Y_2|Y_2>T$ and $Y_2|Y_2\leq T$.
%The probability $\p(Y_2\leq T)=\p(Y_2'\leq T-t)\to 0$ as $t\to \infty$. But the respective optimal value is bounded by~1 and according to the upper bound for the mixture we may assume that $Y_2>T$.
%In this case $v_2\leq \delta/T$, which follows from the representation~\eqref{eq:sol}.
%\ji{This is used in breakdown!}
%Hence indeed $v_2\to 0$ as $t\to\infty$.
\end{proof}

\begin{proof}[Proof of Lemma~\ref{lem:strict}]
Assume $v_1=v_2=v<1$, and let $\lambda_i=h^{-1}_i(\delta)$ as specified in Proposition~\ref{prop:sup}.
Observe that $\lambda_1<\infty$ and
\[\p(Y_1<\lambda_1)\leq v\leq\p(Y_1+\eta\leq\lambda_2).\]
Thus $\lambda_2\geq \lambda_1+\eta$, since otherwise $Y$ has no mass at some small interval $(\lambda_1-\epsilon,\lambda_1)$ and we readily get a contradiction to the definition of~$\lambda_1$. 
Now we write
\begin{equation}\label{eq:contr}h_2(\lambda_2-)=\e(Y_1\ind{Y_1\in(0,\lambda_2-\eta)})+\eta\p(Y_1\in(0,\lambda_2-\eta))+\e(Y_2\ind{Y_2\in(0,\eta]})\geq h_1(\lambda_1-).\end{equation}
For $\lambda_2>\lambda_1+\eta$ we get $h_2(\lambda_2-)>h_1(\lambda_1)\geq \delta$. The latter condradicts the choice of $\lambda_2$ and so we assume $\lambda_2=\lambda_1+\eta$.
According to~\eqref{eq:sol_alt} we thus must have
\[\delta-h_2(\lambda_2-)= \frac{\lambda_2}{\lambda_1}(\delta-h_1(\lambda_1-))\geq 0.\]
Hence either  $h_1(\lambda_1-)=h_2(\lambda_2-)=\delta$ or  $h_2(\lambda_2-)<h_1(\lambda_1-)$, both contradicting~\eqref{eq:contr}.
The former implies $\p(Y_1\in(0,\lambda_1))=0$ and so $\delta=h_1(\lambda_1-)=0$.
\end{proof}

\begin{proof}[Proof of Lemma~\ref{lem:cont}]
We use the representation
\[v=\inf_{\lambda'>0}\{\delta/\lambda'+\e(1-Y^+/\lambda')^+\}.\]
Note that the expression under the $\inf$ is jointly continuous in~$Y,\delta>0,\lambda'>0$, where convergence of $Y$ is understood in the weak sense.
This readily follows by the dominated convergence theorem. %by the Skorokhod's representation theorem and 
Thus we get the required continuity for the modified problem where $\lambda'>0$ is replaced by some compact subinterval $\lambda'\in[a,b]$.

Assume that $\e Y^+>\delta$ so that the optimal $\lambda=h^{-1}(\delta)\in(0,\infty)$.
The optimal $\lambda_n$ corresponding to $Y_n$ and $\delta_n$ may or may not converge to~$\lambda$. % (the latter happens iff $h(\lambda)=\delta<h(\lambda-)$).
Importantly, we have 
\[h^{-1}(\delta)\leq \liminf\lambda_n\leq \limsup\lambda_n\leq h^{-1}(\delta+),\] see~\cite[Prop.\ 0.1]{resnick_extremes}.
Hence for large enough $n$ we may restrict~$\lambda'$ to a compact interval while preserving the infimum. 

Now assume that $\e Y^+\leq  \delta$ and so $v=1$.
We need to prove that $v_n(\delta_n)\to 1$.
For any $\epsilon>0$ we may choose $\delta'\in(0,\delta)$ such that $1-\epsilon<v(\delta')<1$ and hence $\e Y^+>\delta'$, which follows from the left-continuity of~$h^{-1}$ and the right expression in~\eqref{eq:sol}.
But now we can apply the above proven fact to deduce that $v_n(\delta')\to v(\delta')$.
Eventually $\delta_n>\delta'$ and so $v_n(\delta_n)\geq v_n(\delta')$. Hence $\liminf v_n(\delta_n)\geq 1-\epsilon$ and the proof is complete.
\end{proof}

\section{Proofs of properties of the robust halfspace depth}

\begin{proof}[Proof of Proposition~\ref{prop:stdprops}]
P1. 
Since every vector can be represented as $A\x+\bs b$ we have
\[D_\delta(A\z+\bs b|\p_{A\X+\bs b})=\inf_{\norm{u}=1}\sup_{\p':d_W(\p'_{A\X+\bs b},\p_{A\X+\bs b})\leq \delta}\p'(\br{\u,(A\X+\bs b)-(A\z+\bs b)}\geq 0).\]
From the definition of the Wasserstein distance and $\norm{(A\x-\bs b)-(A\y-\bs b)}=\norm{\x-\y}$ %and the fact that $A\X+\bs b\stackrel{d}{=}A\X'+\bs b$ iff $\X\stackrel{d}{=}\X'$ 
we readily deduce that 
\[d_W(\p'_{A\X+\bs b},\p_{A\X+\bs b})=d_W(\p'_{\X},\p_{\X}).\] 
It is left to note that $\br{\u,(A\X+\bs b)-(A\z+\bs b)}=\br{A^\top\u,\X-\z}$, where the vector $A^\top \u$ runs over all unit vectors. The result now follows.

P2. For every $\z\neq 0$ we pick $\u=\z/\norm{\z}$ and observe that 
\[Y=\br{\u,\z-\X}=\norm{\z}-\br{\z,\X}/\norm{\z}\geq \norm{\z}-\norm{\X}.\] 
According to Lemma~\ref{lem:mon} the depth is upper-bounded by the value $\tilde v$ corresponding to $\tilde Y= \norm{\z}-\norm{\X}\to \infty$ a.s.
 Conclude by Lemma~\ref{lem:infinite}(i).

P3. It is sufficient to show that for any $\z_1,\z_2\in \R^d$ and any $t\in [0,1]$ we have
\[D_\delta(t\z_1+(1-t)\z_2|\p)\geq \min(D_\delta(\z_1|\p),D_\delta(\z_2|\p)).\]
Note that $\br{\u,\z_2-\z_1}\geq 0$ implies $\br{\u,\X-\z_2}\leq\br{\u,\X-(t\z_1+(1-t)\z_2)}$, whereas $\br{\u,\z_2-\z_1}\leq 0$ implies  $\br{\u,\X-\z_1}\leq\br{\u,\X-(t\z_1+(1-t)\z_2)}$.
The result now follows from~\eqref{eq:minimax}.
\end{proof}

\begin{proof}[Proof of Proposition~\ref{prop:addprops}]
P4. 
Fix $\z$ and note that $\br{\u,\z-\X}\leq \norm{\z-\X}$ when $\norm{\u}=1$. Choose $c>0$ large and $\epsilon>0$ small such that the later norm is at most $c$ with probability~$\epsilon$. 
Thus for all directions $\u$ we have $\p(Y\leq c)\geq \epsilon$. According to~\eqref{eq:lowerbound} the depth $D_\delta(\z|\p)$ is lower bounded by $\min(\epsilon,\delta/c)>0$.

%We use the first expression in~\eqref{eq:sol} and note that 
%\begin{equation}\label{eq:lowerbound}\delta/\lambda'+\e(1-Y^+/\lambda')^+\geq \max(\delta/\lambda',1-(\e Y^+-\delta)/\lambda').\end{equation}
%Hence it is sufficient to show that $\e Y^+$ stays bounded over all directions~$\u$.
%But $Y\leq \norm{\X-\z}$ and we are done.

P5. Since the space of angles is compact, it is enough to show that the inner supremum in~\eqref{eq:minimax} is jointly continuous in $\u,\z,\delta>0$.
Take $(\u_n,\z_n)\to(\u,\z)$ so that $Y_n=\br{\u_n,\z_n-\X}\to \br{\u,\z-\X}=Y$ a.s.\ and apply Lemma~\ref{lem:cont}.

P6. 
Assume that the upper level set is non-empty for the given~$\a$. By P2 it must be bounded.
By continuity we see that the boundary of the upper level set corresponding to $\a\in(0,1)$ must yield depth~$\a$.
Suppose there is a point $\z_1$ strictly inside this convex set with $D_\delta(\z_1|\p)=\a$. Choose the corresponding optimal direction~$\u$ which must exist by continuity of the inner supremum in~$\u$.
Find $\z_2$ on the boundary such that $\z_2=\z_1+\u \eta$ for some $\eta>0$. Now we have $\br{\u,\z_2-\z_1}=\eta>0$, and so $Y_2-Y_1=\eta$ with $Y_i=\br{\u,\z_i-\X}$.
Furthermore $\a=v_1\leq v_2$, because the direction $\u$ is not necessarily optimal for $\z_2$, but it is for~$\z_1$.
Lemma~\ref{lem:mon} shows that $v_1=v_2$ and then Lemma~\ref{lem:strict} implies that $\a=1$, a contradiction.
\end{proof}

\begin{proof}[Proof of Proposition~\ref{prop:approx}]
It is sufficient to show for any direction $\u$ that 
\[\sup_{d_W(\p',\p)\leq \delta}\p'(\br{\u,\X-\z}\geq 0)-\p(\br{\u,\X-\z}\geq 0)\downarrow 0\]
as $\delta\downarrow 0$.
We may assume that $\p(Y>0)>0$ since otherwise both terms are equal to~1. For $\delta>0$ small enough we have $\e Y^+>\delta$ and so $\lambda_\delta=h^{-1}(\delta)<\infty$.
Now the solution in~\eqref{eq:sol_alt} gives the representation of the above difference:
\begin{equation}\label{eq:approx}\p(Y\in(0,\lambda_\delta])-(h(\lambda_\delta)-\delta)/\lambda_\delta,\qquad (h(\lambda_\delta)-\delta)/\lambda_\delta\in[0,\p(Y=\lambda_\delta)).\end{equation}
This indeed converges to~0 if $\lambda_\delta\to 0$.
Assume that $\lambda_\delta\to\lambda>0$ and note that necessarily $\p(Y\in(0,\lambda))=0$.
If $Y$ has no mass at $\lambda$ then we are done, and if such mass is positive then it must be that $(h(\lambda_\delta)-\delta)/\lambda_\delta$ converges to this mass.
Hence the expression in~\eqref{eq:approx} goes to 0 and the proof is  complete.
\end{proof}

\begin{proof}[Proof of Lemma~\ref{lem:oa}]
Our proof relies on the joint continuity of the depth function, see P5.
Consider $\delta_n\uparrow \delta$ and let $\z$ be a point achieving the maximal depth for~$\delta$.
Now $D_{\delta_n}(\z|\p)\to D_\delta(\z|\p)$ and the latter can not be exceeded for ambiguity radius~$\delta_n$, showing left-continuity of~$\oa$. 
Let $\delta_n\downarrow \delta$ with $\z_n$ being points with maximal depth. These $z_n$ must belong to a compact set, see P2 and use monotonicity of depth in~$\delta_n$. 
Thus $(z_n)$ contains a convergent subsequence establishing right-continuity of $\oa$.
\end{proof}

\begin{proof}[Proof of Lemma~\ref{lem:decay_bound}]
Take $\u=\z/\norm{\z}$ assuming $\z\neq \bs 0$ and reconsider the optimal transport plan $\X^*$ satisfying~\eqref{eq:opt_transport}.
Note that $\br{\u,\X^*-\z}\leq \norm{\X^*}-\norm{\z}$, and so we have the bound
\[D_\delta(\z|\p)\leq \p(\norm{\X^*}-\norm{\z}\geq 0)\leq \frac{\e\norm{\X^*}^p}{\norm{\z}^p},\]
where $\e\norm{\X^*-\X}\leq \delta$; here we have also used Markov's inequality.
Using $(a+b)^p\leq a^p+b^p$ for $p\leq 1$ we find that 
\[\e \norm{\X^*}^p\leq \e(\norm{\X^*-\X}+\norm{\X})^p\leq \e\norm{\X^*-\X}^p+\e\norm{\X}^p.\]
It is left to observe that $\e\norm{\X^*-\X}^p\leq \e(1+\norm{\X^*-\X})\leq (1+\delta)$.
\end{proof}

\begin{proof}[Proof of Lemma~\ref{lem:decay_sublin}]
Choose a constant $r>0$ large enough so that $\p(\norm{\X}\leq r)\geq 1/2$, and let $B_r$ be the ball of radius $r$ centered at the origin.
For any point $\z$ and any direction $\u$ we must have 
\[\br{\u,\z-\x}\leq \norm{\z}+r,\qquad \forall \x\in B_r.\]
Thus $Y=\br{\u,\z-\X}\leq \norm{\z}+r$ with probability $p\geq 1/2$, and according to~\eqref{eq:lowerbound}
we find
\[D_\delta(\z|\p)\geq \min(1/2,\delta/(\norm{\z}+r)).\]
%where $\hat v$ corresponds to $\hat Y=\norm{\z}+r$, which does not depend on the direction. But $\hat v(\delta)=\frac{\delta}{\norm{\z}+r}\wedge 1$ which follows easily from~\eqref{eq:sol_alt}.
This lower bound behaving asymptotically as $\delta/\norm{\z}$ and the result follows.
\end{proof}

\begin{proof}[Proof of Corollary~\ref{cor:lindec}]
Combine Lemma~\ref{lem:decay_bound} and Lemma~\ref{lem:decay_sublin}.
\end{proof}

\begin{proof}[Proof of Lemma~\ref{lem:spheredist}]
Consider $Y=\br{\u,\z-\X}\stackrel{d}{=}\br{\u,\z}-Z$, where $Z$ has the distribution of $\br{\u,\X}$; the same for all $\u$ by assumption.
The largest $\br{\u,\z}$ equals $\norm{\z}$, and according to Lemma~\ref{lem:mon} it gives the smallest value, and hence the depth. According to
\eqref{eq:sol_alt} we then have
\[D_\delta(\z|\p)\geq \p(\norm{\z}-Z< \lambda)=\p(Z>\norm{\z}-\lambda)=(\norm{\z}-\lambda)^{-\a}\ell(\norm{\z}-\lambda),\]
where $\lambda=\lambda(\norm{\z})$ is defined using the random variable~$\norm{\z}-Z$.
For the lower bound it is thus left to show that $\lambda/\norm{\z}\to 0$, see the uniform convergence theorem~\cite[Thm.\ 1.2.1]{BGT}.
For the upper bound we use $\p(\norm{\z}-Z\leq \lambda)$ resulting in the same asymptotic behavior.

On the contrary, assume that for some $\epsilon>0$ we have $\lambda>4\epsilon\norm{\z}$ for some $\z$ with arbitrarily large norm. For simplicity we write $\eta=\norm{\z}$.
From the definition of $\lambda$ we then must have
\[\delta\geq \e\big((\eta-Z)\ind{\eta-Z\in (0,2\epsilon \eta]}\big)\geq \epsilon\eta\p(Z\in [(1-2\epsilon) \eta,(1-\epsilon)\eta)]).\]
The latter probability is asymptotic to $((1-2\epsilon)^{-\a}-(1-\epsilon)^{-\a})\eta^{-\a}\ell(\eta)$, and so for $\a\in(0,1)$ the right hand side in the display increases to $\infty$ as $\eta\to\infty$, which is a contradiction.
\end{proof}

\begin{proof}[Proof of Lemma~\ref{lem:maxmin}]
It is only required to show that the supremum can not exceed $\delta/\sqrt 2$, since by moving $\delta/\sqrt 2$ mass into the origin we readily achieve the value~$\delta/\sqrt 2$.
Note that the optimal direction $\u$ in this case is not unique.

Fix $\p'$ in the ambiguity ball centered at $\p$, choose $n\geq 1$ and consider the angle $\theta=\frac{\pi}{4(n+1)}$.
Next, we consider the masses in the cones 
\begin{align*}
m_i&=\p'\Big(X_2>0,X_2/X_1\in[\tan(\pi/4+i\theta),\tan(\pi/4+(i+1)\theta))\Big), \qquad i=1,\ldots, 2n,\\
\hat m_i&=\p'\Big(X_2<0,-X_2/X_1\in[\tan(\pi/4+i\theta),\tan(\pi/4+(i+1)\theta))\Big), \qquad i=1,\ldots, 2n,\\
m_0&=\p'\Big(X_1=X_2=0\text{ or }X_1<0,X_2/X_1\in[\tan(\pi/2+n\theta),\tan(-\pi/2-n\theta)]\Big),
\end{align*}
see Figure~\ref{fig:maxmin}.
%\begin{align*}
%&a_1=\p'(0<X_1\leq X_2/2), &a_2=\p'(0<X_1\leq -X_2/2),\\
%&b_1=\p'(-X_2/2<X_1\leq 0), &b_2=\p'(X_2/2<X_1\leq 0),\\
%&c=\p'(X_1\leq -|X_2/2|),
%\end{align*}
%see Figure~\ref{fig:maxmin}.
\begin{figure}
\centering
\begin{tikzpicture}
\filldraw (1,1) circle(1pt) node[right]{$\frac{1}{2}$};
\filldraw (1,-1) circle(1pt) node[right]{$\frac{1}{2}$};
\filldraw[red!30] (0,0) -- (1,2.4) -- (0,2.4);
\filldraw[blue!30] (0,0) -- (1,-2.4) -- (0,-2.4);
\filldraw[red!60] (0,0) -- (0,2.4) -- (-1,2.4);
\filldraw[blue!60] (0,0) -- (0,-2.4) -- (-1,-2.4);
\filldraw[magenta] (0,0) -- (-1,2.4) -- (-1,-2.4);
\draw[->] (0,0) -- (2.5,0);
\draw[->] (0,0) -- (0,2.5);
\draw (0.4,2) node{$m_1$};
\draw (0.4,-2) node{$\hat m_1$};
\draw (-0.4,2) node{$m_2$};
\draw (-0.4,-2) node{$\hat m_2$};
\draw (-0.5,0) node{$m_0$};
\draw (1,1)--(0.5,1.207);
\draw (1,1)--(0,1);
\draw (1,1) --(-0.207,0.5);
\end{tikzpicture}
\caption{Illustration of the cones in the proof of Lemma~\ref{lem:maxmin} for $n=1$.}
\label{fig:maxmin}
\end{figure}
Observe that 
\begin{align}
& \inf_{\u\in\R^d}\p'(\br{\u,\X}\geq 0)\leq \label{eq:maxmin_min} \\ 
& \min\{m_1+\cdots +m_{2n}+m_0,\,m_2+\cdots +m_{2n}+m_0+\hat m_{2n},\ldots,m_0+\hat m_{2n}+\cdots+\hat m_{1}\}, \nonumber
\end{align}
where the minimum runs over $2n+1$ terms corresponding to certain halfspaces.
Moreover, the distance from the point $(1,1)$ to the $i$th halfspace (the closest point therein) is $\sin(i\theta)\sqrt 2$, and the distances are reversed for the point $(1,-1)$.
Furthermore,  $d_W(\p',\p)\leq \delta$ implies the constraint
\begin{equation}\label{eq:maxmin_constr}
(m_1+\hat m_1)\sin(\theta)+\cdots+(m_{2n}+\hat m_{2n})\sin(2n\theta)+m_0\sin((2n+1)\theta)\leq \delta/\sqrt 2.
\end{equation}

We now maximize the minimum in~\eqref{eq:maxmin_min} under the constraint~\eqref{eq:maxmin_constr}, which must yield an upper bound on the supremum of interest. 
Observe that the solution must make all the $2n+1$ terms in~\eqref{eq:maxmin_min} equal and so $m_i=\hat m_{2n+1-i}$ for $i=1,\ldots,2n$. 
This is so since otherwise the mass can be redistributed in such a way that smallest terms become slightly larger. 
The underlying argument requires a little thought which is left for the reader. 

Letting $\tilde m_i=m_i+m_{2n+1-i}=\hat m_{2n+1-i}+\hat m_i$ for $i=1,\ldots, n$ we observe that we need to maximize $\tilde m_1+\cdots+\tilde m_n+m_0$ subject to 
\[\tilde m_1(\sin(\theta)+\sin(2n\theta))+\cdots+\tilde m_n(\sin(n\theta)+\sin((n+1)\theta))+m_0\sin((2n+1)\theta)\leq \delta/\sqrt 2.\]
By concavity of $\sin$ on $[0,\pi/2]$ we see that $\sin((2n+1)\theta)$ is the smallest coefficient and so the maximum is obtained by taking $m_0=\frac{1}{\sin((2n+1)\theta)}\frac{\delta}{\sqrt 2}$ and all other variables being zero.
It is left to note that 
\[\sin((2n+1)\theta)=\sin((2n+1)\pi/(4(n+1)))\uparrow \sin(\pi/2)=1\]
as $n\to\infty$.
\end{proof}

\section{Proofs for finite-sample version and affine invariance}

\begin{proof}[Proof of Proposition~\ref{prop:consistency}]
Firstly, the case $\delta=0$ is classical~\citep{donoho_gasko}.
According to $|\inf f-\inf g|\leq \sup|f-g|$ and the representation~\eqref{eq:sol} of the depth for $\delta>0$, it is sufficient to show that 
\[\sup_{\z,\u\in\R^d}|\e_n(1-\br{\u,\z-\X}^+)^+-\e(1-\br{\u,\z-\X}^+)^+|\to 0,\]
where the division by $\lambda'$ is incorporated into~$\u$. 

Now the result follows from empirical process theory~\cite[Thm.\ II.24 and Lem.\ II.25]{pollard}. It is only needed to observe that the graph of the function $f_{\u,\z}(\x)=(1-\br{\u,\z-\x}^+)^+\in[0,1]$ can be constructed from three half-spaces in $\R^{d+1}$ using union and intersection operations, and so the respective class has polynomial discrimination~\cite[Lem.\ II.15]{pollard}.
\end{proof}

\begin{proof}[Proof of Corollary~\ref{cor:convdelta}]
Consider the bound
\[|D_{\delta_n}(\z|\p_n)- D_0(\z|\p)|\leq |D_{\delta_n}(\z|\p_n)- D_{\delta_n}(\z|\p)|+|D_{\delta_n}(\z|\p)- D_0(\z|\p)|\]
and observe that both terms converge to 0 according to Proposition~\ref{prop:consistency} and Proposition~\ref{prop:approx}.
\end{proof}

\begin{proof}[Proof of Lemma~\ref{lem:astar}]
According to~\eqref{eq:inc_v} the function $p\oa(\delta/p)$ is non-decreasing in $p\in[0,1]$, which follows by considering all points $\z$ and all directions~$\bs u$; it is assumed to be 0 at~0.
Hence $p\oa(\delta/p)-(1-p)$ is strictly increasing with values $-1$ and $\oa(\delta)>0$ at the end points.
This function inherits continuity from $\oa$, see Lemma~\ref{lem:oa}.
Hence there is a unique root, proving the first claim. The equivalence of inequalities follows from the strict monotonicity.
Finally, $\a^*=1/[1+1/\oa(\delta/(1-\a^*))]$ which must be in $(0,1/2]$ since $\oa\leq 1$, and the value $1/2$ corresponds to $\oa(\delta/(1-1/2)=1$.
The fact that $\a^*(\delta)$ is non-decreasing is inherited from the same property of $\oa(\delta)$.
Note also that $(1-\oa)\oa(\delta/(1-\oa))\leq \oa$ and hence $\a^*\leq \oa$.
\end{proof}

\begin{proof}[Proof of Proposition~\ref{prop:bpastar}]
The proof follows the same ideas as in~\cite[Lem 3.1]{donoho_gasko}, but here we rely on the properties of the optimal value~$v$ (for a fixed direction) established in Section~\ref{sec:sol}. 
Let $m$ be the smallest integer satisfying $m/(n+m)\geq \a$, that is, $m=\lceil n\a/(1-\a)\rceil$.
First, we show that $m$ new points are sufficient for the breakdown.
We place all $\y_i$ at the same location $\z=t\u$ and let $t\to\infty$. The value at $\z$ along any direction is at least~$m/(m+n)$ according to Lemma~\ref{lem:mix}.
Hence the depth at this far removed $\z$ is at least $\a$ and we are done.

Next, suppose $k<m$ points are sufficient to achieve infinite distance while having non-empty upper level set.
Let $\z$ be a point in this new upper level set. We denote its distance to the convex hull of $\x_i$ by $d$. By assumption we may choose $\z$ such that $d$ is arbitrarily large.
Take the direction $\u$ along which this distance is computed, and such that it points in the direction of~$\z$.
The value at $\z$ along the direction $\u$ converges to $k/(n+k)<\a$ as $d\to\infty$ according to Lemma~\ref{lem:infinite}.
This yields a contradiction.

Now we show that for any $k<m$ additional points the upper level set is non-empty; the original level set is non-empty since $\a^*\leq \oa$ according to Lemma~\ref{lem:astar}. %; hence also the original level set is non-empty. 
Let $p=n/(n+k)$ and take a point $\z$ with the maximal depth $\oa(\delta/p)$ in the original data set for an inflated ambiguity ball.
For any direction $\u$ and arbitrary $\y_i$ the new depth of $\z$ is larger than $p\oa(\delta/p)$, see Lemma~\ref{lem:mix}.
It is left to check that this lower bound is not smaller than~$\a$. By assumption $k<n\a/(1-\a)$ and thus $p>1-\a$.
According to~\eqref{eq:inc_v} we indeed have
\[p\oa(\delta/p)\geq (1-\a)\oa(\delta/(1-\a))\geq \a,\]
where the latter follows from Lemma~\ref{lem:astar} and the assumption $\a\leq\a^*(\delta,\p_n)$.

The final statement, in view of Lemma~\ref{lem:astar}, only requires checking that $\oa(\delta,\p_n)\to \oa(\delta,\p)$ a.s.\ for any $\delta\geq 0$, which is a consequence of the uniform consistency result in Proposition~\ref{prop:consistency}, see also the proof of Lemma~\ref{lem:oa}.
\end{proof}

\begin{proof}[Proof of Proposition~\ref{prop:breakdown}]
Suppose there exists $\z$ outside of the convex hull of $\x_i$, and such that it belongs to the new median region.
Its depth is no more than $p+\delta/d$ with $p=m/(n+m)$, where $d$ is the distance from the convex hull, see Lemma~\ref{lem:mix} as well as Lemma~\ref{lem:mon} with~\eqref{eq:det}; alternatively one may use Lemma~\ref{lem:infinite}.
The points in the old median region have new depth at least $(1-p)\oa(\delta/(1-p),\p_n)$, and hence we must have 
\[p\geq (1-p)\oa(\delta/(1-p),\p_n).\]
According to Lemma~\ref{lem:astar} we must have $p\geq \a^*_n$ and then $m\geq n\a^*_n/(1-\a^*_n)$.

%The equivalence of $p^*=1/2$ and $q^*(2\delta)=1$ follows readily from the definition of~$p^*$.
Let us now show that $m=n$ new points are sufficient for the breakdown, so that the breakdown point is exactly $1/2$ when $\a_n^*=1/2$.
Suppose breakdown does not occur, and so the new median must be contained in some inflation of the given convex hull. Now we copy the given $n$ points and shift them sufficiently far away, so that the two inflations do not intersect.
But we may regard the new set of points as the original data set, which readily leads to a contradiction.
The final statement follows from the convergence $\a^*_n\to \a^*(\delta,\p)$ a.s.
\end{proof}

\begin{proof}[Proof of Proposition~\ref{prop:affinv}]
	The proof is straightforward noting $\V=A\X + \bs b$ and noticing that, for each $\u\in\R^d$, $\|u\|=1$, the following holds
	\begin{equation*}
		\br{\Lambda_{V}\u,\V} = (\Lambda_{V}\u)^\top\V = \bigl((A^\top)^{-1}\Lambda_X\u\bigr)^\top(A\X) = \u^\top\Lambda_X^\top A^{-1}AX = \br{\Lambda_X\u,X},
	\end{equation*}
	and the same for $\z$; the second equality holds because $(\Sigma_V)^{-1}=(A^\top)^{-1}\Lambda_X\Lambda_X^\top A^{-1}$.
\end{proof}

\section{Proofs of properties of the depth-trimmed regions}

\begin{proof}[Proof of Proposition~\ref{prop:outer}]
Let $\z^*$ be the unique point in $H$ such that $\norm{\z-\z^*}=d(\z,H)\geq \delta n$. 
Consider the direction $\u=(\z-\z^*)/\norm{\z-\z^*}$ and observe that this direction maximizes $\min_i \br{\u,\z-\x_i}$.
This can be easily seen by examining the vertices of the face containing~$\z^*$.
For such $\u$ consider the solution to the inner problem as described in Corollary~\ref{cor:sample}.
Note that $m=n$ and $y^{(1)}=d(\z,H)\geq \delta n$. Hence we get the solution~$\delta/d(\z,H)\leq 1/n$. Any other direction yields a larger value, since either $y^{(1)}$ is smaller or $m<n$. This proves the first statement.

Since $\delta/\a\geq \delta n$ we see from the first statement that $\{\z\in\R^d:d(\z,H)= \delta/\a\}$ must have depth~$\a$.
%The points $\z$ on the outside of this closed curve must have a smaller depth. 
The second result now follows from the property P6 of the depth.

For the third result we let $c$ be the maximal distance of the points in $H$ from the origin. Then $\norm{\z}-c\leq d(\z,H)\leq \norm{\z}+c$ showing that $d(\z,H)\sim \norm{\z}$, and the result follows easily from the first statement.
\end{proof}

\begin{proof}[Proof of Lemma~\ref{lem:geometry}]
Note that we optimize over a compact set and that the objective function is continuous in~$\u$. Thus the $\rm{ argmax}$ is non-empty.

We employ the induction in~$n$. For $n=1$ the result is obvious with $I=\{1\}$ and $\u=\vv_1/\norm{\vv_1}$.
Suppose that our result is proven for some $n\geq 1$, and consider $n+1$ non-zero vectors~$\vv_i$.
Let $\u^*$ be any of the optimal directions. If $\br{\u^*,\vv_k}\leq 0$ for some $k\leq n+1$ then 
\[\sum_{i\leq n+1}\br{\u^*,\vv_i}^+=\sum_{i\leq n+1,i\neq k}\br{\u^*,\vv_i}^+\leq \sum_{i\leq n+1,i\neq k}\br{\u^*_k,\vv_i}^+\leq \sum_{i\leq n+1}\br{\u^*_k,\vv_i}^+,\]
where $\u^*_k$ is any of the optimal directions for the problem with $\vv_k$ excluded. 
By maximality of the left hand side we must have equalities. % and thus $\br{\u^*_k,\vv_k}\leq 0$. 
Thus $\u^*$ is an optimal direction for the problem with $\vv_k$ excluded, and by the inductive assumption $\u^*$ has the stated form with $I\subset\{1,\ldots,n\}\backslash\{k\}$.
It is left to recall that $\br{\u^*,\vv_k}\leq 0$, and so the form is as stated.

Finally, we assume that there is an optimal direction \[\u^*\in S^+=\{\u\in\R^d:\br{\u,\vv_i}>0\,\forall i\leq n+1\}.\]
One readily checks that $S^+$ is a convex set with a non-empty interior by assumption.
This leads to the convex optimization problem 
\[\sup_{\u\in S^+,\norm{\u}\leq 1}\sum_{i=1}^{n+1}\br{\u,\vv_i}\]
with a linear objective function, where the constraint $\norm{\u}\leq 1$ is equivalent to $\norm{\u}=1$.
Clearly there is a feasible solution $\u$ in the interior of $S^+$ which has $\norm{\u}<1$.
By standard theory~\cite[Thm.\ 28.3]{rockafellar} every maximizer is obtained by maximizing the Lagrangian
\[\sum_{i=1}^{n+1}\br{\u,\vv_i}-\lambda(\u^\top \u-1)\]
in $\u\in S^+$ for some $\lambda\geq 0$. In the case $\lambda>0$ the candidate solutions are of the form $\u=\sum_{i\leq n+1} \vv_i/(2\lambda)$, because the boundary of $S^+$ does not have common points with $S^+$.
This solution must be in $S^+$ and so $\br{\u,\vv_i}>0$ for all $i$, yielding the stated form with $I=\{1,\ldots,n+1\}$.
The case $\lambda=0$ results in $\sum \vv_i=0$, but then the optimal value is 0, contradicting our assumption on~$\u^*$.

Considering the optimal value, we observe for $I$ satisfying~\eqref{eq:I} that
\[\sum_{j\leq n}\br{\sum_{i\in I}\vv_i/c,\vv_j}^+=\sum_{j\in I}\br{\sum_{i\in I}\vv_i,\vv_j}/c=\norm{\sum_{i\in I}\vv_i},\]
where $c=\norm{\sum_{i\in I}\vv_i}$. Now take an arbitrary subset $J$ and let $\u=\sum_{i\in J}\vv_i/c$ with $c$ such that $\norm{\u}=1$.
Note that 
\[\sum_j\br{\u,\vv_j}^+\geq \sum_{j\in J}\br{\u,\vv_j}=\norm{\sum_{i\in J}\vv_i},\]
because the second sum ignores some positive contribution and adds up some negative contribution. But the left hand side is upper bounded by the optimal value, and thus adding additional terms $\norm{\sum_{i\in J}\vv_i}$ to the maximum does not change the result.
\end{proof}

\begin{proof}[Proof of Proposition~\ref{prop:median}]
From Corollary~\ref{cor:sample} we see that $D_\delta(\z|\p_n)=1$ holds iff 
$\sum_i\br{\u,\x_i-\z}^+\leq \delta n$ for all directions~$\u$ (we take the reverse sign here for convenience). Thus according to Lemma~\ref{lem:geometry} the depth is~1 iff $\norm{\sum_{i\in I}\vv_i}\leq \delta n$ with $\vv_i=\x_i-\z$ for all $I\neq \emptyset$ satisfying~\eqref{eq:I}, and then also for all non-empty~$I$.
This inequality can be rewritten as  
\[\norm{\z-\sum_{i\in I}\x_i/|I|}\leq \delta n/|I|.\] %for all $I\neq\emptyset$ satisfying~\eqref{eq:I}, or alternatively for all non-empty subsets of $\{1,\ldots,n\}$.
%The problem with~\eqref{eq:I} is that it depends on~$\z$, whereas the number of all subsets can be prohibitively large.
Recall that $I$ satisfying~\eqref{eq:I} for some $\z$ is necessarily such that $\{\x_i:i\in I\}$ can be separated from the rest by a hyperplane, and the stated form of $M_\delta$ follows.

By comparing the radius $\delta$ of the ball corresponding to all the observations and the radius $\delta n/m$ of the ball corresponding to some $I$ of size~$m$.
The difference of such radii is $\delta(n/m-1)\to \infty$ as $\delta\to\infty$ for any $m\in\{1,\ldots,n-1\}$. Hence the intersection in the representation of $M_\delta$ will eventually become the first ball.
\end{proof}

\end{document}